\documentclass[onefignum,onetabnum]{siamonline220329}

\usepackage{subcaption}

\usepackage{lipsum}
\usepackage{amsfonts}
\usepackage{graphicx}
\usepackage{epstopdf}
\usepackage{todonotes}
\usepackage{mathtools}
\usepackage{amssymb}
\usepackage{cleveref}
\usepackage{booktabs}

\usepackage{subcaption}

\usepackage{algpseudocode}
\usepackage{algorithm}

\ifpdf
  \DeclareGraphicsExtensions{.eps,.pdf,.png,.jpg}
\else
  \DeclareGraphicsExtensions{.eps}
\fi

\newcommand{\N}{{\mathbb N}}

\newcommand{\R}{{\mathbb R}}

\newcommand{\Z}{{\mathbb Z}}
\newcommand{\E}{{\mathbb E}}

\newcommand{\norm}[1]{\left\lVert#1\right\rVert}

\DeclarePairedDelimiterX{\inp}[2]{\langle}{\rangle}{#1, #2}

\DeclareMathOperator*{\argmin}{arg\,min}
\DeclareMathOperator*{\arginf}{arg\,inf}

\usepackage{color,soul}
    \usepackage{mathtools}
\definecolor{darkred}{rgb}{.7,0,0}
\definecolor{darkgreen}{rgb}{.3,0.7,0}
\definecolor{purple}{rgb}{0.5, 0, 0.5} 
\definecolor{lightblue}{rgb}{0.68, 0.85, 0.9}
\definecolor{lightred}{rgb}{1.0, 0.6, 0.6}

\crefname{section}{Section}{Sections}
\crefname{subsection}{Subsection}{Subsections}

\newtheorem{assumption}[theorem]{Assumption}
\newtheorem*{assumptionn}{Assumption}
\renewcommand{\theassumptionn}{P}

\usepackage{enumitem}
\setlist[enumerate]{leftmargin=.5in}
\setlist[itemize]{leftmargin=.5in}
\newlist{assumpenum}{enumerate}{5}
\setlist[assumpenum]{leftmargin=4.3em, font={\bfseries}, label=\arabic*., itemsep=0.5em}
\crefname{assumpenumi}{Assumption}{Assumptions}

\newlist{contrib}{enumerate}{5}
\setlist[contrib]{leftmargin=4em, font={\bfseries}, label=\arabic*., itemsep=0.5em}
\crefname{contribi}{Contribution}{Contributions}

\newsiamremark{remark}{Remark}
\newsiamremark{hypothesis}{Hypothesis}
\crefname{hypothesis}{Hypothesis}{Hypotheses}
\newsiamthm{claim}{Claim}

\title{
Solving Roughly Forced Nonlinear PDEs\\ via Misspecified Kernel Methods and Neural Networks}

\author{
    \hspace{0.5em}Ricardo Baptista\footnote{Computing and Mathematical Sciences, California Institute of Technology, Pasadena, CA 91125.}\and
    Edoardo Calvello\footnotemark[1]\and
    Matthieu Darcy\footnotemark[1] \thanks{Corresponding author: \url{mdarcy@caltech.edu}.}\and 
    \vspace{0.001em} \\
    Houman Owhadi\footnotemark[1]\and
    Andrew M. Stuart\footnotemark[1]\and
    Xianjin Yang\footnotemark[1]
}

\usepackage{amsopn}

\ifpdf
\hypersetup{
  pdftitle={Solving Roughly Forced Nonlinear PDEs
via Misspecified Kernel Methods and Neural Networks},
  pdfauthor={Matthieu Darcy}
}
\fi

\begin{document}

\maketitle
\nocite{*}

\begin{abstract}
We consider the use of Gaussian Processes (GPs) or Neural Networks (NNs) to numerically approximate the solutions to nonlinear partial differential equations (PDEs) with rough forcing or source terms, which commonly arise as pathwise solutions to stochastic PDEs. 
Kernel methods have recently been generalized to solve nonlinear PDEs by approximating their solutions as the maximum a posteriori estimator of GPs that are conditioned to satisfy the PDE at a finite set of collocation points. The convergence and error guarantees of these methods, however, rely on the PDE being defined in a classical sense and its solution possessing sufficient regularity to belong to the associated reproducing kernel Hilbert space.
We propose a generalization of these methods to handle roughly forced nonlinear PDEs while preserving convergence guarantees with an oversmoothing GP kernel that is misspecified relative to the true solution's regularity.
This is achieved by conditioning a regular GP to satisfy the PDE with a modified source term in a weak sense (when integrated against a finite number of test functions). This is equivalent to replacing the empirical $L^2$-loss on the PDE constraint by an empirical negative-Sobolev norm. 
We further show that this loss function can be used to extend physics-informed neural networks (PINNs) to stochastic equations, thereby resulting in a new NN-based variant termed Negative Sobolev Norm-PINN (NeS-PINN).
\end{abstract}


\section{Introduction}
\label{sec:intro}

\subsection{Setting and Literature Review}
\label{ssec:setting}

In this paper, we aim to develop a Gaussian Process (GP)-based kernel method for solving roughly forced partial differential equations (PDE) of the form
\begin{equation}\label{eq: rpde}
\begin{split}
 \mathcal{P}(u) &= \xi, \quad x \in\Omega,\\
 u &= g, \quad x \in \partial\Omega
\end{split}
\end{equation}
where $\mathcal{P}$ is a (possibly nonlinear) differential operator\footnote{It is possible to generalize the setting to a non-identity boundary operator so that $\mathcal{B}(u)=g$ for $x\in\partial\Omega$; doing so is beyond the scope of this work.}.
In the context of this article, this PDE is said to be roughly forced when $\xi$ is highly irregular. For example, $\xi$ can be less regular than an $L^2(\Omega)$ function or belong to a space where the evaluation functionals $\delta_x$ are not continuous, and hence $\xi$ is not well-defined pointwise. As a result, the solution to  PDE~\eqref{eq: rpde} is also irregular and is not a classical solution. Consider as an illustrative example, the following semi-linear elliptic PDE with homogeneous boundary conditions and $\xi \in H^{-1}(\Omega)$:
\begin{subequations}
\begin{align*}
     -\Delta u + f(u) &= \xi, \quad x \in\Omega,\\
     u &= 0, \quad x \in \partial\Omega.
\end{align*}
\end{subequations}
Here, $f$ is a pointwise Nemitskii operator that is defined by lifting a nonlinear scalar function $f: \R \to \R$ to act on functions. This PDE does not have a classical solution because
$u$ does not belong to $C^2(\Omega)$.

Roughly forced PDEs arise in the context of stochastic partial differential equations (SPDE) \cite{lord2014introduction} where $\xi$ is a sample of a stochastic forcing term such as space-time white noise. SPDEs serve as useful models in many domains of application. These include: mathematical biology, where the Nagumo equation models the voltage in the axon of a neuron \cite{walsh1981stochastic, stochasticneuroscience}; physics, where the Allen-Cahn SPDE provides a model phase transitions with random fluctuations  \cite{AllenCahn, garcia1999noise}; stochastic filtering, where the Kushner-Stratonovich SPDE governs the evolution of the conditional distribution for the state given related observations \cite{fundamentalsstochasticfilter}; and probability, where Langevin equations provide a mechanism to sample measures on pathspace~\cite{stuart2004conditional,hairer2005analysis,hairer2007analysis}.
Computational methods for solving time-dependent SPDEs include finite element Galerkin methods \cite{yangalerkinfem}, spectral Galerkin methods \cite{Grecksch_Kloeden_1996}, wavelet-based methods \cite{gyongy2009rate}, and finite differences \cite{gyongy2009rate}. In the context of PDEs with random coefficients, methods include stochastic Galerkin methods \cite{ghanem1991stochastic}, Monte Carlo finite elements \cite{barth2011multi, cliffe2011multilevel}, polynomial chaos expansions \cite{xiuspdes,zhang2017numerical}, and dynamic bi-orthogonal methods based on Karhunen–Lo\`{e}ve expansions \cite{CHENG2013753, CHENG2013843}.

Machine learning PDE solvers, which trace back to statistical inference approaches for numerical approximation \cite{larkin1972gaussian}, as reviewed in \cite{owhadi2019statistical} and \cite[Chap.~20]{owhadi_scovel_2019}, 
have become increasingly prevalent due to their ability to solve a wide class of PDEs in a flexible manner. In particular, these approaches are compatible with generic engineering pipelines and have the ability to handle inverse problems and the
incorporation of previously collected data, e.g., computed solutions for different initial conditions. There are two main categories of such solvers: methods based on artificial neural networks, such as physics-informed neural networks (PINNs)~\cite{PINNs, PINNstraining}, and those based on kernels and Gaussian processes~\cite{ batlle2023error, chen2021solving}. Although GP-based methods have originally been introduced for rough linear PDEs \cite{owhadi2015bayesian, owhadi2017multigrid, owhadi2017gamblets} and have been studied in the context of probabilistic numerics \cite{cockayne2017PN-PDE, pfoertner2022PN-PDE, Wang2021PN-PDE}, their convergence guarantees required the solution to belong to the RKHS of the underlying kernel and the generalization to nonlinear PDEs in~\cite{ batlle2023error, chen2021solving} require the solution also to be classically defined.
Similarly, the original version of PINNs assumes that the solution is classically defined and while more recent ANN-variants such as the DeepRitz method \cite{DeepRitz} that is based on a variational approach, and wPINNs \cite{wPINNs} that computes weak solutions of hyperbolic PDEs, have started relaxing those regularity assumptions, an ANN-based framework for roughly forced arbitrary nonlinear PDEs is still lacking. 

In this work, we are motivated by the advantages of GP-based methods, such as simple and transparent computational and theoretical 
guarantees, e.g., their near-linear complexity \cite{chen2024sparse} and convergence results on the statistical and approximation error \cite{batlle2023error, chen2021solving}. As a result, we place our initial focus on developing a GP-based approach for solving roughly forced arbitrary nonlinear PDEs, such as SPDEs, and subsequently extend the ensuing framework to ANN-based methods, such as PINNs. 

\subsection{Contributions and Outline}

In the context of problem~\eqref{eq: rpde}, we consider the operator  $\mathcal{P} \times \mathrm{Tr} : H^{t}(\Omega) \rightarrow H^{-s}(\Omega) \times H^{t-\frac{1}{2}}(\partial \Omega)$ 
that maps the solution $u$ to the right hand side, where $\mathcal{P}$ is a nonlinear differential operator and $\mathrm{Tr}$ is the trace operator. The forcing term $\xi$ belongs to $H^{-s}(\Omega)$ and the boundary term $g$ belongs to $H^{t-\frac{1}{2}}(\partial\Omega)$. In this setting, we assume the existence of a unique solution $u^* \in H^t(\Omega)$.
To solve the PDE we first reformulate \eqref{eq: rpde} as an infinite-dimensional optimal recovery problem within a reproducing kernel Hilbert space (RKHS). This space will be potentially misspecified in the sense that the true solution $u^*$ is not necessarily contained in it. This reformulation, which we term the \textit{continuum problem}, seeks to identify an optimal approximation to the solution in the RKHS that accounts for the misspecification by incorporating additive singular noise in the source term of \eqref{eq: rpde}.
We then formulate a \textit{discretized problem} by (1) approximating the negative Sobolev norms using test functions and (2) solving the corresponding optimal recovery problem through Gauss-Newton iteration. This process reduces to solving a sequence of weakly defined linearized PDEs with the kernel corresponding to the underlying misspecified RKHS. In this section, we describe the continuum and discretized problem followed by the contributions and outline of this work.

\paragraph{Continuum problem}
Consider a kernel $K: \overline{\Omega} \times \overline{\Omega} \rightarrow\R$ and its associated RKHS, which we denote by $\mathcal{H}_K$. Given a regularization parameter $\gamma >0$, we approximate the solution to PDE \eqref{eq: rpde} by a misspecified element of $\mathcal{H}_K$ that is obtained from a minimizer of the following  \textit{continuum loss} in~\eqref{eq: infinite data}. The loss relaxes the PDE constraint using a negative Sobolev $H^{-s}$-norm rather than the usual (possibly empirical) $L^2$ norm that is commonly employed in Scientific Machine Learning. The resulting continuum optimization problem is given by
\begin{subequations}\label{eq: infinite data}
    \begin{align}
        \inf_{u \in \mathcal{H}_K} &\norm{\mathcal{P}(u) - \xi}^2_{H^{-s}} + \gamma\norm{u}^2_{\mathcal{H}_K} \\
        & \text{s.t. } u(x) = g(x) \quad x \in \partial\Omega.
\end{align}
\end{subequations}
For ease of presentation, we assume the boundary condition is regular enough to be satisfied by an element of the RKHS\footnote{The ensuing methodology can be extended to the setting where the boundary condition only applies in a trace sense by relaxing the boundary condition as well.}.
We show under mild assumptions that for any $\gamma>0$, problem \eqref{eq: infinite data} has a minimizer $u^{\gamma}\in \mathcal{H}_K$ and that any minimizing sequence has a convergent subsequence to a minimizer in  $\mathcal{H}_K$. Furthermore, $u^{\gamma}$  converges to the solution of the PDE \eqref{eq: rpde} as the regularization parameter $\gamma$ goes to $0$. The use of a negative Sobolev $H^{-s}$-loss, which is weaker than the usual $L^2$ loss, was already explored in the context of linear elliptic PDEs in \cite{bonito2024convergence}.

\paragraph{Discretized problem}

Problem \cref{eq: infinite data} is an optimization problem over the space $\mathcal{H}_K$ and cannot be directly solved computationally. We introduce a \textit{discretized} approximation of this problem based on an $N$-dimensional test space $\Phi^N :=\text{span} \{ \varphi_i\}_{i=1}^N \subset H^s_0(\Omega)$ and a set of $M$ collocation points on the boundary $\partial\Omega$. That is, the discretized problem is given by
\begin{subequations}
\label{eq: constrained problem approx}
\begin{align}
    \inf_{u\in \mathcal{H}_K} &|\mathcal{P}(u) - \xi|^2_{\Phi^N} + \gamma \|u\|^2_{\mathcal{H}_K}, \\ 
    &\text{s.t. } u(x_j) = g(x_j) \quad\text{for }j=1,\dots,M.
\end{align}
\end{subequations}
For simplicity, we have assumed $g$ to be continuous on the boundary\footnote{As for the continuum problem, the proposed approach can naturally be extended to the setting where the boundary condition has to be defined in a weak/trace sense.}.
We will show that the seminorm $|\cdot|_{\Phi^N}$ approximates the negative Sobolev norm $\| \cdot \|_{H^{-s}}$, and that the minimizer $u^{\gamma, N, M}$ of~\cref{eq: constrained problem approx} converges to the true solution of the PDE \eqref{eq: rpde} as the regularization $\gamma \to 0$  and as the number of measurements (given by test functions and collocation points) satisfy $(M, N) \to \infty.$  We then formulate a Gauss-Newton iteration to solve the discretized problem, which leads to a computational method to solve roughly forced PDEs of the form in \eqref{eq: rpde}. The formal derivation, theoretical analysis and numerical implementation of this algorithm will be the object of this paper. The main contributions are outlined next.

\paragraph{Contributions} In this work we provide the following theoretical contributions:
\begin{contrib}[label=(T\arabic*)]
    \item 
    \label{contrib:norm}
    We propose a loss based on a negative Sobolev norm, which is weaker than the $L^2$ norm and adapted to the problem. Minimizing this loss can be interpreted as solving a weak form of the PDE and can be applied to general machine-learning solvers, including kernel methods and PINNs. 
    \item 
    \label{contrib:convergenceT}
    We show that under mild assumptions on the problem, the proposed kernel method is asymptotically convergent as the regularization parameter $\gamma$ goes to $0$, hence allowing the recovery of solutions that are very irregular. We thus extend the guarantees present in \cite{chen2021solving} beyond the setting of classical solutions. This approach can accommodate both linear and nonlinear PDEs.
    \item We introduce a new concept for the solution of PDEs defined by irregular data, allowing regularization of each term in the PDE via a Gaussian process, and resulting in an optimization problem which includes the proposed Sobolev loss as a special case.
    \label{contrib:new}
\end{contrib}
From a computational perspective, we provide the following contributions:
\begin{contrib}[label=(C\arabic*)]
    \item 
    \label{contrib:norm_discretization}    
    We propose an efficient numerical discretization of the Sobolev norm, which inherits the complexity of state-of-the-art solvers for dense kernel matrices \cite{owhadi_scovel_2019}.
    \item 
    \label{contrib:Gauss_Newton}
    We provide an efficient numerical method to minimize the resulting loss based on the Gauss-Newton algorithm.
    \item \label{contrib:convergence_guarantees}
    We provide asymptotic convergence guarantees for the resulting numerical method,
    extending the convergence results in \cite{chen2021solving, batlle2023error} to the case where the solution does not belong to the RKHS underlying the optimization problem.
    \item 
    \label{contrib:kernel}
    We demonstrate the competitiveness of the kernel-based framework for roughly forced PDEs on a range of numerical examples. 
    \item 
    \label{contrib:PINN}
    We demonstrate the flexibility of the ideas presented by using the loss in the context of PINNs applied to a range of roughly forced PDEs.
\end{contrib}

\paragraph{Outline} The remainder of the  paper is structured as follows. In  \cref{subsec:notation} we introduce notation that will be used throughout. \cref{sec: infinite data problem}  presents the formal kernel method in the continuum setting for minimization problem \eqref{eq: infinite data} and investigates the theoretical properties of its minimizer. Hence, this section relates to \cref{contrib:norm,contrib:convergenceT}. In \cref{sec: numerical method} we derive a discretization of the continuum problem, corresponding to \eqref{eq: infinite data}, resulting in an implementable numerical method with asymptotic convergence guarantees. Hence, this section relates to \cref{contrib:norm_discretization,contrib:Gauss_Newton,contrib:convergence_guarantees}. In \cref{sec: numerical experiments} we present numerical experiments demonstrating the performance of the proposed kernel method  and a PINN variant using the proposed equation-adapted norm on several nonlinear PDEs; we compute empirical convergence rates and compare the methodology to finite element and spectral Galerkin methods, thereby
addressing \cref{contrib:kernel,contrib:PINN}. Our code to reproduce the numerical experiments is available in a public repository \footnote{\url{https://github.com/MatthieuDarcy/Kernel-NeuralNetworks-SPDEs}}. Finally, using the same kernel-based framework we have developed, \cref{sec:new_concept} introduces a 
\textit{generalized continuum problem} involving regularization of nonlinearities and terms under the action of differential operators; this section concerns \cref{contrib:new}. 

\subsection{Notation}
\label{subsec:notation}

We denote the positive integers and non-negative integers
respectively by $\N=\{1,2,\cdots\}$ and $\Z^+=\{0,1,\cdots\}$
and use the notation $\R=(-\infty,\infty)$ for the reals. 
For a set $D$, we use  $\overline{D}$ to denote the closure of the set, i.e. the union of the set itself and the set of all its limit points. 
Throughout,  $\Omega\subset\R^{d}$ will be a bounded open domain with Lipschitz boundary $\partial\Omega$. We let $C(\overline{\Omega})$ denote the infinite-dimensional Banach space of continuous functions taking inputs in $\Omega$. The space is endowed with the supremum norm. Similarly $C^r(\overline{\Omega})$ for $r\geq 0$ will denote the space of $r$-times continuously differentiable functions over $\Omega$. For $p\in [1,\infty)$ we let $L^p(\Omega)$ denote the infinite dimensional space of  $p$-integrable functions over $\Omega$. The space is endowed with the $L^p$ norm.  For $p\in [1,\infty)$ and $t>0$, $W^{t,p}(\Omega)$ denote the (fractional, if $t\notin\N$) Sobolev space of functions. We recall that setting $p=2$ defines a Hilbert space which we denote by $H^t(\Omega)$. We define the space $H^t_0(\Omega)$ as the closure under the $H^t(\Omega)$ norm of the set of smooth compactly supported functions on $\Omega$. We use the notation $H^t_g(\Omega)$ to denote the space of functions $u\in H^t(\Omega)$ such that applying the trace operator gives $\mathrm{Tr}(u) = g\in H^{t-1/2}(\partial \Omega)$. For ease of presentation, we will keep $g$ fixed and only consider $\mathcal{P}: H^t_g(\Omega) \rightarrow H^{-s}(\Omega)$, keeping the trace operator implicit. Throughout the article, we will use the calligraphic notation $\mathcal{H}_K$ to denote the reproducing kernel Hilbert space (RKHS) associated to the kernel $K$.

For notational convenience, we will drop  $\Omega$ for each of the introduced function spaces as this is implied throughout. Furthermore, we let $\norm{\cdot}_{\mathcal{U}}$ denote the norm that is endowed with each space $\mathcal{U}$, and similarly let $\langle \cdot,\,\cdot\rangle_{\mathcal H}$ denote the associated inner product to each inner product space $\mathcal{H}$. Convergence of a sequence $(u_n){n\in \N}\subset \mathcal{U}$ in the strong topology given by the norm on $\mathcal{U}$ to an element $u\in\mathcal{U}$ will be denoted by $u_n\to u$, while weak convergence will be denoted by $u_n\rightharpoonup u$. For a normed space $\mathcal{U}$, we will let $B_{\mathcal{U}}(f, \delta)\subset \mathcal{U}$ denote the open ball in $\mathcal{U}$ of radius $\delta$ centered at $f\in \mathcal{U}$ i.e.,
\[
    B_{\mathcal{U}}(f, \delta) := \{ g \in H^{-s} : \norm{f-g}_{\mathcal{U}} < \delta\}.
\]
For function spaces $\mathcal{U}$ and $\mathcal{V}$ we write $\mathcal{U}\hookrightarrow\mathcal{V}$ to denote $\mathcal{U}$ being compactly embedded in $\mathcal{V}$.
\section{The Continuum Problem}
\label{sec: infinite data problem}

We begin our discussion by recalling the PDE \eqref{eq: rpde} where $\mathcal{P}: H^{t}_g\to H^{-s}$ for $t > \frac{1}{2}, s > 0$ is a possibly nonlinear differential operator. We assume that there exists a unique $u^*\in H^{t}_g$ satisfying the PDE \eqref{eq: rpde}. 
\begin{remark} \label{rem:refer}
    In general, the boundary condition should be interpreted in a trace sense. If $t > \frac{d}{2}$, however, then the boundary condition can be interpreted pointwise. We will make this assumption, but we note that the methodology can also be extended to the case where the boundary condition is defined in a trace sense by enforcing the condition in a weak sense.
\end{remark}
We recall that $\mathcal{H}_K$ is an RKHS, with associated kernel $K: \overline{\Omega} \times \overline{\Omega} \rightarrow\R$. For $\gamma \in \R^+$, we consider the regularized minimization problem 
\begin{subequations}
\label{eq: nonlinear reg problem}
    \begin{align}
        \inf_{u \in \mathcal{H}_K} &\norm{\mathcal{P}(u) - \xi}^2_{H^{-s}} + \gamma\norm{u}^2_{\mathcal{H}_K},\\
        \qquad\text{s.t. } & u = g \text{ on }\partial\Omega.
    \end{align}
\end{subequations}
We introduce the following assumptions on the differential operator $\mathcal{P}$ and the RKHS $\mathcal{H}_K$.

\begin{assumptionn}[Assumption P]
\label{assumptionn:P}
The differential operator $\mathcal{P}$ satisfies:
\vspace{0.5em}
    \begin{assumpenum}[label=(P\arabic*)]
        \item 
        \label{assumption: cont} 
        The operator $\mathcal{P}: H^t_g \to H^{-s}$ is continuous.
        \item
        \label{assumption: local stability}
   The solution operator associated to $\mathcal{P}$ is locally stable. That is, there exists a $\delta\in\R^+$ and a constant $C(\delta) \in \R^+$ such that for any  $u \in H^{t}_g$ with $\mathcal{P}(u) \in B_{H^{-s}}(\xi, \delta)$, then 
   \begin{align*}
       \norm{u - u^*}_{H^t_0} \leq C(\delta)\norm{\mathcal{P}(u) - \xi}_{H^{-s}}.
   \end{align*}
    \end{assumpenum}
    
\end{assumptionn}
   
\renewcommand{\theassumptionn}{H}
\begin{assumptionn}[Assumption H]
\label{assumptionn:H}
The kernel $K: \overline{\Omega} \times \overline{\Omega} \rightarrow\R$ is such that the RKHS $\mathcal{H}_K$ satisfies:
\vspace{0.5em}
    \begin{assumpenum}[label=(H\arabic*)]
        \item
        \label{assumption: compact_embedding_HK}
        The embeddings $\mathcal{H}_K  \hookrightarrow  C(\overline{\Omega})$ and $\mathcal{H}_K\hookrightarrow H^{t}$ are compact.
        \item 
        \label{assumption: density RKHS}
    $\mathcal{H}_K$ contains a dense subset of  $H^{t}_g$ \footnote{This assumption can be replaced by the more restrictive one that $\mathcal{H}_K$ is such that $C^{\infty}  \subset \mathcal{H}_K$. This assumption is often more easily checked and is the case for example with the Mat\'{e}rn class of kernels \cite{matern_rkhs}.}.
    \end{assumpenum}
\end{assumptionn}

The following discussion will focus on identifying the assumptions under which the recovery of $u^*$ can be achieved using a scheme based on the regularized problem in \eqref{eq: nonlinear reg problem}. Specifically, in Subsection \ref{subsec: nonlinear - existence}, we demonstrate that under \cref{assumption: cont,assumption: compact_embedding_HK}, a solution to \eqref{eq: nonlinear reg problem} exists for any $\gamma \in \mathbb{R}^+$. In Subsection \ref{subsec: nonlinear - convergence} we show, under \cref{assumptionn:P} and \cref{assumptionn:H}, $(u_n)_{n\in\mathbb{N}} \subset \mathcal{H}_K$ solving \eqref{eq: nonlinear reg problem} for $\gamma=\gamma_n \in \mathbb{R}^+$, converges to $u^* \in H^{t}_g$ solving \eqref{eq: rpde} as the sequence $(\gamma_n)_{n\in\mathbb{N}} \subset \mathbb{R}^+$ decreases to zero.

\subsection{Existence of a Minimizer}
\label{subsec: nonlinear - existence}
In this subsection we prove the existence of a solution to the regularized optimization problem in \eqref{eq: nonlinear reg problem} for any $\gamma\in\R^+$. Define the objective
\begin{align}\label{eq: objective functional}
    J(u ; \gamma) \coloneqq \norm{\mathcal{P}(u) - \xi}^2_{H^{-s}} + \gamma  \norm{u}^2_{\mathcal{H}_K}, 
\end{align}
and consider the infimization problem 
\begin{subequations}
\begin{align}\label{eq: nonlinear unconstrained reg problem}
        &\inf_{u \in \mathcal{H}_K}  J(u; \gamma) \\
        &\text{s.t. } u = g \text{ on }\partial\Omega.
\end{align}
\end{subequations}

We assume that there is at least one $u\in \mathcal{H}_K$ which satisfies the boundary condition (note that $\mathcal{H}_K  \hookrightarrow  C(\overline{\Omega})$ in \cref{assumption: compact_embedding_HK} implies that  $\mathcal{H}_K$ is  of sufficient regularity to enforce the boundary condition pointwise).

We are now able to state and prove the following existence theorem.
\begin{theorem}[Existence of a Minimizer]\label{thm: existence of a minimizer}
   Let \cref{assumption: cont} be satisfied and let the kernel $K$ be such that \cref{assumption: compact_embedding_HK} is satisfied.   Let $m$ be the optimal value of problem \eqref{eq: nonlinear unconstrained reg problem}, i.e.,
   \begin{subequations}
   \begin{align*}
       m \coloneqq   &\inf_{u\in\mathcal{H}_K} J(u; \gamma)\\
        &\text{s.t. } u = g \text{ on }\partial\Omega.
   \end{align*}
   \end{subequations}
    Then it holds that $0 \leq  m  < \infty$ and that there exists  $u_\infty \in \mathcal{H}_K$ such that $m = J(u_\infty;\gamma)$. Moreover, if $(u_n)_{n\in\N}$ is a minimizing sequence, then there is a subsequence which converges strongly to $\bar{u}_\infty$ (possibly distinct from $u_\infty$) in $\mathcal{H}_K$ such that $m=J(\bar{u}_\infty;\gamma)$. Therefore, if the minimizer of $J$ is unique then the entire sequence $u_n$ converges to $u_\infty$.
\end{theorem}

\begin{proof}
    To prove the existence of a minimizer, we first show the boundedness of the minimizing sequence and then proceed by showing weak lower semi-continuity of the objective, thus proving the claim. 
    We reformulate our optimization problem as 
\begin{align*}
        \inf_{u \in \mathcal{H}_K} L(u; \gamma) := \inf_{u \in \mathcal{H}_K} J(u; \gamma) +  \chi(\sup_{x \in \partial\Omega} |u(x) -g(x) |)
\end{align*}
where
\begin{equation*}
    \chi(a) = \begin{cases}
        0 &\text{ if } a= 0, \\
        \infty &\text{ otherwise.}
    \end{cases}
\end{equation*}
    
    \textit{Step 1: Existence of weak limit}.
    We first observe that $L(u; \gamma) \in [0,\infty)$  for all $u \in \mathcal{H}_K$ satisfying the boundary condition. Therefore $0 \leq m < \infty$. Consider some minimizing sequence $(u_n)_{n\in\N}$ so that
    \begin{equation*}
        L(u_n; \gamma) \to m \quad \text{as } n \to \infty. 
    \end{equation*}
    By the definition of $L$, the minimizing sequence is bounded, i.e., there exists $C\in\R^+$ so that
    \begin{equation*}
        \norm{u_n}_{\mathcal{H}_K} \leq C,
    \end{equation*}
    for any $n\in\N$. By the Banach–Alaoglu theorem \cite{alaoglu1940weak}, since $\mathcal{H}_K$ is reflexive, there exists $u_\infty\in \mathcal{H}_K$ and a subsequence $(u_{n_k})_{k\in\N}$ such that $u_{n_k} \rightharpoonup u_\infty$ in $\mathcal{H}_K$ as $k\to\infty$. 

    \textit{Step 2: Weak lower semicontinuity of the objective.} By \cref{assumption: compact_embedding_HK} and by \cite[p.171]{brezis_2010}
    we deduce that $u_{n_k} \to u_\infty$ in $g+{H}_{0}^t$ as $k\to\infty$. Therefore, by \cref{assumption: cont}, we obtain
    \begin{equation*}
        \norm{\mathcal{P}(u_{n_k}) - \xi}^2_{H^{-s}} \to \norm{\mathcal{P}(u_\infty) - \xi}^2_{H^{-s}} \quad \text{as}\;k\to\infty.
    \end{equation*}
    We note that by the same argument $L$ is also weakly continuous on $\mathcal{H}_K$. Likewise, by \cref{assumption: compact_embedding_HK}, it follows that $u_{n_k} \to u_\infty$ in $C(\overline{\Omega})$ as $k\to\infty$. Since $u_{n_k}$ satisfies the boundary conditions for all $k$, this implies that $u_\infty = g$ for $x\in \partial\Omega$ and 
    \begin{equation*}
       \chi(\sup_{x \in \partial\Omega} |u_{n_k}(x) -g(x) |) = \chi(\sup_{x \in \partial\Omega} |u_\infty(x) -g(x) |) = 0.
    \end{equation*} 
  By the same argument this term is also lower-semicontinuous. 
    Since norms are lower-semicontinuous \cite[p. 61]{brezis_2010}, it follows that $L: \mathcal{H}_K \to \R$ is lower-semicontinuous. Therefore
    \begin{equation*}
        u_{n_k} \rightharpoonup u_\infty \Rightarrow \liminf_{n \to \infty} L(u_{n_k}; \gamma) \geq L(u_\infty; \gamma),
    \end{equation*}
    and consequently, since $\lim_{n \to \infty} L(u_n; \gamma) = m$, we have
    \begin{equation*}
        m = L(u_{\infty}; \gamma).
    \end{equation*}
    
    \textit{Step 3: Strong convergence of the minimizing sequence.} We now turn our attention to showing the existence of a subsequence of $(u_n)_{n\in\N}$ that converges strongly to a minimizer. For ease of exposition, let 
    \begin{equation*}
        \Phi(u) \coloneqq \norm{\mathcal{P}(u) - \xi}^2_{H^{-s}} + \chi(\sup_{x \in \partial\Omega} |u(x) -g(x) |),
    \end{equation*}
    so that 
    \begin{equation*}
        L(u) =  \Phi(u) + \gamma \norm{u}_{\mathcal{H}_K}^2.
    \end{equation*}
  Using the polarization equality, we have that
    \begin{align*}
        \frac{\gamma}{4} \norm{u_l - u_k}^2 &=   \gamma \Big( \frac{1}{2}\norm{u_l}^2 + \frac{1}{2}\norm{u_k}^2 - \norm{\frac{1}{2}\left(u_l + u_k \right)}^2\Big) \\
        &= \frac{1}{2}L(u_k) + \frac{1}{2}L(u_l) - L\bigg(\frac{1}{2}\big(u_l + u_k \big)\bigg) - \frac{1}{2}\Phi(u_l) - \frac{1} {2}\Phi(u_k) + \Phi\bigg(\frac{1}{2}\big(u_l + u_k \big)\bigg), 
    \end{align*}
    for any $l,k\in\N$. Letting $l,k \geq N_1(\delta)$ so that $L(u_k)\leq m+\delta$ and $L(u_l)\leq m+\delta$, we obtain
    \begin{subequations}
        \begin{align*}
             \frac{\gamma}{4} \norm{u_l - u_k}^2 &\leq m + \delta - m- \frac{1}{2}\Phi(u_l) - \frac{1} {2}\Phi(u_k) + \Phi\bigg(\frac{1}{2}\big(u_l + u_k \big)\bigg) \\
             & = \delta - \frac{1}{2}\Phi(u_l) - \frac{1} {2}\Phi(u_k) + \Phi\bigg(\frac{1}{2}\big(u_l + u_k \big)\bigg).
        \end{align*}
    \end{subequations}
    Since $u_l$, $u_k$ and $\frac{1}{2}\big(u_l + u_k \big)$ converge to $u_\infty$ in $H^t$ and in $C(\overline{\Omega})$, and since $\Phi$ is continuous on those two spaces, we deduce that there is integer $N_2(\delta)$ such that for $k,l \geq N_2(\delta)$:
    \begin{align*}
        \frac{\gamma}{4} \norm{u_l - u_k}_{\mathcal{H}_K}^2 &\leq \delta +\delta = 2\delta.
    \end{align*}
    Therefore, the sequence $(u_n)_{n\in\N}$ is Cauchy in $\mathcal{H}_K$ and as a consequence converges strongly, i.e.,
    $u_n \rightarrow u_\infty$  in $\mathcal{H}_K$ as $n\to\infty$. 
\end{proof}

\subsection{Convergence as Regularizer is Removed}
\label{subsec: nonlinear - convergence}
We now turn our attention to showing that the sequence  $(u_n)_{n\in\N}\subset \mathcal{H}_K$ where $u_n$ is the solution to \eqref{eq: nonlinear reg problem} for $\gamma=\gamma_n \in \R^+$, converges to $u^*\in H^{t}_0$ for any sequence $(\gamma_n)_{n\in\N}\subset \R^+$ decreasing to $0$. This key result lays the groundwork for the development of an implementable numerical method for the solution of \eqref{eq: nonlinear reg problem}. To prove the convergence of the numerical method to the true solution, previous work on the use of Gaussian processes to solve PDEs \cite{chen2021solving} assumes that the true solution belongs to the RKHS, i.e., $u^* \in \mathcal{H}_{K}$. In this work we treat the case where $u^* \notin \mathcal{H}_K$; extending to this setting is natural for roughly forced PDEs where, due to the lack of regularity of the right-hand side, the solution may live outside the RKHS defined by the kernel. 
\begin{theorem}[Convergence]
\label{thm:convergence to solution}
Let \cref{assumptionn:P} and 
\cref{assumptionn:H} hold and let $\xi\in H^{-s}$. For any $n\in\N$, let $u_{n}$ denote a minimizer of \eqref{eq: nonlinear reg problem} with regularizer $\gamma=\gamma_n$. Then, for all sequences $(\gamma_n)_{n\in\N} \subset \R^+$ decreasing to 0, the sequence $(u_{n})_{n\in\N} \subset \mathcal{H}_K$ converges to $u^*\in H^{t}_g$ in $H^t$ and we have 
    \begin{equation*}
        \norm{\mathcal{P}(u_n) -\xi }^2_{H^{-s}} + \gamma_{n}\norm{u_{n}}^2_{\mathcal{H}_K} \to 0\,.
    \end{equation*}
\end{theorem}
\begin{proof}
We will need the following lemma whose proof is straightforward.
\begin{lemma}[Density]\label{lemma: density}
    Under~\cref{assumption: cont,assumption: density RKHS}, the image of $\mathcal{H}_K$ under $\mathcal{P}$ denoted as $\mathcal{P}(\mathcal{H}_K)$ is dense in $\mathcal{P}(g + H_0^t)$. 
\end{lemma}
From \cref{thm: existence of a minimizer} we can 
  define $u_{n}$ as a minimizer of \eqref{eq: nonlinear reg problem} with regularizer $\gamma_n$.
    \cref{lemma: density} implies that there exists a sequence $(\bar{w}_n)_{n\in\N}$ in $\mathcal{H}_K$ such that $\mathcal{P}(\bar{w}_n) \to \xi$ in $H^{-s}$ and $\bar{w}_n = g$ on $\partial\Omega$.
    For $n\geq 0$, let 
\begin{equation*}
a_n:=\max\{m \leq n \mid \|\bar{w}_m\|_{\mathcal{H}_K} \leq \|\bar{w}_0\|_{\mathcal{H}_K}+\gamma_m^{-\frac{1}{2}}  \}.
\end{equation*}
 Note that $\lim_{n\rightarrow \infty }\gamma_n= 0$ implies  $\lim_{n\rightarrow \infty }a_n= \infty$ and $\lim_{n\rightarrow \infty}\gamma_n \|\bar{w}_{a_n}\|_{\mathcal{H}_K}=0$. Therefore, if we define $w_n:=\bar{w}_{a_n}$, we deduce that $(w_n)_{n\in\N}$ is a sequence in $\mathcal{H}_K$ such that $\mathcal{P}(w_n) \to \xi$ in $H^{-s}$, $w_n = g$ on $\partial\Omega$ and
    \begin{equation*}
        \norm{\mathcal{P}(w_n) -\xi }^2_{H^{-s}} + \gamma_{n}\norm{w_{n}}^2_{\mathcal{H}_K} \to 0.
    \end{equation*}
Now, since 
\begin{equation*}
    \norm{\mathcal{P}(u_n) -\xi }^2_{H^{-s}} + \gamma_{n}\norm{u_{n}}^2_{\mathcal{H}_K} \leq \norm{\mathcal{P}(w_n)-\xi }^2_{H^{-s}} + \gamma_{n}\norm{w_{n}}^2_{\mathcal{H}_K},
\end{equation*}
we have that 
\begin{equation*}
    \norm{\mathcal{P}(u_n)-\xi}^2_{H^{-s}} \to 0.
\end{equation*}
Therefore, \cref{assumption: local stability} and the fact that $u_n$ satisfies the boundary conditions implies  
\begin{equation*}
    \norm{u_{n} - u^*}_{H^t_0} \to 0.
\end{equation*}
\end{proof}

It is of general interest to investigate the rate at which the solution to the continuum problem \eqref{eq: nonlinear reg problem} converges to the true solution $u^*$ as the regularizer $\gamma\to 0$. We provide insight into identifying an explicit convergence rate in the following remark.

\begin{remark}[Convergence Rate for Poisson Equation]
    For the purposes of this discussion consider the  Poisson equation 
\begin{subequations}
\label{eq:Poisson}
    \begin{align}
        -\Delta u &= \xi, \quad x \in \Omega \\
        u &= 0, \quad x \in \partial\Omega\,.
    \end{align}
\end{subequations}
Assume the solution $u^*$ of \eqref{eq:Poisson} to live in  $H^t_0$, and  $\xi$ in $H^{-s}$. We may make the specific choice of RKHS $\mathcal{H}_K = H^q$ (e.g., by letting $K$ be a Mat\'{e}rn kernel), for some $q>t$, so that the solution space is not contained in the RKHS. In~\cref{app: convergence rates} we demonstrate the convergence of $u_\gamma$ to $u^*$ as $\gamma\to 0$ in the $H^\alpha$ norm for any $\alpha \leq t = 2-s$. The proof of~\cref{prop:conv_linear} for the convergence of the solution to the Poisson equation, and the constants therein, lead to the conjecture that $u_\gamma$ converges to $u^*$ as $\gamma\to 0$ in $H^{2-s}$ at the rate $\gamma^2$. We provide further details in~\cref{app: convergence rates}.
\end{remark}


\section{The Numerical Method}
\label{sec: numerical method}
In this section, we present an algorithm for solving the  constrained optimization problem: \begin{subequations}
\label{eq: constrained problem1}
\begin{align}
    \inf_{u\in \mathcal{H}_K} &\|\mathcal{P}(u) - \xi\|^2_{H^{-s}} + \gamma \|u\|^2_{\mathcal{H}_K}, \\ 
    &\text{s.t. } u = g \text{ on }\partial\Omega.
\end{align}
\end{subequations}

The section is structured as follows. First, in \cref{sec: approx sobolev norm}, we introduce a discretization of the negative Sobolev norm $\| \cdot\|_{H^{-s}}$, which reduces \eqref{eq: constrained problem1} to a nonlinear least-squares problem. This approximation of the norm  
relies on the choice of a test space. 
In \cref{sec: GN algorithm}, we relax the boundary condition to be satisfied at only a finite set of collocation points, 
and we employ a variant of the Gauss-Newton algorithm to solve the resulting infinite-dimensional minimization problem. The representer theorem allows us to obtain a finite-dimensional representation of the solution for each iteration of Gauss-Newton. Finally, in section \cref{sec: numerical convergence} we present asymptotic convergence results for the resulting numerical method. The derivation of this numerical method and its proof of convergence will be based on the following density assumption on the test space and collocation points (\cref{assumptionn:D}) and regularity assumption on the differential operator (\cref{assumption: regularity of f 2}).

\renewcommand{\theassumptionn}{D}
\begin{assumptionn}[Assumption D]
\label{assumptionn:D}
The $N$-dimensional test space $\Phi^N := \text{span} \{ \varphi_i\}_{i=1}^N \subset H^s_0$ and set of collocation points on the boundary $\{ x_j\}_{j=1}^M \subset \partial \Omega$ satisfy:
\vspace{0.5em}
    \begin{assumpenum}[label=(D\arabic*)]
        \item
        \label{assumption: density of test space}
    It holds that $\Phi^N$ becomes dense in $H^s_0$ as $N \to \infty$, i.e.,
    \[
        \overline{\lim_{N \to \infty} \text{span}  \{\varphi_i\}_{i=1}^N} = H^{s}_0.
    \]

        \item 
        \label{assumption: boundary discrete}
The fill-distance on the boundary is such that
\[
    \sup_{x \in \partial \Omega}\,\, \inf_{x_j \in \{ x_i\}_{i=1}^M}|x - x_j|_{\partial \Omega} \,\,\rightarrow 0, \quad\text{as } M \rightarrow \infty.
\]
    \end{assumpenum}
\end{assumptionn}

\renewcommand{\theassumptionn}{R}
\begin{assumptionn}[Assumption R]
\label{assumption: regularity of f 2}
    Assume that $\xi\in H^{-s+\delta}$ for some $\delta\in\R^+$, and assume $\mathcal{P}:H^{t}_g\to H^{-s+\delta}$.
\end{assumptionn}

Throughout this section we consider the assumptions of \cref{sec: infinite data problem} to also be satisfied. In particular we consider that
\cref{assumptionn:P} and \cref{assumptionn:H} are satisfied, that $g$ is of sufficient regularity to enforce the boundary condition pointwise
and that there is at least one $u\in \mathcal{H}_K$ which satisfies the boundary condition.

\subsection{Approximating The Negative Sobolev Norm}
\label{sec: approx sobolev norm}
We now derive a method for approximating the $H^{-s}$ norm in problem~\eqref{eq: constrained problem1}. First,  recall from~\cite{henry2006geometric} that one can define a norm (equivalent to the usual $\norm{\cdot}_{H^s}$ norm) on the space $H^{s}_0$ as 
\begin{align}\label{eq: energy norm}
    \norm{u}_{H^s_0}^2 = \int u (-\Delta)^s u, 
\end{align}
where $\Delta$ is the Laplacian operator with homogeneous Dirichlet boundary conditions. One can additionally consider fractional values of $s$ through the spectral expansion of the negative Laplacian. The negative Sobolev norm $\norm{\cdot}_{{H^{-s}}}$ is then defined as the dual norm to the space $H^{s}_0$ with respect to the $L^2$ duality
\begin{equation} \label{Hms_norm}
    \norm{f}_{{H^{-s}}}:= \sup_{u \in H^s_0}\frac{\int f u}{\norm{u}_{H^s_0}}.
\end{equation}

This supremum in~\eqref{Hms_norm} cannot be directly computed as it involves solving an optimization problem over an infinite-dimensional space. Our strategy is to replace the infinite-dimensional space $H^{s}_0$ by a finite-dimensional test space: 
\begin{equation*}
    \Phi^N := \text{span} \{ \varphi_i\}_{i=1}^N \subset H^{s}_0.
\end{equation*}
Appropriate choices of $ \Phi^N$ will be discussed later. We further assume that the set $\{\varphi_i\}_{i=1}^N$ is comprised of linearly independent functions. We may now define the seminorm
\begin{equation} \label{eq:semi-norm}
    |f|_{\Phi^N} :=  \sup_{u \in \Phi^N}\frac{\int f u}{\norm{u}_{H^s_0}}.
\end{equation} 
We make use of the vector notation 
\begin{equation*}
    [f, \varphi]:= \Big(\int f \varphi_1, \int f \varphi_2, \dots, \int f \varphi_n\Big)^\intercal \in \R^N
\end{equation*}
and define the matrix $A\in \R^{N \times N}$ element-wise as
\begin{equation*}
    A_{i,j} := \int \varphi_i (-\Delta)^{s} \varphi_j,
\end{equation*}
for $1\leq i,j \leq N$. The following result provides an explicit expression for the seminorm in~\eqref{eq:semi-norm}. 
\begin{proposition}\label{prop: expression n norm}
For any $f \in H^{}-s$, the seminorm $|\cdot|_{\Phi^N}$ satisfies
\begin{equation}\label{eqnsobn}
    |f|_{\Phi^N} =  \sqrt{[f, \varphi]^\intercal A^{-1}[f, \varphi]}.
\end{equation} 
\end{proposition}

\begin{proof}
    See \cref{app: numerical method}.
\end{proof}

\begin{remark}
    Another view of our approximation to the norm comes from considering the dual of our test space $\Psi^N := \{ (-\Delta)^{s}\varphi_i\}_{i=1}^N \subset H^{-s}$. One can 
    consider the $H^{-s}$ projection of $f$ onto $\Psi^N$, denoted as $\Pi^N f$.
    Then, one can show
    \begin{align*}
        |f|_{\Phi^N} = \norm{\Pi^N f}_{H^{-s}}\leq \|f\|_{H^{-s}}
    \end{align*}
    Hence, our seminorm can be seen as the norm of the $H^{-s}$ projection of $f$ onto  $(-\Delta)^{s}\Phi^N$.
\end{remark}
The next proposition ensures the consistency of our seminorm as the number of test functions $N \rightarrow \infty$. We first make an assumption on our choice of test space.
\begin{proposition}\label{prop: norm consistency}
   Under \cref{assumption: density of test space}, it holds that for all $f \in H^{-s}$
    \begin{equation*}
        \lim_{N \to \infty} |f|_{\Phi^N} \to \norm{f}_{H^{-s}}.
    \end{equation*}
\end{proposition}
\begin{proof}
    See \cref{app: numerical method}. 
\end{proof}
Writing $\Phi^{N, \perp}$ for the orthogonal complement of $\Phi^{N}$ in $H^{-s}$, 
we may also define the orthogonal seminorm as
\begin{align*}
     |f|_{\Phi^{N, \perp}}&:=  \sup_{u \in \Phi^{N, \perp}}\frac{\int f u}{\norm{u}_{H^s_0}}\,.
\end{align*}
Since $|f|_{\Phi^{N, \perp}}$ is the $H^{-s}$-norm of the $H^{-s}$-projection of $f$ onto $\Phi^{N, \perp}$, we have
 the orthogonal decomposition
\begin{align}\label{eq: norm orthogonal decomp}
    \norm{f}^2_{H^{-s}} = |f|^2_{\Phi^N} + |f|^2_{\Phi^{N, \perp}}
\end{align}
Combining this identity with \cref{prop: norm consistency} leads to the following corollary.
\begin{corollary}
\label{cor: norm consistency}
For $f\in H^{-s}$, the family of seminorms $|\cdot|_{\Phi^N}$ satisfy
    \begin{equation*}
        \lim_{N \to \infty} |f|_{\Phi^{N, \perp}} \to 0.
    \end{equation*}
\end{corollary}

\begin{remark} (Choice of test space $ \Phi^N$).
A standard choice for the basis functions of $ \Phi^N$
is to use polyharmonic splines \cite{wendland2004scattered}. With this selection, the stiffness matrix, its inverse, and the projected norms can be computed with near-linear complexity using gamblets \cite{owhadi_scovel_2019}. However, this approach is not explored in the present work.
\end{remark}

\subsection{Nonlinear Least-squares Problem and the Gauss-Newton Algorithm}
\label{sec: GN algorithm}
We now consider the approximation in \eqref{eqnsobn} for the negative Sobolev norm, and relax the boundary condition to be satisfied only at a finite set of collocation points $\{ x_j\}_{j=1}^M \subset \partial \Omega$. This leads to the  following minimization problem for the regularized PDE solution:
\begin{subequations}\label{eq: loss finite-dimensional}
    \begin{align} 
    \inf_{u\in \mathcal{H}_K} &|\mathcal{P}(u) - \xi|^2_{\Phi^N} + \gamma \|u\|^2_{\mathcal{H}_K}, \\ 
    &\text{s.t. } u(x_j) = g(x_j) \quad\text{for }j=1,\dots,M,
\end{align}
\end{subequations}
which by \cref{prop: expression n norm}, is equivalent to 
\begin{subequations}\label{eq: loss expression}
\begin{align}
     \inf_{u\in \mathcal{H}_K} &[\mathcal{P}(u) - \xi, \varphi]^\intercal A^{-1} [\mathcal{P}(u) - \xi,\varphi]  + \gamma \norm{u}^2_{\mathcal{H}_K}\\
     &\text{s.t. } u(x_j) = g(x_j) \quad\text{for }j=1,\dots,M.
\end{align}
\end{subequations}
Similarly to \cite{chen2021solving}, we solve~\eqref{eq: loss expression} via a variant of the Gauss-Newton algorithm, which successively linearizes the nonlinear least-squares problem~\cite{Nocedal2006}. Writing $u_{n+1} = u_n + \delta u$ for the iterations of the Gauss-Newton algorithm, the increment $\delta u$ is chosen such that 
\begin{subequations}
\begin{align*}
    \delta u  = \arginf_{v\in \mathcal{H}_K} & \,[\mathcal{P}(u_n) +\mathcal{P}'(u_n)v  - \xi, \varphi]^\intercal A^{-1} [\mathcal{P}(u_n) +\mathcal{P}'(u_n)v - \xi, \varphi] + \gamma \norm{u_n + v}^2_{\mathcal{H}_K}\\
    &\text{s.t. } u_n(x_j) + \delta u(x_j) = g(x_j) \quad\text{for }j=1,\dots,M.
\end{align*}
\end{subequations}
Equivalently, we can identify $u_{n+1}$ directly by solving
\begin{subequations}
\label{eq:infinite_dim_problem}
\begin{align}
   u_{n+1}  = \arginf_{v\in \mathcal{H}_K} &\,[\mathcal{P}(u_n) +\mathcal{P}'(u_n)(v - u_n)  - \xi, \varphi]^\intercal A^{-1} [\mathcal{P}(u_n) +\mathcal{P}'(u_n)(v - u_n) - \xi, \varphi] + \gamma \norm{v}^2_{\mathcal{H}_K}\\
   &\text{s.t. } v(x_j) = g(x_j) \quad\text{for }j=1,\dots,M,
\end{align}
\end{subequations}
The next theorem provides a representer formula for the minimizer of \eqref{eq:infinite_dim_problem}.
\begin{theorem}[Representer Theorem]
    \label{thm: representer} Let $K(\cdot, \chi_i)$ denote the function $$K(\cdot, \chi_i) := \int \bigl(\mathcal{P}'(u)(x)K(\cdot, x)\bigr)\varphi_i(x)dx.$$
    The infimizer 
\begin{subequations}
\label{eq:infinite_dim_problem thm}
\begin{align}
   \widehat{u}:=  \arginf_{v\in \mathcal{H}_K} &\,[\mathcal{P}(u) +\mathcal{P}'(u)(v - u)  - \xi, \varphi]^\intercal A^{-1} [\mathcal{P}(u) +\mathcal{P}'(u)(v - u) - \xi, \varphi] + \gamma\norm{v}^2_{\mathcal{H}_K}\\
   &\text{s.t. } v(x_j) = g(x_j) \quad\text{for }j=1,\dots,M,
\end{align}
\end{subequations}
admits a representation of the form
\begin{equation*}
    \widehat{u} = \sum_{i=1}^N\alpha_iK(\cdot,\chi_i) + \sum_{j=1}^M\beta_jK(\cdot,x_j),
\end{equation*}
for some scalars $\{\alpha_i\}_{i=1}^N, \{\beta_j\}_{j=1}^M\subset \R$. 
\end{theorem}
\begin{proof}
   See \cref{app: numerical method}.
\end{proof}
Using \cref{thm: representer} it is possible to compute the solution to  \eqref{eq:infinite_dim_problem}. The result is given in the next corollary. To describe this, we introduce the following notation: denote $K(\cdot, \varphi_i) = K(\cdot, \chi_i)$ when $1 \leq i \leq N$ and $K(\cdot, \varphi_i) = K(\cdot, x_{i-N})$ when $N+1\leq i \leq M+N$. We will also adopt the following matrix notation:
\begin{align*}
    &K(\varphi, \varphi) \in \R^{(N +M) \times (N+M)} \quad &&\text{with entries }  K(\varphi, \varphi)_{i,j} = K(\varphi_i, \varphi_j) \\
    &K(\chi, \varphi) \in \R^{N \times (N+M)} \quad &&\text{with entries }  K(\chi, \varphi)_{i,j} = K(\chi_i, \varphi_j) \\
    &K(X, \varphi) \in \R^{M \times (N+M)} \quad &&\text{with entries }  K(X, \varphi)_{i,j} = K(x_i, \varphi_j) \\
    & g(X) \in \R^{M} \quad &&\text{with entries } g(X)_i = g(x_i) 
\end{align*}

\begin{proposition}
\label{prop: representer, explicit solution}
 Let $r_n := \xi - \mathcal{P}(u_n) +  \mathcal{P}'(u_n)u_n$. Then, at each step $n$, the solution to minimization problem~\eqref{eq:infinite_dim_problem} is given by
    \begin{align*}
        u_{n+1} = \sum_{i=1}^N\alpha_iK(\cdot,\chi_i) + \sum_{j=1}^M\beta_jK(\cdot,x_j)
    \end{align*}
    where $c = (\alpha_1, \dots, \alpha_N, \beta_1, \dots, \beta_M)\in \R^{N+M}$ arises from the solution to the linear system
\begin{align*}
\begin{bmatrix}
 K(\chi, \varphi)^\intercal A^{-1} K(\chi, \varphi) + \gamma K(\varphi, \varphi) & K(X, \varphi)  \\
K(X, \varphi) &  0
\end{bmatrix}
\begin{bmatrix}
    c \\
    \nu
\end{bmatrix}
= 
\begin{bmatrix}
    K(\chi, \varphi)^\intercal A^{-1}[r_n, \varphi]  \\
    g(X)
\end{bmatrix}
\end{align*}
with $\nu\in \R^{M}$.
\end{proposition}
\begin{proof}
Equation \ref{eq:infinite_dim_problem} is a standard quadratic program with linear constraints; this is used
to facilitate the proof. See \cref{app: numerical method}.
\end{proof}

\subsection{Convergence of the Numerical Method}
\label{sec: numerical convergence}
We now establish the asymptotic convergence of the minimizer of \eqref{eq: loss finite-dimensional} to the true solution. To achieve this, we impose the following assumptions: \cref{assumption: regularity of f 2} ensures the regularity of $\xi$; \cref{assumption: density of test space} guarantees that, as the dimension of the test space increases, the true norm of $\xi$ can be accurately recovered; and \cref{assumption: boundary discrete} ensures that the boundary condition is satisfied in the asymptotic limit.
The following result demonstrates that, in the limit, the solution of the numerical problem \eqref{eq: loss finite-dimensional} converges to the true solution in $H^t_0$.
\begin{theorem}[Convergence to the True Solution]\label{thm: convergence numerical}
    Let \cref{assumptionn:D} and \cref{assumption: regularity of f 2} be satisfied. Let $u^{\gamma, N, M}\in\mathcal{H}_K$ be the minimizer of \eqref{eq: loss finite-dimensional}, then
    \begin{equation*}
       \lim_{\gamma\to 0} \limsup_{N\to\infty }\limsup_{M\to\infty}\left\|u^{\gamma, N, M}-u^*\right\|_{H^t}=0.
    \end{equation*}
\end{theorem}
\begin{proof}
Since there is at least  one $u\in \mathcal{H}_K$ which satisfies the boundary condition, for a fixed $N$ and $\gamma>0$, 
\eqref{eq: loss finite-dimensional} implies that 
$\|u^{\gamma, N, M}\|_{\mathcal{H}_K}$ is uniformly bounded in $M$. 
By the Banach–Alaoglu theorem \cite{alaoglu1940weak},  $u^{\gamma, N, M}$converges weakly in $\mathcal{H}_K$ along a subsequence (which we denote using the same index $M$) to some element $u^{\gamma, N, \infty}$. By \cref{assumption: compact_embedding_HK}, this subsequence converges strongly in $C(\overline{\Omega})$ and in $H^{t}$. Using \cref{assumption: boundary discrete}, since $u^{\gamma, N, M}$ satisfies the boundary condition on an increasingly dense set of collocation points and remains continuous on the boundary uniformly with respect to the number of those points, 
  we deduce 
  \cite{interpolation2007} that $u^{\gamma, N, \infty}|_{\partial \Omega} := \lim_{M \rightarrow \infty} u^{\gamma, N, M}|_{\partial \Omega} = g$. 

We first define 
 \begin{align*}
     J^N(u; \gamma) &:= |\mathcal{P}(u) - \xi|^2_{\Phi^N} + \gamma  \norm{u}^2_{\mathcal{H}_K}  \\
     R^N(u) &:= |\mathcal{P}(u^{\gamma, N}) - \xi|_{\Phi^{N,\perp}}^2.
 \end{align*} 
Recall the definition in equation \eqref{eq: objective functional} of the functional $J(u; \gamma)$ and observe that by the orthogonal decomposition \cref{eq: norm orthogonal decomp}, $J(u;\gamma) = J^N(u; \gamma) + R^N(u)$ for any choice of $N$ and $\gamma$.
 By using \cref{assumption: cont} and the lower-semicontinuity of norms in a similar manner as in the proof of \cref{thm: existence of a minimizer}, we obtain
\begin{align*}
J^N(u^{\gamma, N, \infty}; \gamma) \leq \liminf_{M \to \infty}  J^N(u^{\gamma, N, M}; \gamma).
\end{align*}
By the definition of $u^{\gamma,N, M}$, for any $u \in \mathcal{H}_K$ such that $u|_{\partial \Omega} = g$, we have
\begin{align*}
   J^N(u^{\gamma, N, M}; \gamma) \leq J^N(u; \gamma),
\end{align*}
from which it follows that 
\begin{subequations}
\begin{align*}
    J^N(u^{\gamma, N, \infty}; \gamma) \leq \inf_{u\in \mathcal{H}_K}& J^N(u; \gamma) \\
    & \text{s.t. } u = g \text{ on } \partial \Omega.
\end{align*}
\end{subequations}
Because $u^{\gamma, N, \infty}$ satisfies the boundary conditions, we deduce that $u^{\gamma, N, \infty}$ is a minimizer of the optimization problem \cref{eq: loss finite-dimensional}.

For ease of presentation, we will now drop the last superscript and use the notation $u^{\gamma, N}:= u^{\gamma, N, \infty}$. Recall the definition in equation \eqref{eq: objective functional} of the functional $J(u; \gamma)$ and let $u^\gamma$ be any minimizer of this functional as defined by~\cref{thm: existence of a minimizer}. Observe that
\begin{equation}
\label{eq:big ineqq}
    J^N(u^{\gamma, N}; \gamma)\leq J^N(u^{\gamma}; \gamma)\leq J(u^{\gamma};\gamma)\leq J(u^{\gamma,N};\gamma) = J^N(u^{\gamma,N}; \gamma) - R^{N}(u^{\gamma, N}).
\end{equation}
Now, by \eqref{eq:big ineqq} and the inequality
\begin{equation*}
    \gamma \|u^{\gamma, N}\|^2_{\mathcal{H}_K}\leq J(u^{\gamma}; \gamma).
\end{equation*}
we deduce that the sequence $\bigl(u^{\gamma,N}\bigr)_{N\in\N}\subset\mathcal{H}_K$ is uniformly bounded in $N$ for fixed $\gamma$, hence by the assumed continuity of $\mathcal{P}$, the sequence $\bigl(\mathcal{P}(u^{\gamma,N})\bigr)_{N\in\N}$ is uniformly bounded by some $C(\gamma)\in\R^+$ in the norm of $H^{-s+\delta}$. Then, $H^{-s+ \delta} \hookrightarrow H^{-s}$ implies that $\bigl(\mathcal{P}(u^{\gamma,N})\bigr)_{N\in\N}$ converges along a 
 subsequence  to some limit  in $H^{-s}$. Therefore, there exists a subsequence (which we denote using the same index $N$) such  that
\begin{equation*}
    h^{N} := \mathcal{P}(u^{\gamma, N}) - \xi \to h^{\infty} \in H^{-s}.
\end{equation*}
We note that
\begin{align*}
     |h^N|_{\Phi^{N, \perp}}^2 &\leq  |h^N - h^\infty|_{\Phi^{N, \perp}}^2 +  |h^\infty|_{\Phi^{N, \perp}}^2\\
     &\leq   \|h^N - h^\infty\|_{H^{-s}}^2 +  |h^\infty|_{\Phi^{N, \perp}}^2.
\end{align*}
Therefore, the fact that $h^N \to h^\infty$ in $H^{-s}$ and  Corollary \ref{cor: norm consistency} imply that 
\begin{equation*}
     |h^N|_{\Phi^{N, \perp}}^2 \to 0 \quad \text{as } N \to \infty.
\end{equation*}
Thus, it holds that $|\mathcal{P}(u^{\gamma, N}) - \xi|_{\Phi^{N,\perp}}^2\to 0$ as $N\to \infty$ and we conclude that $R^{N}(u^{\gamma,N})\to0$ as $N\to\infty$. 
Using \eqref{eq:big ineqq}, we therefore deduce that 
\begin{equation*}
\lim_{N\to\infty}J(u^{\gamma,N}; \gamma)=J(u^{\gamma}; \gamma).
\end{equation*}
Therefore, applying the result from Theorem \ref{thm:convergence to solution}, we deduce that
\begin{equation*}
\lim_{\gamma\to0}\lim_{N\to\infty}J(u^{\gamma, N}; \gamma)=0.
\end{equation*}
By the local stability assumption on $\mathcal{P}$, given in \cref{assumption: local stability}, we finally conclude that 
\begin{align*}
    \lim_{\gamma \to 0}\lim_{N \to \infty} \left\|u^{\gamma, N}-u^*\right\|_{H^{t}}=0.
\end{align*}
Therefore, $u^{\gamma, N, M}$ converges along a subsequence to $u^*$ in $H^{t}$. Since the limit $u^*$ is independent of the considered subsequence, we may therefore conclude that the entire sequence converges to $u^*$. For contradiction, suppose that this were not the case: there would be a $\varepsilon$ and a subsequence $u^{\gamma, N_k, M_k}$ such that $\left\|u^{\gamma, N_k, M_k}-u^*\right\|_{H^{t}} \geq \varepsilon$ for all $N_k, M_k$. Then, applying the previous reasoning to this subsequence would provide a further subsequence that does converge to $u^*$, thus yielding a contradiction. We therefore conclude that 
\begin{equation*} \lim_{\gamma\to0}\limsup_{N\to\infty}\limsup_{M\to\infty}\left\|u^{\gamma, N, M}-u^*\right\|_{H^{t}}=0.
\end{equation*}
where the limit in $N, M$ is now taken over the entire sequence.
\end{proof}


\section{Numerical Experiments}
\label{sec: numerical experiments}
We now investigate our proposed methodology on a variety of spatial and time-dependent PDEs.
\cref{ssec:deet} contains some overarching details of the computational methodology that apply widely throughout this section. In \cref{sec: spatial pde}, we apply our method to PDEs in one and two spatial dimensions, utilizing both a kernel method and an artificial neural network (NN) trained using a physics-informed negative Sobolev loss, which we denote by NeS-PINN. We compare the performance of the methods trained using the negative Sobolev norm against the same methods trained using the usual pointwise loss,
which correspond to the GP-PDE~\cite{chen2021solving} and PINN~\cite{PINNs} methods\footnote{This corresponds to using an $L^2$ loss, which is approximated using pointwise evaluations of the functions.}. 
 In \cref{sec: time dependent spdes} we show how to apply our method to solve time-dependent stochastic PDEs.

\subsection{Details of Computational Methodology} \label{ssec:deet}
For all numerical experiments, we use the Mat\'{e}rn kernel with smoothness parameter $\nu =\frac{5}{2}$, which is defined as 
\begin{align*}\label{eq: matern kernel} 
&K_{\frac{5}{2}}(x,y) = \Big(1 + \frac{\sqrt{5}||x -y ||}{l} + \frac{5||x -y ||^2}{3l^2}\Big)\exp\Big(- \frac{\sqrt{5}||x -y ||}{l}\Big).
\end{align*}
We find that selecting a length scale \( l \) within the range \( [0.1, 2.0] \) produces favorable results in our experiments. We note that hierarchical methods to learn the length scale may be needed in general for other applications.

Our standard neural network architecture comprises four layers: a non-trainable random Fourier layer embedding as proposed in \cite{rff_pinn}, two trainable hidden layers, and a trainable outer layer. We use the random embedding layer as our experiments show that it improves accuracy and helps the network capture the high-frequency behavior of the function. Following standard implementations, we use the hyperbolic tangent activation function~\cite{PINNstraining}. A detailed table of all architectures can be found in \cref{sec: NN architecture}.
We train the ANNs for $10^5$ iterations with the ADAM optimizer \cite{adam} using the learning rate $10^{-3}$ and an exponential decay schedule. Our code to reproduce the experiments is available in the public repository: \url{https://github.com/MatthieuDarcy/Kernel-NeuralNetworks-SPDEs}.

\subsection{Rough Spatial PDEs}
\label{sec: spatial pde}
In this section, we evaluate our method on two spatial PDEs with rough right-hand sides: a one-dimensional linear PDE in subsubsection~\ref{sec: elliptic} and a two-dimensional semilinear PDE in subsubsection~\ref{sec: semilinear_PDE}. We compare the performance of the kernel method and the PINN when they are trained using the negative Sobolev norm or the standard pointwise loss.

\subsubsection{1D Linear Elliptic PDE}
\label{sec: elliptic}
We first consider the elliptic equation with homogeneous Dirichlet boundary condition
\begin{align} \label{eq: Poisson equation}
\begin{split}
 -\nu\Delta u  + u&= \xi, \quad  x \in \Omega = (0,1), \\
    u &= 0, \quad x \in \partial \Omega,
\end{split}
\end{align}
where $\nu = 10^{-2}$. We sample the right-hand side $\xi$ using the following expansion:
\begin{equation}\label{eq: forcing term 1}
    \xi \sim \sum_{j=1}^\infty \xi_j  \sqrt{2}\sin(\pi j x) ,\quad \xi_j \sim \mathcal{N}(0,1) \; \text{i.i.d.}
\end{equation}
A fixed realization of $\xi$ belongs to $H^{-s}$ for any $s >\frac{1}{2}$  a.s.. The solution of PDE~\eqref{eq: Poisson equation} is given in closed form by 
\begin{align*}
    u^*(x) =\sum_{j=1}^\infty \frac{\xi_j}{\nu j^2\pi^2 + 1}\sqrt{2}\sin( \pi j x),
\end{align*}
where $u^* \in H^{t}$ for $t >\frac{3}{2}$ a.s.. In practice, we truncate series~ \eqref{eq: forcing term 1} to \( L = 2^{14} \) terms. The number of terms $L$ is set to be much larger than the number of measurements $N+M$ to ensure that the truncation error is negligible.

\begin{figure}[htp]
    \centering
    \begin{subfigure}[b]{0.31\columnwidth}
        \centering
        \includegraphics[width=\textwidth]{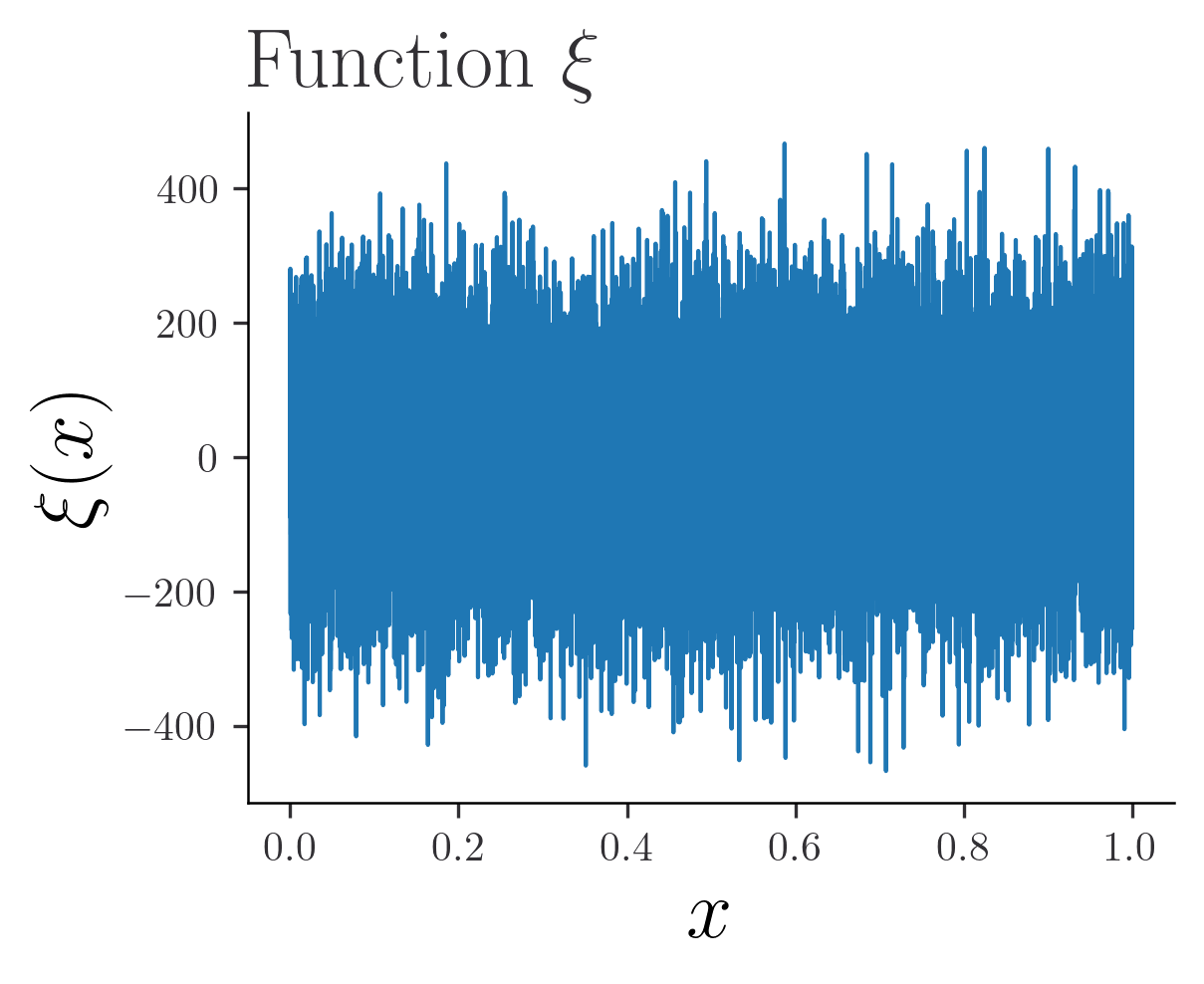}
        \caption{\footnotesize Forcing term}
    \end{subfigure}%
    \hspace{0.05\columnwidth}
    \begin{subfigure}[b]{0.31\columnwidth}
        \centering
        \includegraphics[width=\textwidth]{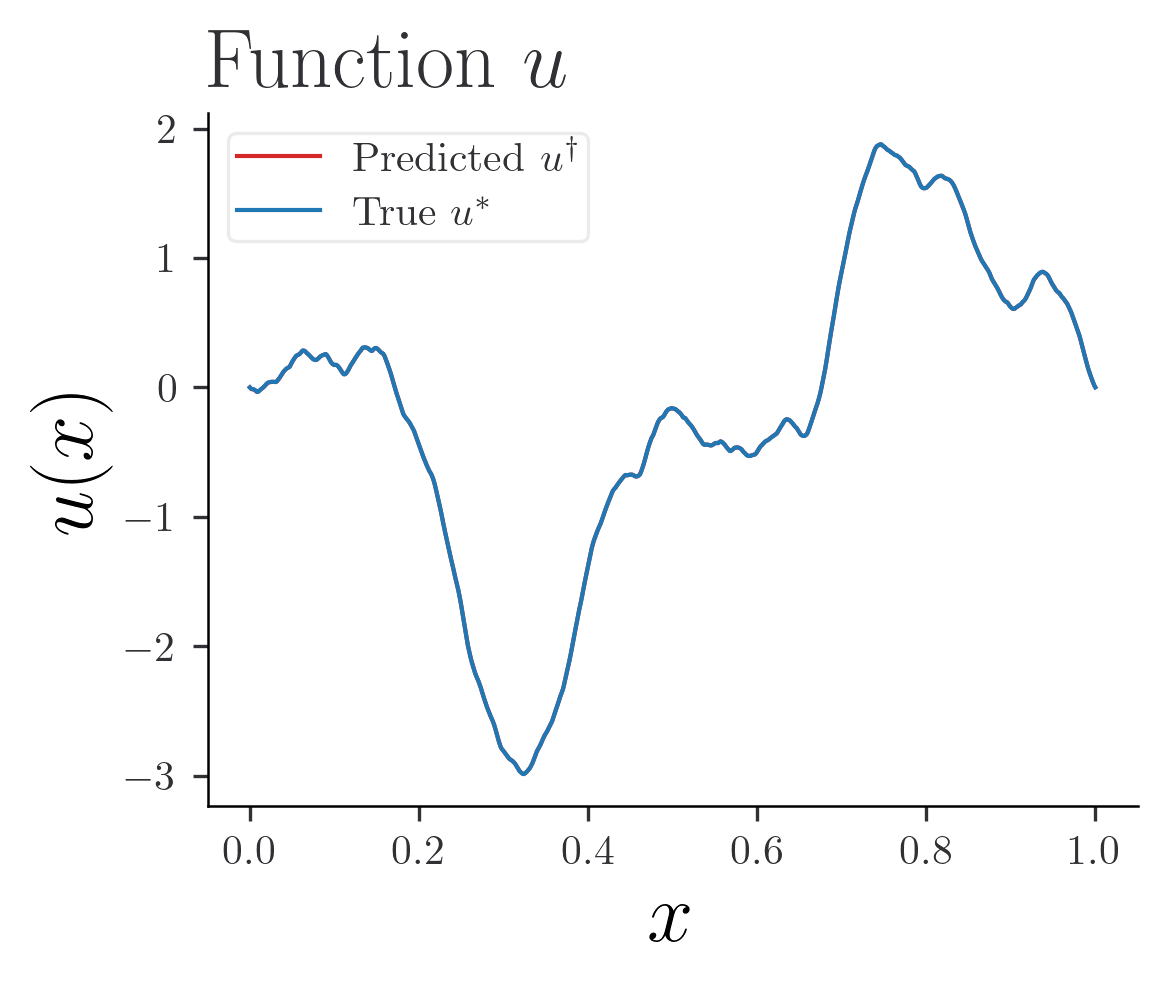}
        \caption{\footnotesize Solution}
    \end{subfigure}%

    \begin{subfigure}[b]{0.31\columnwidth}
        \centering
        \includegraphics[width=\textwidth]{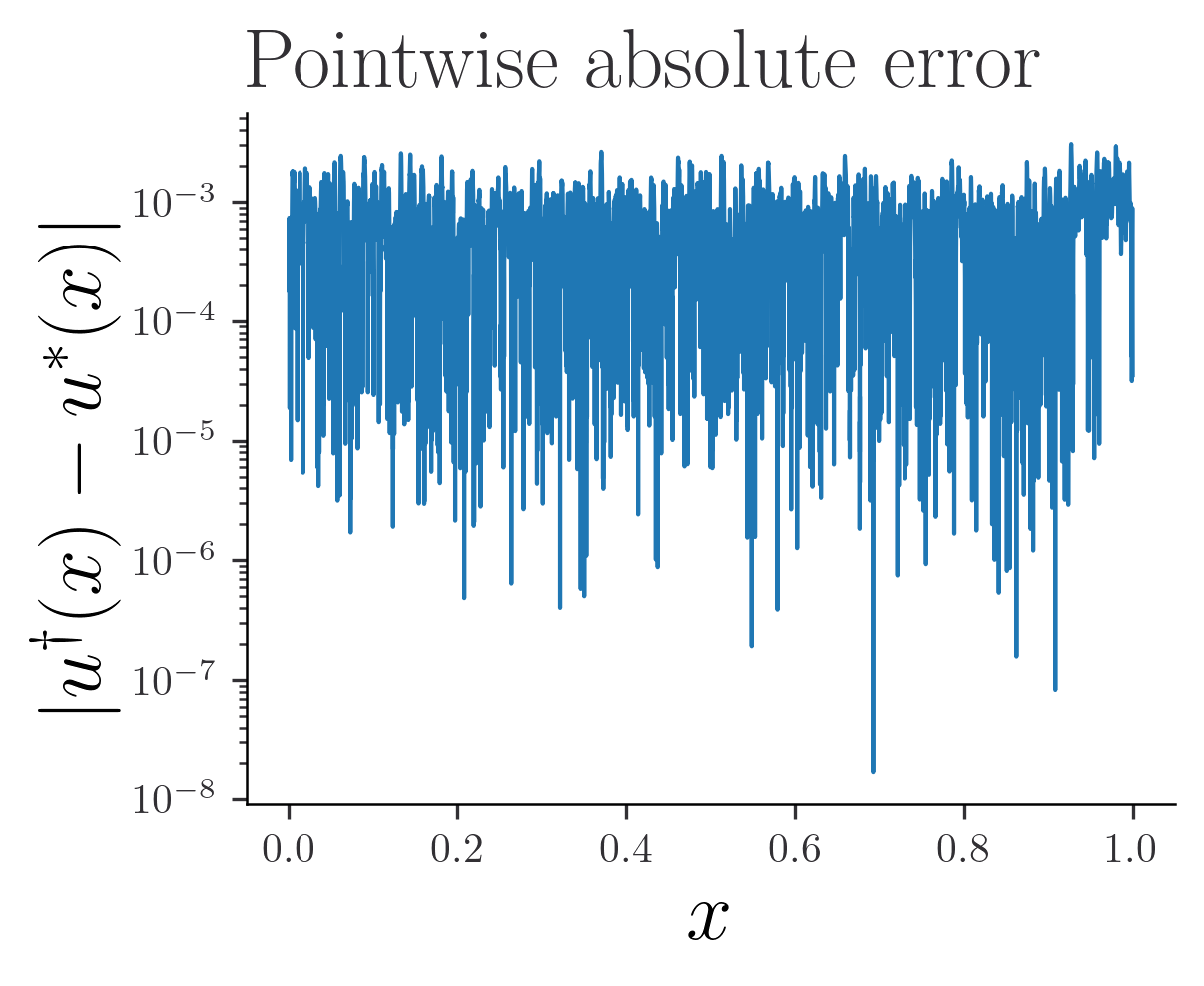}
        \caption{\footnotesize Pointwise error}
    \end{subfigure}%
    \hspace{0.05\columnwidth}
    \begin{subfigure}[b]{0.31\columnwidth}
        \centering
        \includegraphics[width=\textwidth]{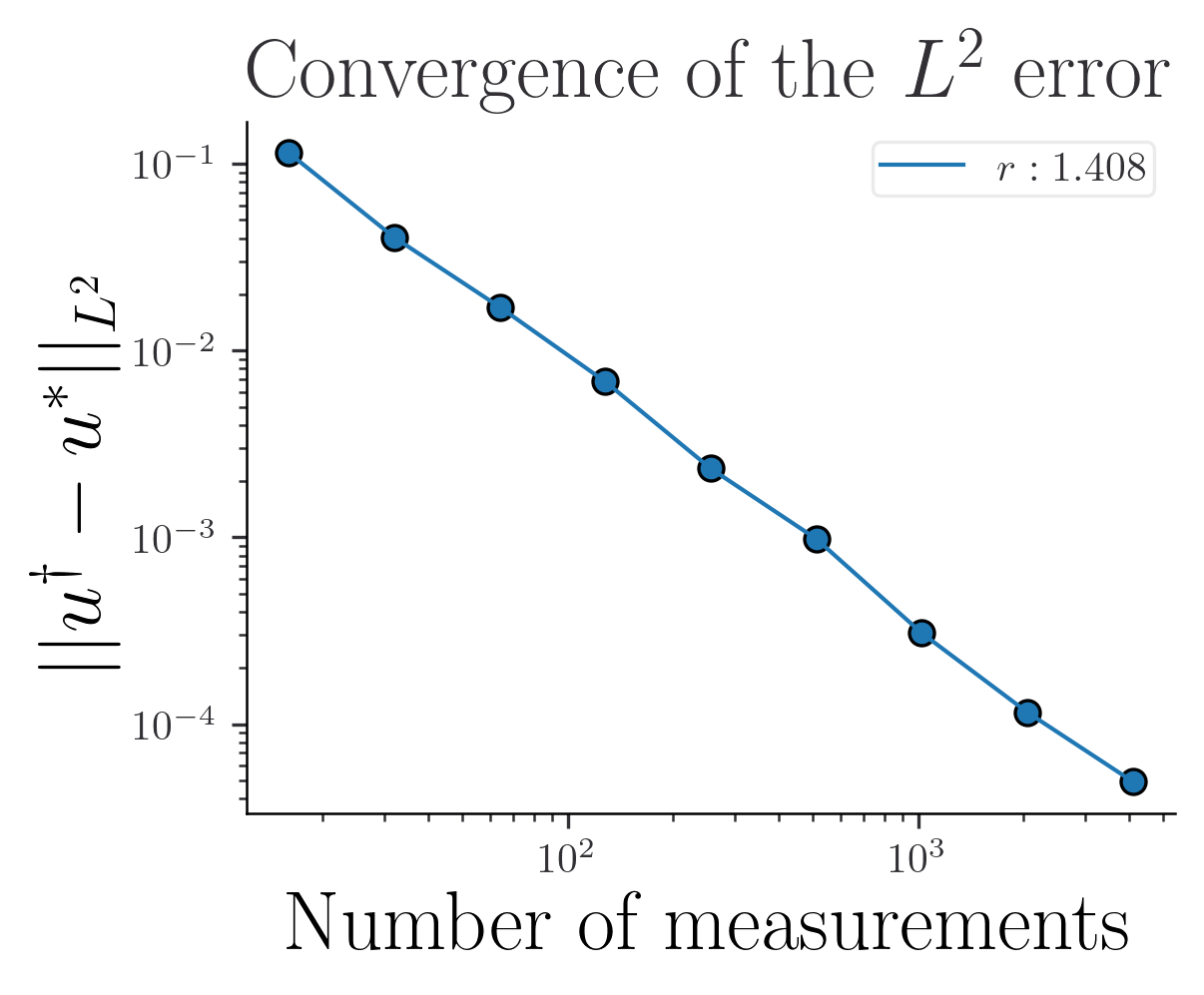}
        \caption{\footnotesize $L^2$ error}
    \end{subfigure}%
    
    \caption{Numerical results for the one-dimensional linear elliptic PDE \eqref{eq: Poisson equation}.}
    \label{fig:1d_linear_elliptic}
\end{figure}

For the choice of $t$ in the objective norm, we choose the $H^{-1}$ norm and we choose our test space to be spanned by the first $N = 4096$ basis functions $\Phi^N = \text{span}\{\varphi_j\}_{j=1}^N$ with $\varphi_j = \sin(\pi j x)$. Each measurement of \( \xi \) or \( u \) corresponds to the projection of \( \xi \) or \( u \) onto the basis function \( \varphi_j \), with the coefficients given by \( \xi_j = \langle \xi, \varphi_j \rangle \) or \( u_j = \langle u, \varphi_j \rangle \), respectively. We use the fast discrete sine transform to compute series \eqref{eq: forcing term 1} and to compute integrals against the test functions $\varphi_j$. We compute the \( L^2 \) and pointwise errors of the solution for both the kernel method and the PINN, and compute the $L^2$ error as a function of the number of measurements \( N \) for the kernel method. Our results are summarized in~\cref{tab: spatial results} and illustrated in~\cref{fig:1d_linear_elliptic}. We observe that both the kernel method and the neural network show significant improvement in the errors when trained with the negative Sobolev loss compared to the pointwise loss.

\begin{table}[h]
    \centering
    \begin{tabular}{|l|c|c|c|c|}
    \hline
    &Kernel ($H^{-1}$)& Kernel (pointwise) &  NN ($H^{-1}$) &NN (pointwise)  \\ 
         \hline 
        1D linear  & $3.95 \times 10^{-5}$ & 1.68& $0.0200 $&$1.49$ \\
        \hline
         2D semi-linear  & $0.0232$&  1.00 & $0.0203$ & 3.40  \\
         \hline 
    \end{tabular}
    \caption{Relative $L^2$ error of the estimated solutions for the spatial PDEs with $N = 4096$ measurements. The errors are reported for both the kernel and the NN approximation that are trained with either the $H^{-1}$ norm or the pointwise loss (in brackets).}
    \label{tab: spatial results}
\end{table}

\subsubsection{2D Semilinear Elliptic PDE} \label{sec: semilinear_PDE}
Next, we consider the semi-linear elliptic PDE with homogeneous Dirichlet boundary condition:
\begin{align}\label{eq: 2d semi-linear pde}
\begin{split}
 -\nu\Delta u + u + \sin(\pi u) &= \xi \quad  x \in \Omega = (0,1) \times (0,1), \\
    u &= 0 \quad x \in\partial \Omega.
\end{split}
\end{align}
We take $\nu = 10^{-1}$.
We consider a manufactured solution that is sampled from the series: 
\begin{equation}\label{eq: 2d semi-linear solution}
    u \sim \sum_{i=1}^\infty\sum_{j=1}^\infty \frac{u_{ij}}{(i^2 + j^2)^{1 + \varepsilon}}2\sin(i\pi x)\sin(j\pi y) ,\quad u_{ij} \sim \mathcal{N}(0,1) \; \text{i.i.d.}
\end{equation}
for some $\varepsilon >0$. In this case, $u^* \in H^{t}_0$ for $t < 1+ 2\varepsilon $ a.s.\thinspace which implies that $\xi \in H^{-s}$ a.s.\thinspace for $s > 1 - 2\varepsilon$. 

\begin{figure}[htp]%
    \centering
    \begin{subfigure}[b]{0.31\columnwidth}
        \centering
        \includegraphics[width=\textwidth]{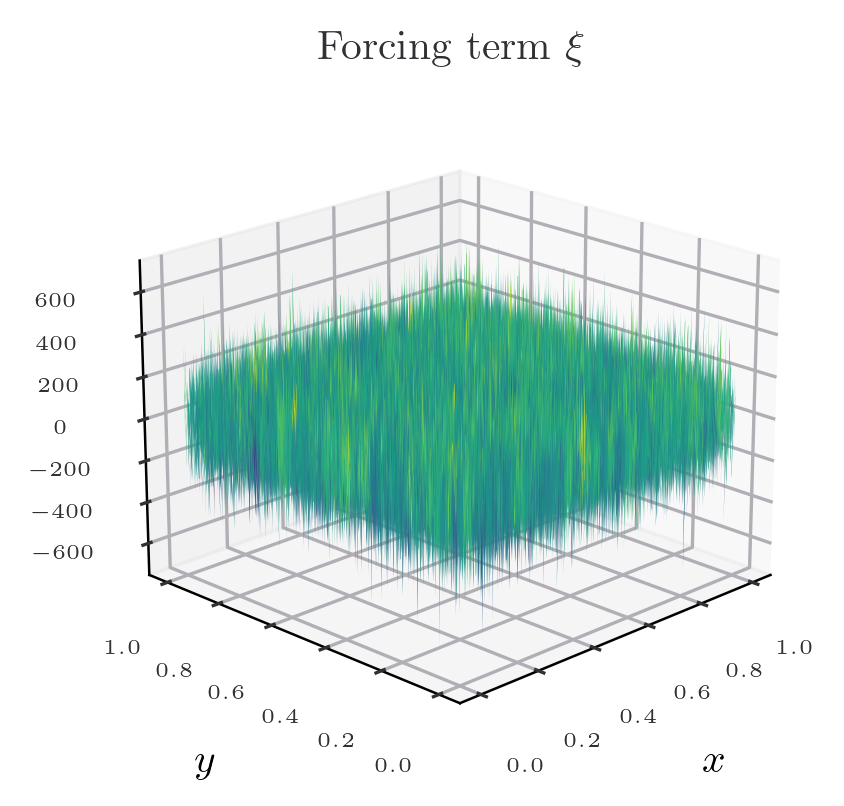}
        \caption{\footnotesize Forcing term}
    \end{subfigure}%
    \hspace{0.03\columnwidth}
    \begin{subfigure}[b]{0.31\columnwidth}
        \centering        \includegraphics[width=\textwidth]{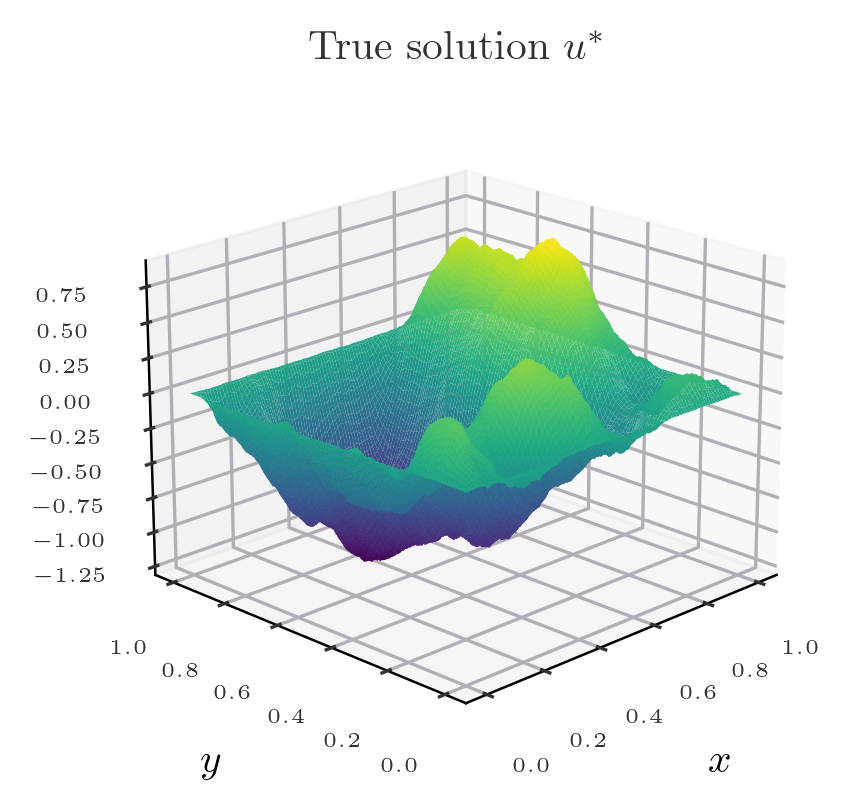}
        \caption{\footnotesize True solution}
    \end{subfigure}%

    \begin{subfigure}[b]{0.31\columnwidth}
        \centering
\includegraphics[width=\textwidth]{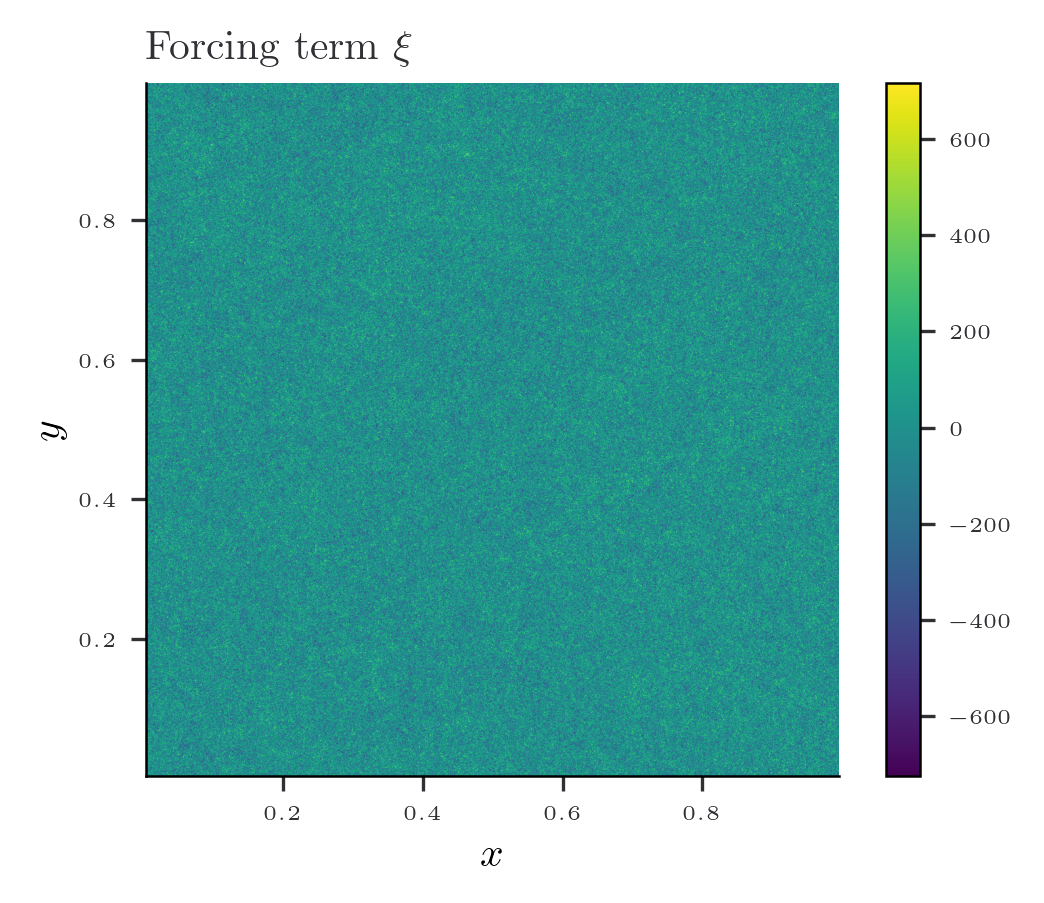}
        \caption{\footnotesize Forcing term }
    \end{subfigure}%
    \hspace{0.03\columnwidth}
    \begin{subfigure}[b]{0.31\columnwidth}
        \centering
        \includegraphics[width=\textwidth]{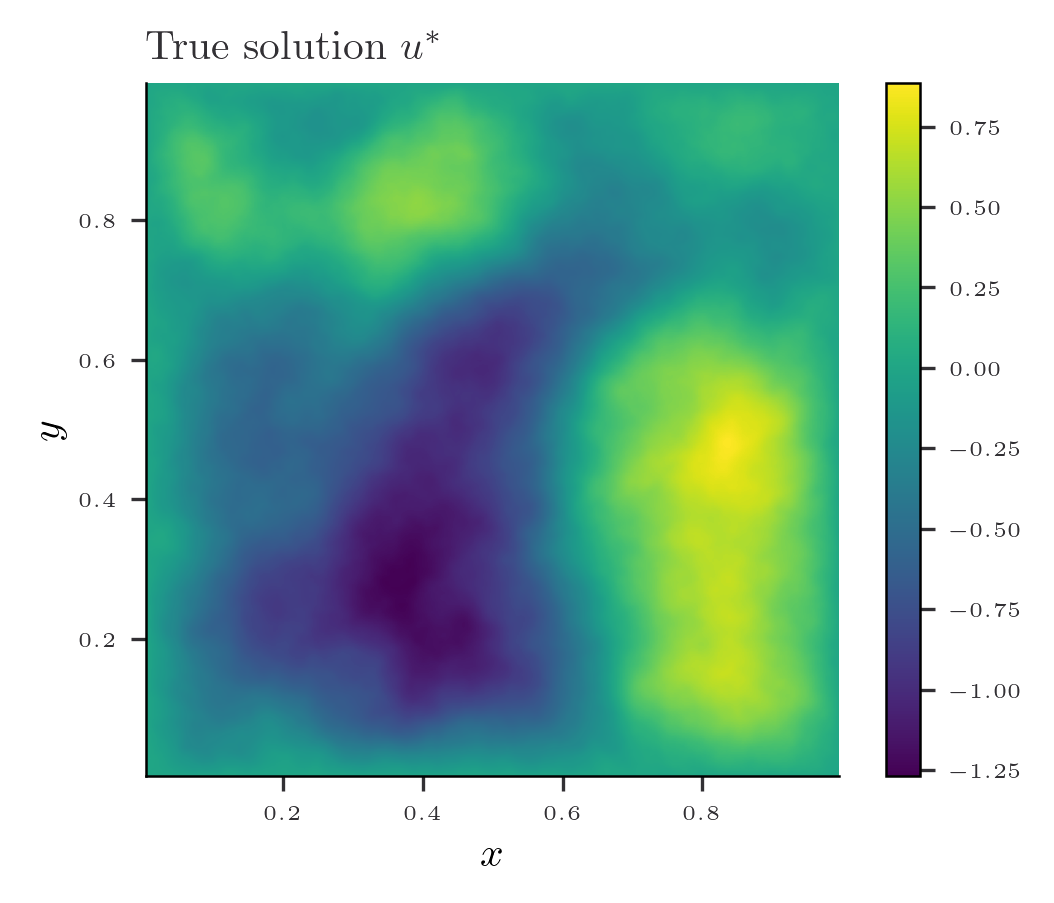}
        \caption{\footnotesize True solution}
    \end{subfigure}%
    \caption{One example solution to the two-dimensional semi-linear elliptic PDE \eqref{eq: 2d semi-linear pde}.}
    \label{fig: 2d semi-linear elliptic}%
\end{figure}

Again, we consider a truncated series expansion in~\eqref{eq: 2d semi-linear solution} with $1 \leq i,j \leq L$ and set $L= 2^{10}$ to be much larger than the maximum number of measurements. We generate each realization of the manufactured solution using the discrete sine transform. A sample of the forcing term and the corresponding solution are illustrated in \cref{fig: 2d semi-linear elliptic}.

\begin{figure}[htp]
    \centering
    \begin{subfigure}[b]{\textwidth}
        \centering
        {{\includegraphics[width=0.31\columnwidth]{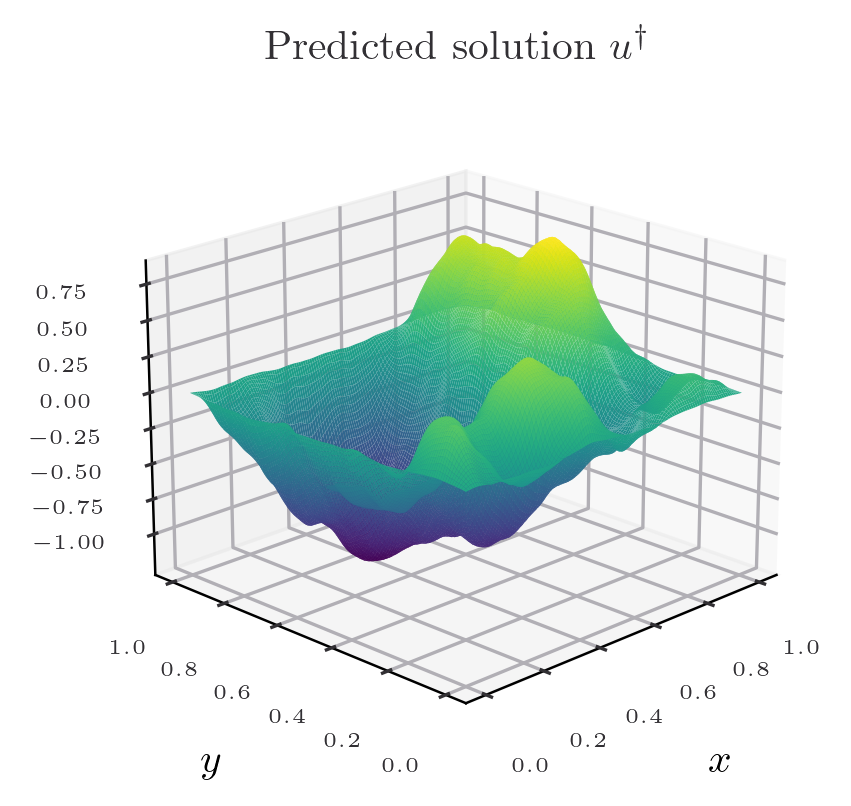} }}%
       {{\includegraphics[width=0.31\columnwidth]{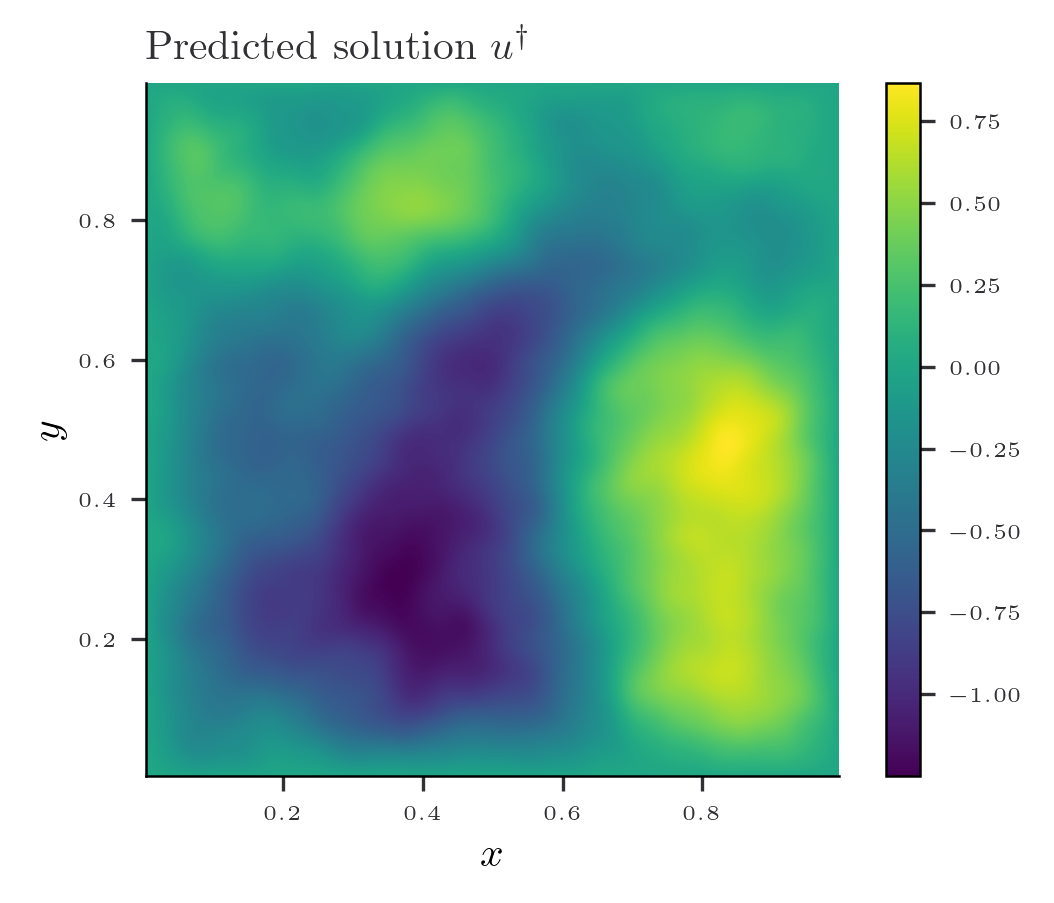} }}%
       {{\includegraphics[width=0.31\columnwidth]{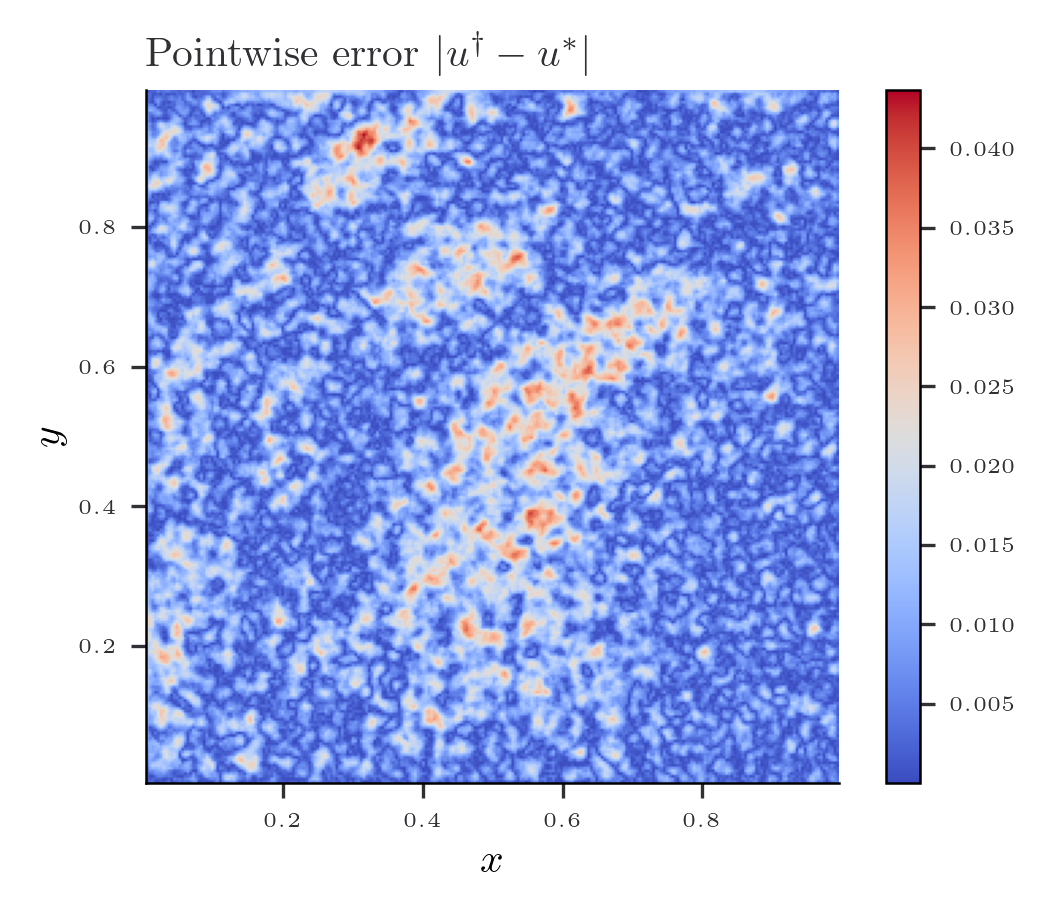} }}%
        \caption{Kernel: $H^{-1}$ loss}
    \end{subfigure}

        \begin{subfigure}[b]{\textwidth}
        \centering
        {{\includegraphics[width=0.31\columnwidth]{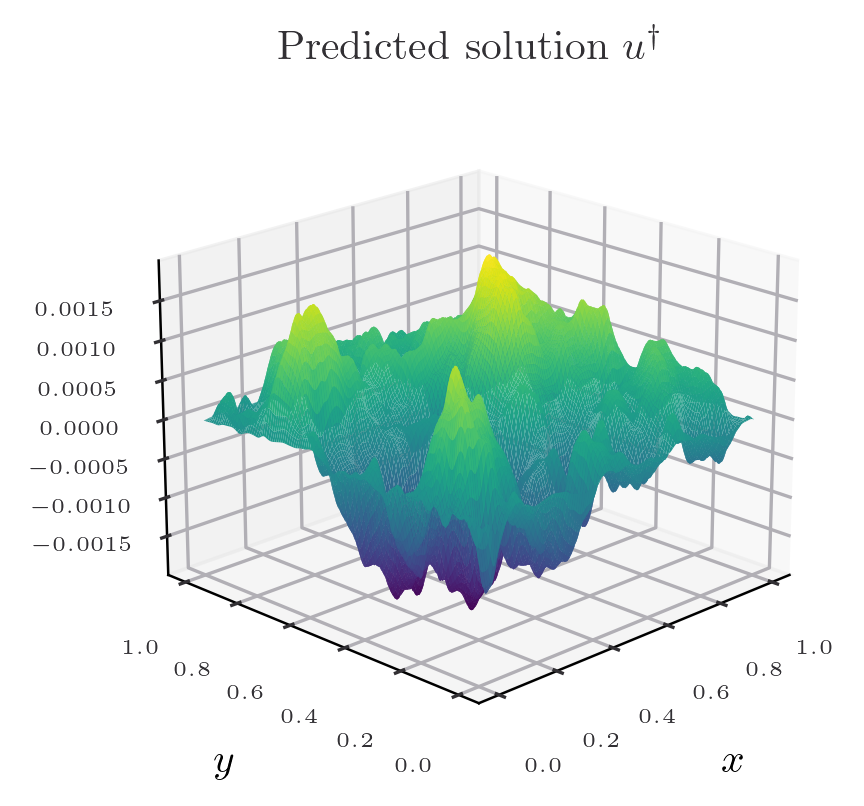} }}%
       {{\includegraphics[width=0.31\columnwidth]{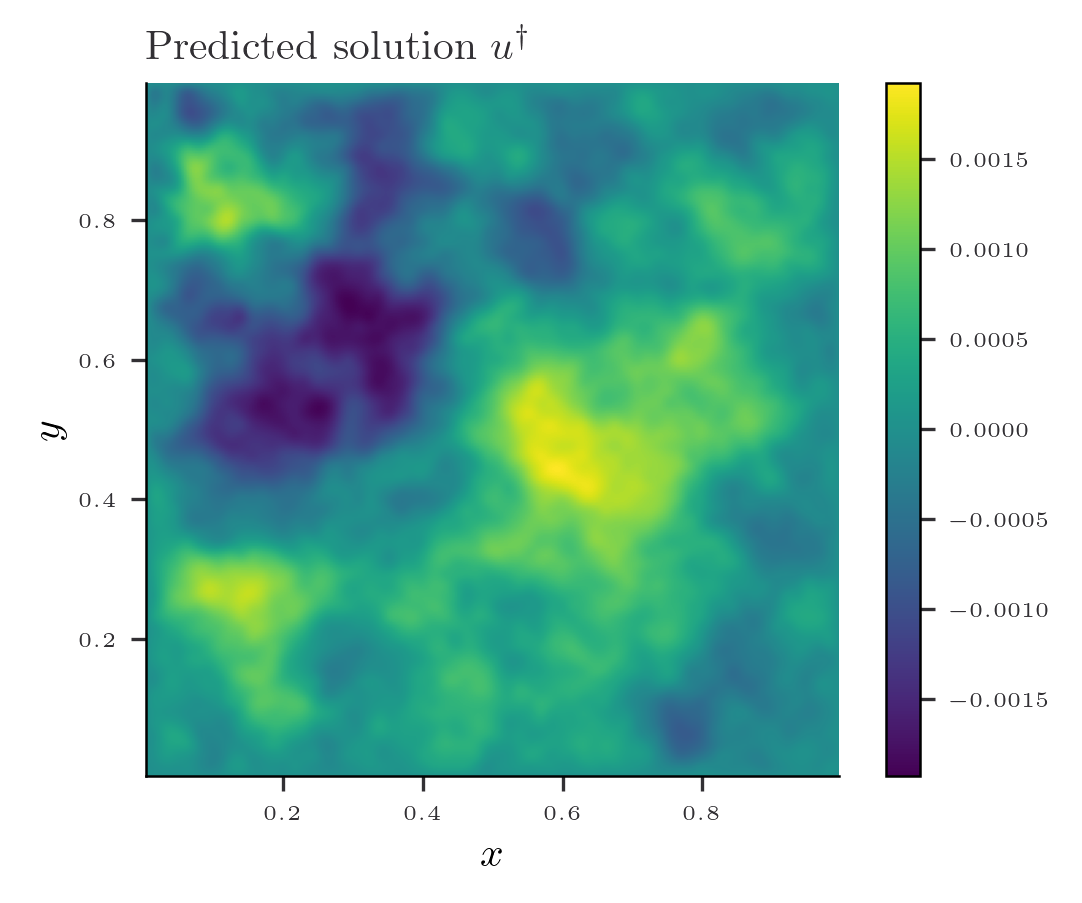} }}%
       {{\includegraphics[width=0.31\columnwidth]{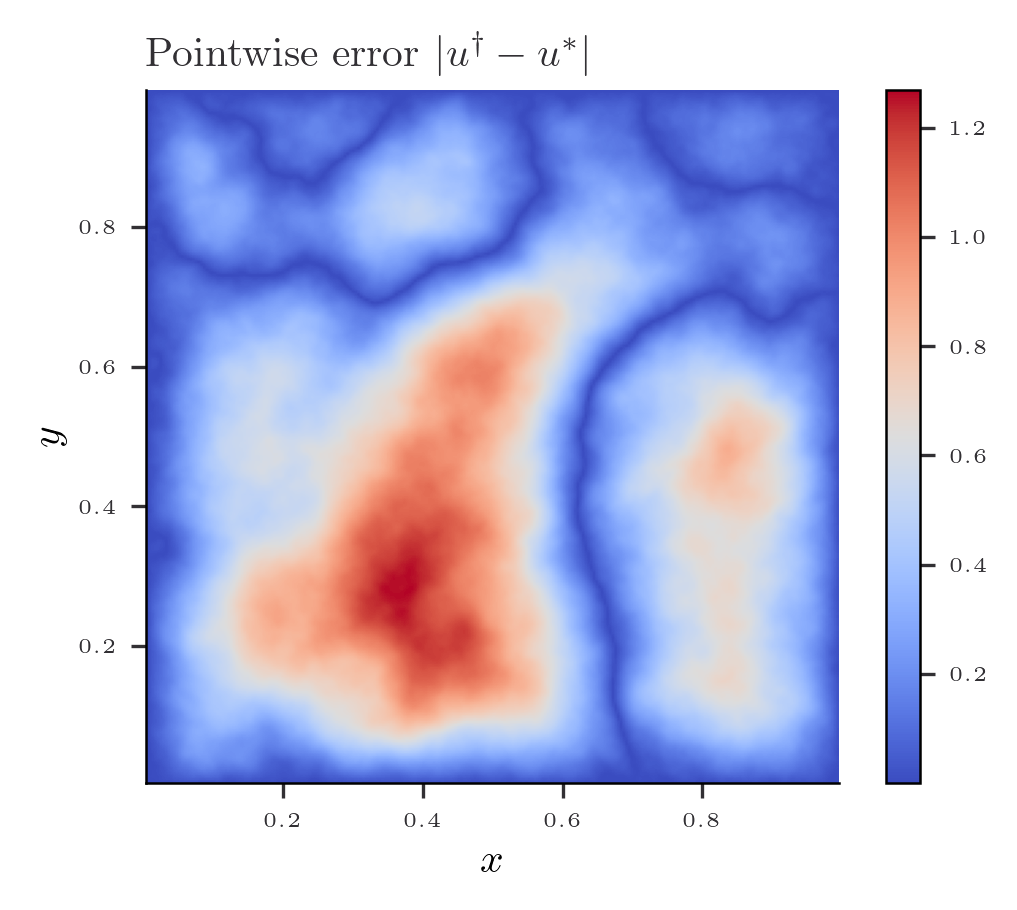} }}%
        \caption{Kernel: pointwise loss}
    \end{subfigure}

    \begin{subfigure}[b]{\textwidth}
        \centering
       {{\includegraphics[width=0.31\columnwidth]{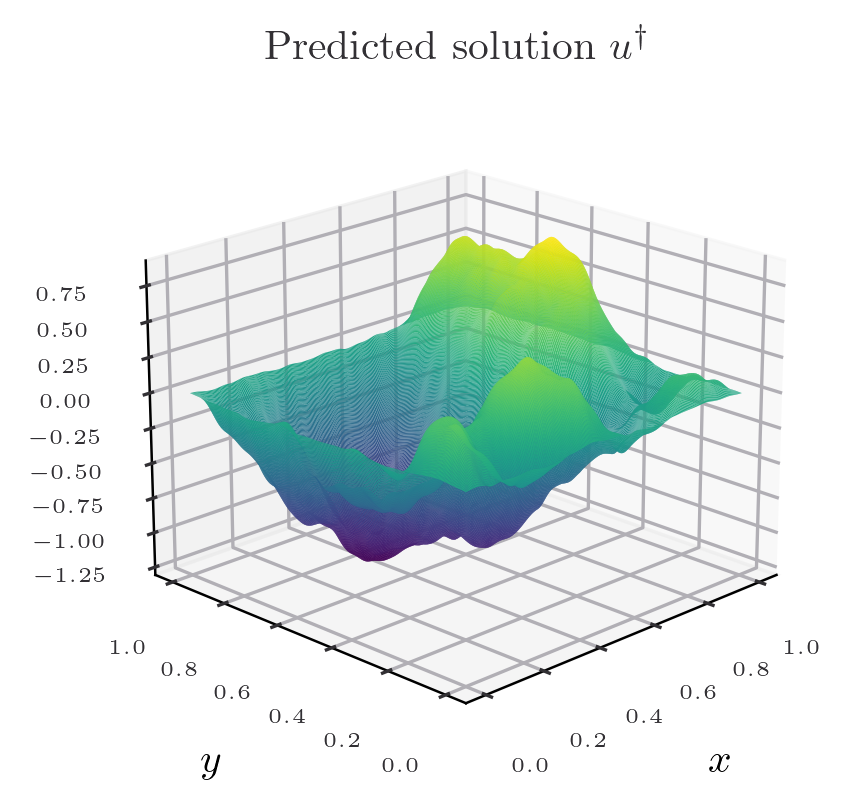} }}%
{{\includegraphics[width=0.31\columnwidth]{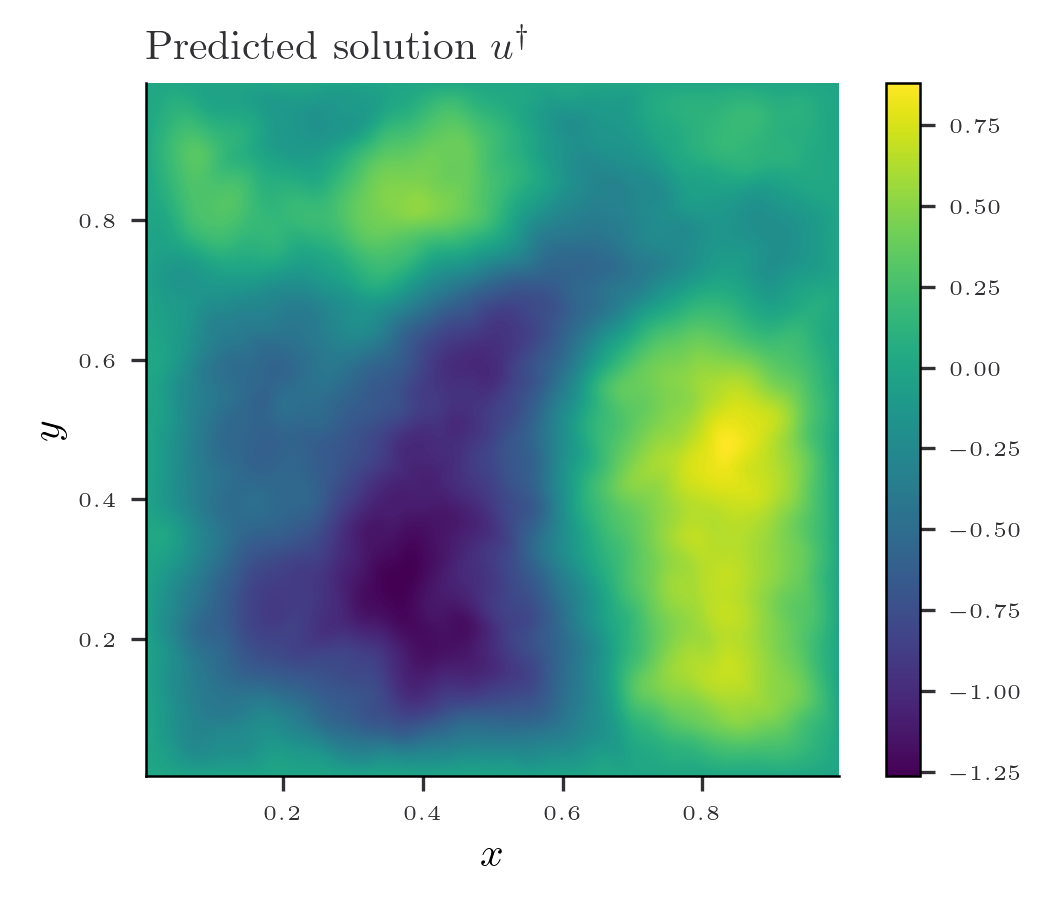} }}%
       {{\includegraphics[width=0.31\columnwidth]{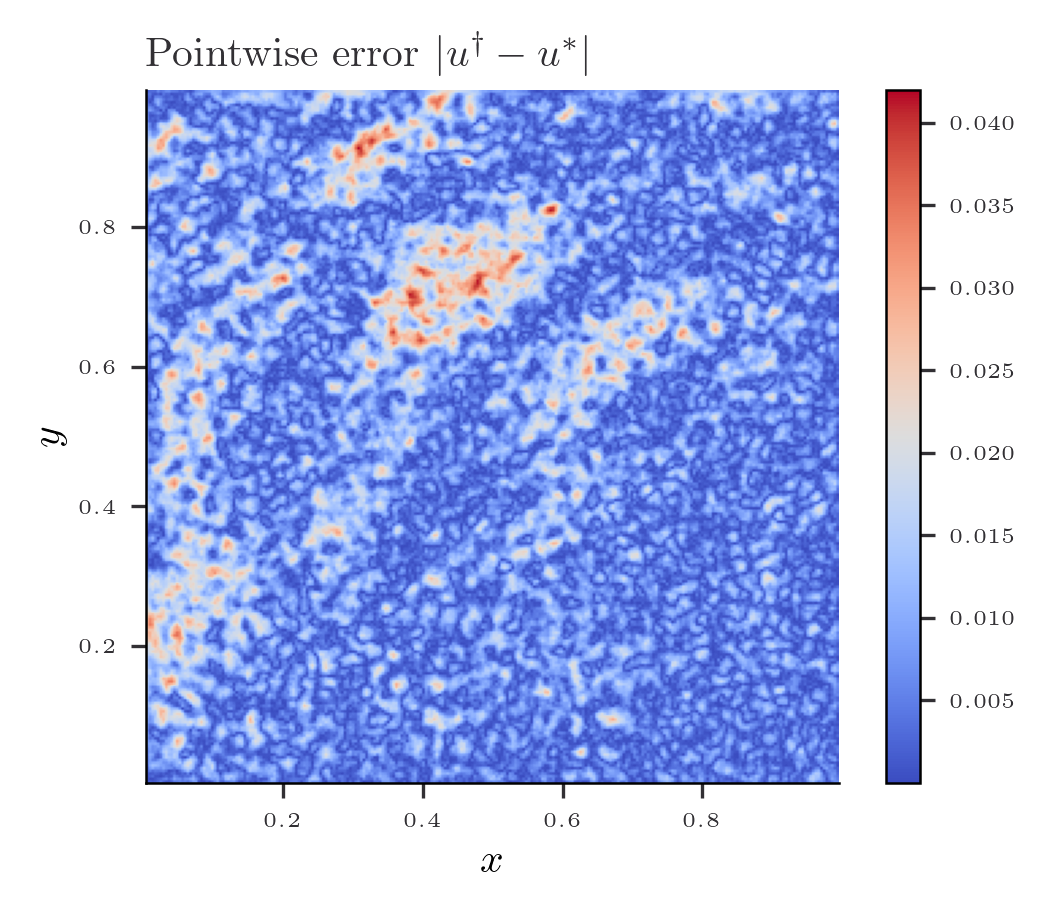} }}%
        \caption{PINN: $H^{-1}$ loss}
    \end{subfigure}

    \begin{subfigure}[b]{\textwidth}
        \centering
       {{\includegraphics[width=0.31\columnwidth]{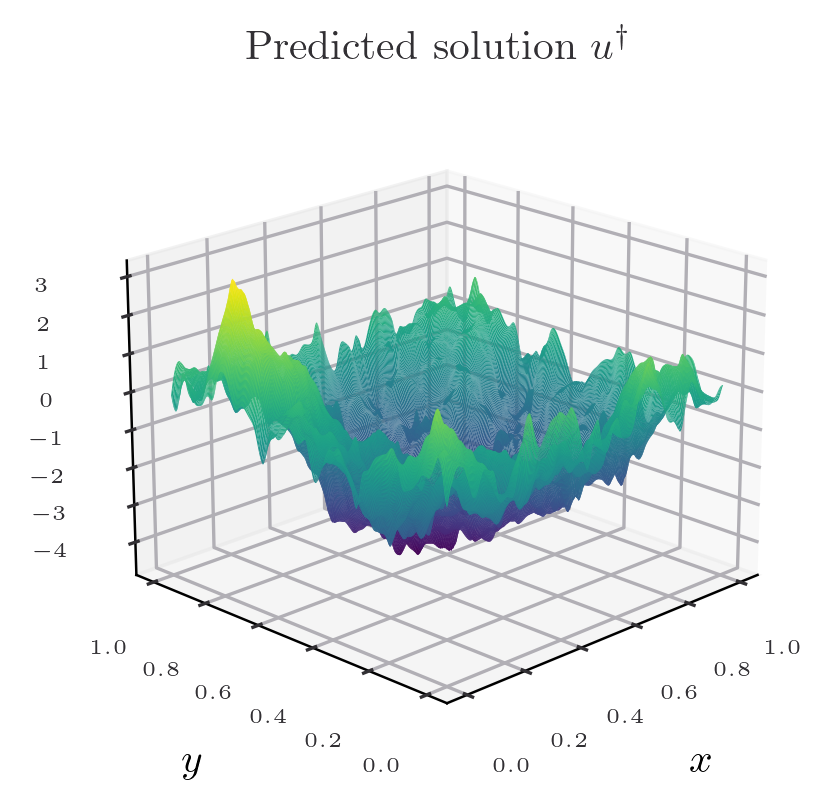} }}%
       {{\includegraphics[width=0.31\columnwidth]{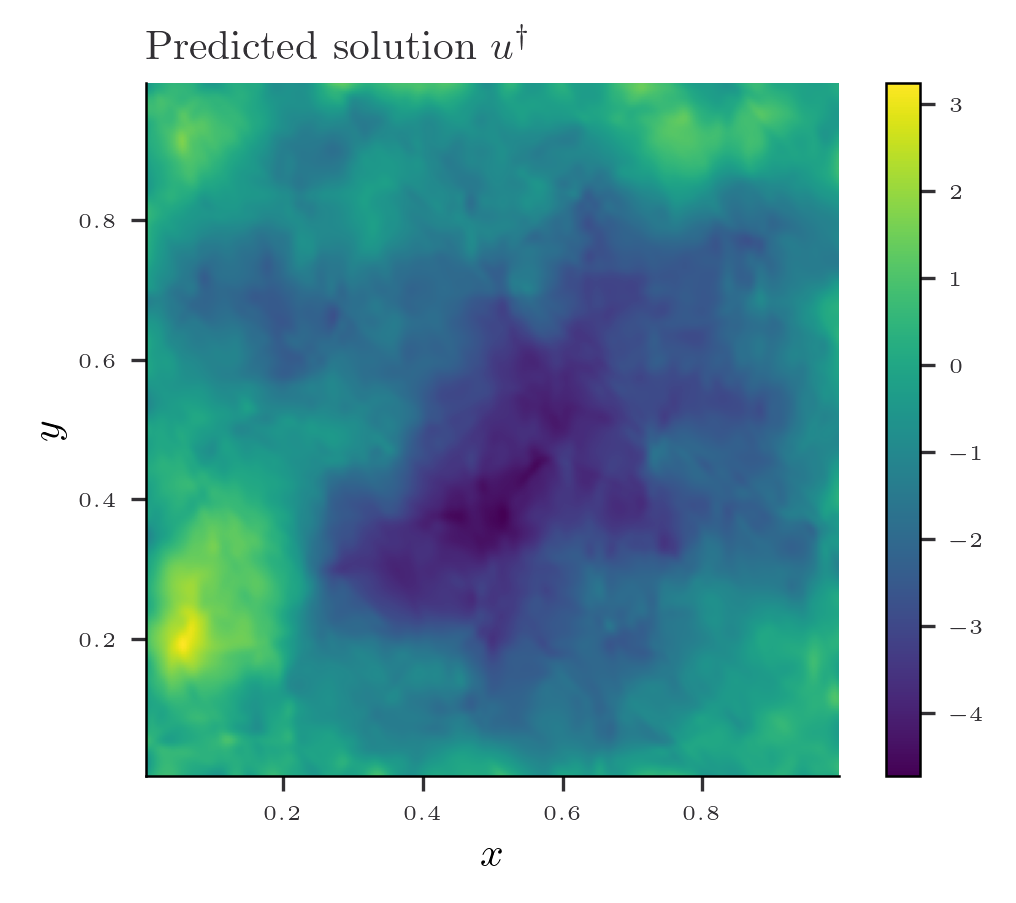} }}%
      {{\includegraphics[width=0.31\columnwidth]{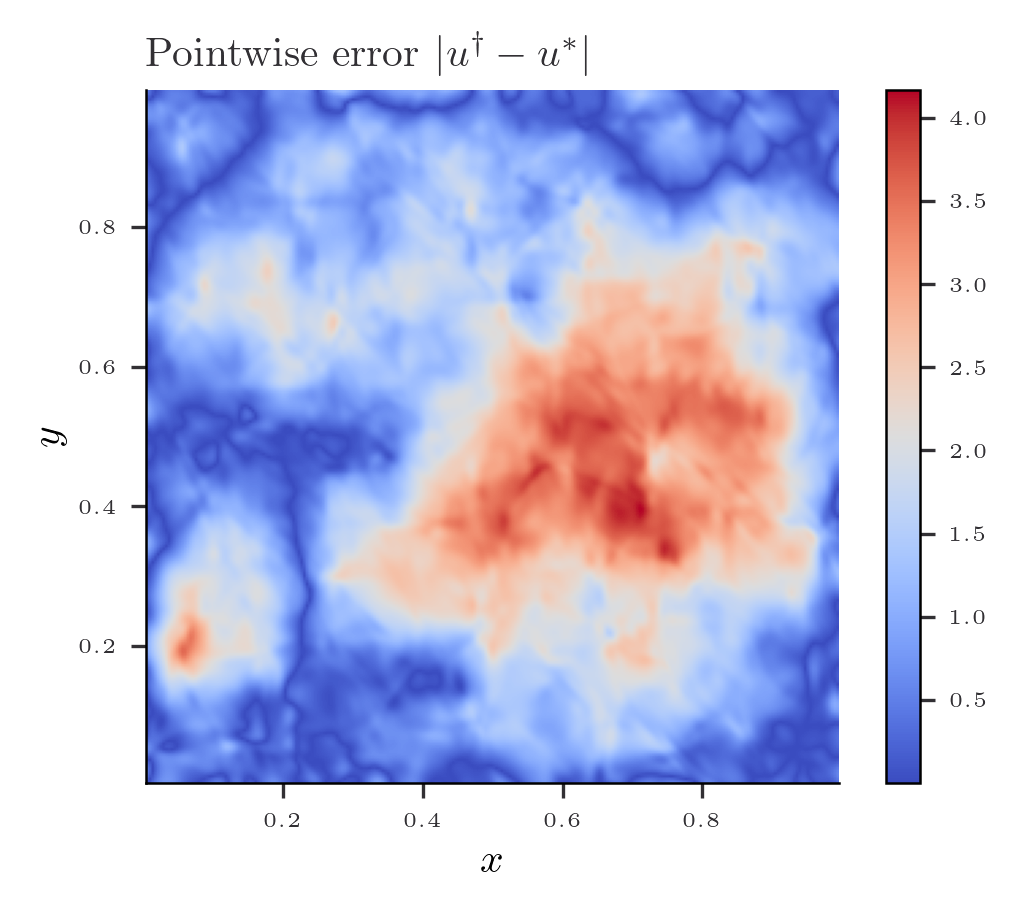} }}%
        \caption{PINN: pointwise loss}
    \end{subfigure}

    \caption{Numerical solutions to the two-dimensional semi-linear elliptic PDE \eqref{eq: 2d semi-linear pde} using the kernel and NN approximations with either the negative Sobolev loss or the pointwise loss.}
    \label{fig: 2d semi-linear elliptic solutions}
\end{figure}

We set \( \varepsilon = 0.15 \) and, as in the previous section, use the \( H^{-1} \) norm as the objective norm. The measurements are chosen as projections onto the basis functions \( \varphi_{ij}(x, y) = \sin(\pi i x)\sin(\pi j y) \). We train a kernel method and a NeS-PINN using $1 \leq i,j \leq N$ with $N =64$ for a total of 4096 measurements. We also train a PINN using the standard pointwise loss on the same grid to assess the impact of using the weak norm compared to the classical pointwise loss. The numerical results are summarized in \cref{tab: spatial results} and illustrated in \cref{fig: 2d semi-linear elliptic solutions}. We again observe significant improvement for both methods when using the negative Sobolev loss as opposed to the pointwise loss.


\subsubsection{The Impact of Norm Selection}
We now explore the effect of the selection of the $H^{-s}$ norm, for different choices of $s$, 
on the accuracy of the recovered solutions. Our approach differs from the standard physics-informed loss in two ways: it employs measurements as integrals against test functions rather than pointwise evaluations, and it formulates the loss as a negative Sobolev norm based on these measurements. The previous examples demonstrate that using weak measurements yields significantly better performance as compared to pointwise evaluations of the PDE residual.

\begin{figure}[htb]%
    \centering
    \begin{subfigure}[b]{0.31\columnwidth}
        \centering
        \includegraphics[width=\textwidth]{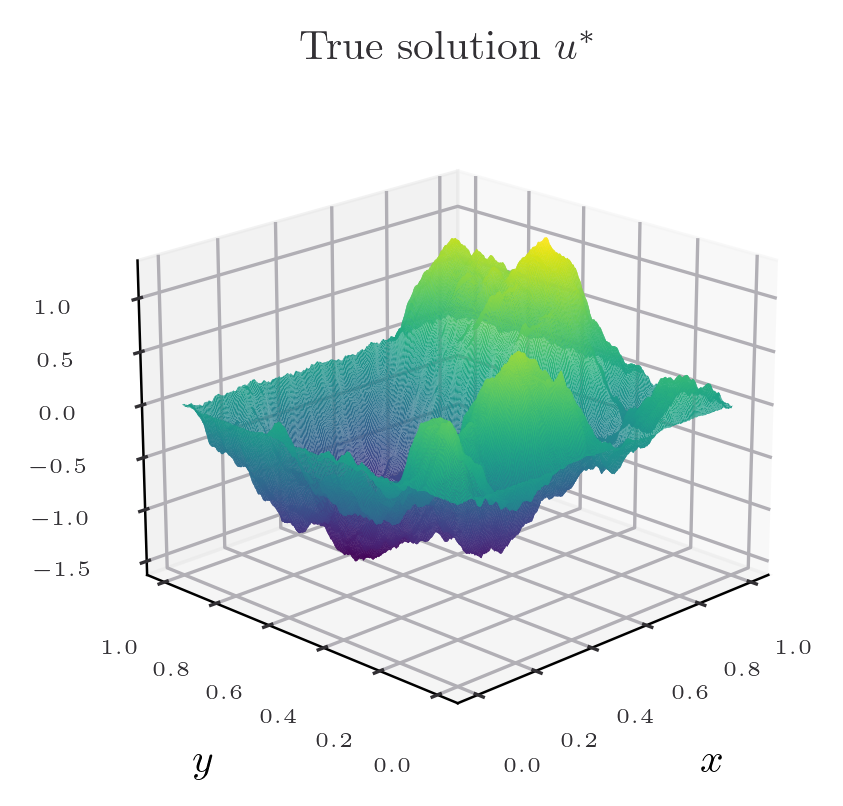}
        \caption{\footnotesize True solution}
    \end{subfigure}%
    \hspace{0.03\columnwidth}
    \begin{subfigure}[b]{0.31\columnwidth}
        \centering
        \includegraphics[width=\textwidth]{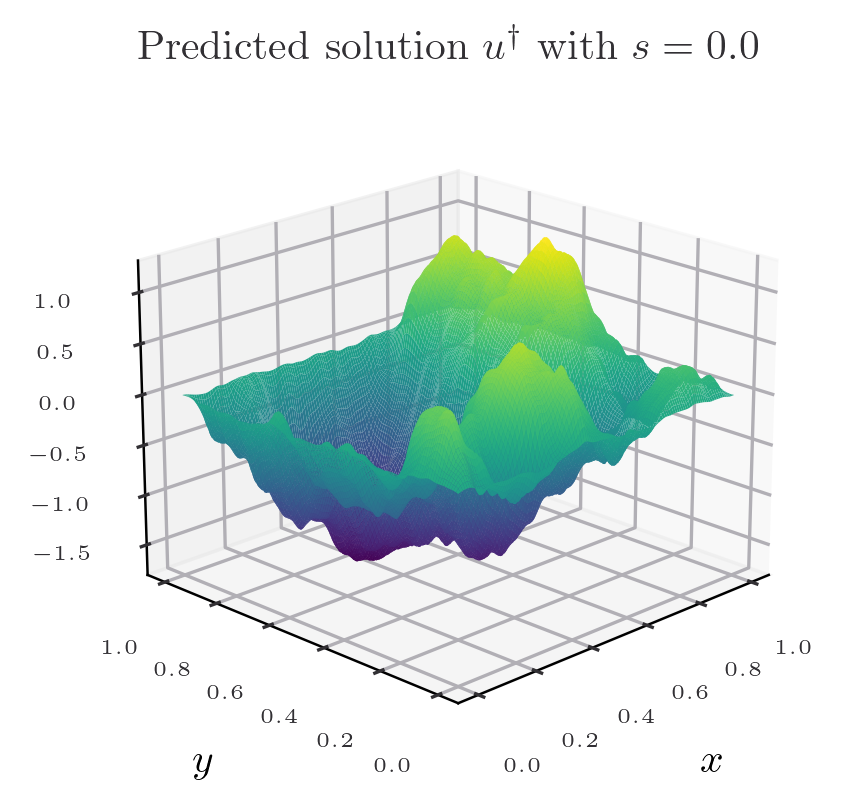}
        \caption{\footnotesize Unstable ($s=0$)}
    \end{subfigure}%
    
    \begin{subfigure}[b]{0.31\columnwidth}
        \centering
\includegraphics[width=\textwidth]{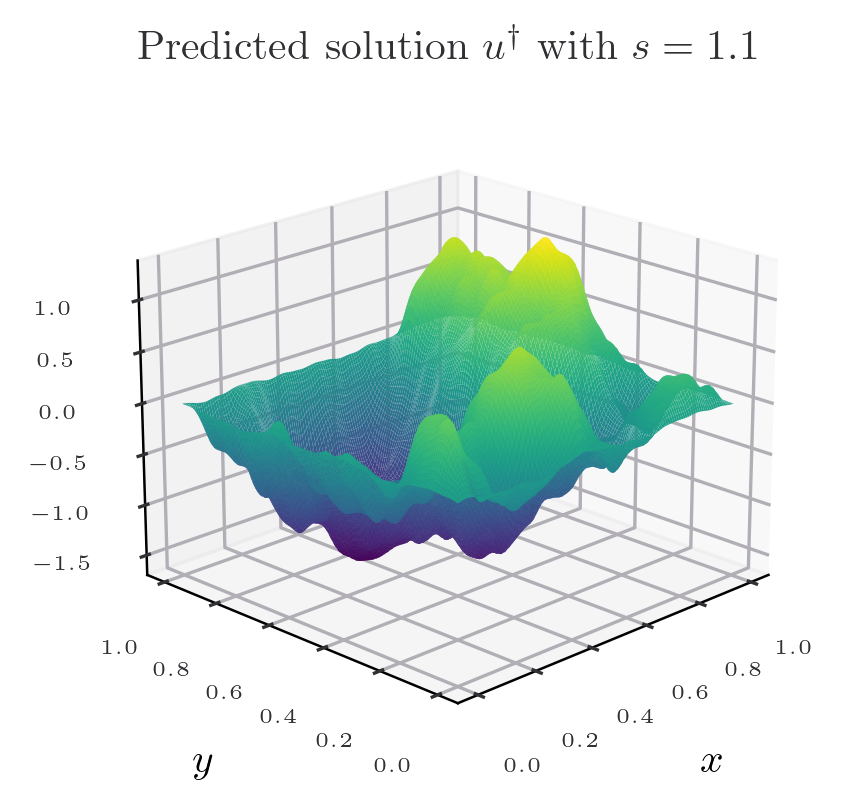}
        \caption{\footnotesize Optimal ($s=1.1$)}
    \end{subfigure}%
    \hspace{0.03\columnwidth}
    \begin{subfigure}[b]{0.31\columnwidth}
        \centering\includegraphics[width=\textwidth]{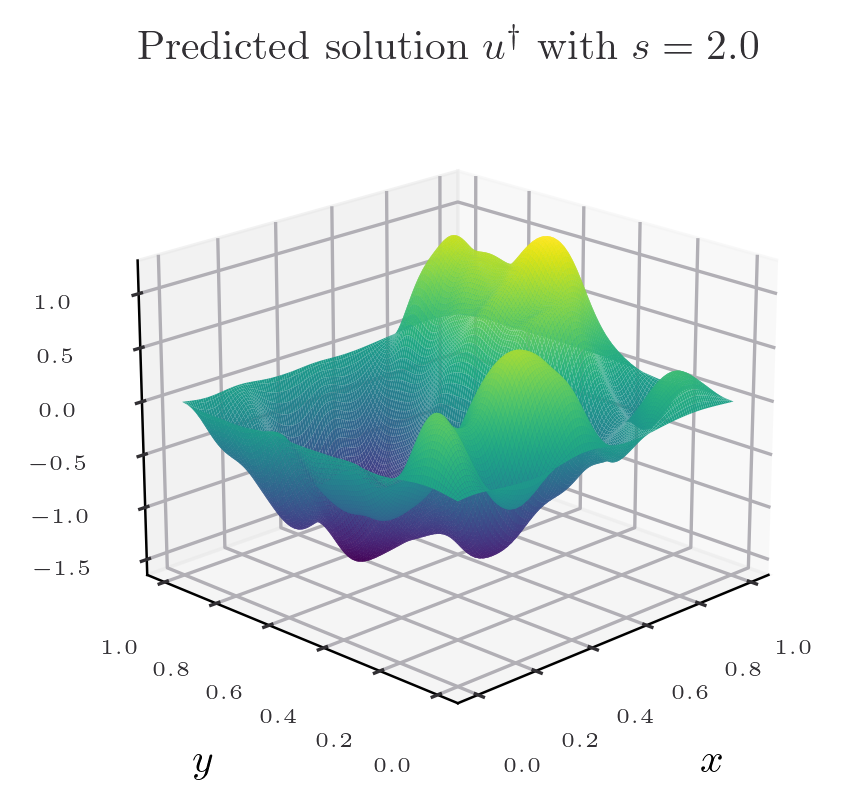}
        \caption{\footnotesize Oversmoothed ($s=2$)}
    \end{subfigure}%
    
    \caption{Impact of norm selection on the recovered solution.}
    \label{fig:effect_of_the_norm}%
\end{figure}

To demonstrate the effect of the Sobolev norm choice, we consider the two-dimensional PDE in~\eqref{eq: 2d semi-linear pde} with \( \varepsilon = 0 \). In this case, the solution belongs to $H^{t}$ for  \( t < 1 \) and the forcing term \( \xi \) belongs to \( H^{-s} \) for \( s > 1 \) (note the strict inequality). We recover the true solution using several \( H^{-s} \) norms: \( s = 0\footnote{The case \( s=0 \) corresponds to \( H^0 = L^2\).}, 0.5, 1.0, 1.1, 2.0 \), which all computed with the same test functions. When the norm is misspecified relative to the true solution's regularity (i.e., \( s=0, 0.5, 1.0 \)), the neural network underperforms as compared to cases where the norm is correctly specified (i.e., \( s=1.1, 2 \)). Furthermore, we observe a notable difference in performance at the threshold between \( s=1 \) and \( s=1.1 \), even when truncating the expansion \eqref{eq: 2d semi-linear solution}. When using the \( L^2 \) norm, the neural network's training becomes unstable. The best results are achieved when the norm is well-specified with \( s=1.1 \). However, when using a weaker norm \( s=2.0 \), the error increases again due to over-smoothing effects on the higher frequencies. This effect can also be observed with the kernel method, although its performance is less sensitive to the choice of the norm than the neural network approximation. Based on the results, we deduce that both the use of test functions (as opposed to a pointwise loss) and the choice of loss are crucial to achieving high accuracy on PDEs with rough right-hand sides. These observations are summarized in \cref{tab: norm comparison results} and illustrated in \cref{fig:effect_of_the_norm}. 
Moreover, since the right-hand side $\xi$ belongs to the space $H^{-s}$ for any $s >1$, the choice of $s = 1.1$ is in theory not the only possible one. However, we observe that setting $s$ too close to $1$ (such as $s=1.01$) can also lead to a decrease in performance due to instabilities in training. These observations suggest the following strategy to construct machine learning approximations: at the beginning of training, one may begin the training process  with a weaker norm (i.e., large $s$) to promote an over-smoothed solution, and then gradually decrease $s$ throughout training until convergence has been achieved.

\begin{table}[h]
    \centering
    \begin{tabular}{|l|c|c|c|c|c|}
    \hline
    Norm $H^{-s}$& $s =0.0$& $s=0.5$ & $s=1.0$ & $s=1.1$& $s=2.0$  \\ 
         \hline 
       Kernel method & $0.121 $&$0.116$ & $0.046$& $0.041$&$0.079$\\
        \hline
        Neural network   &$0.514$ & $0.434$& $0.116 $ & $0.046$ & $0.082 $\\
         \hline 
    \end{tabular}
    \caption{Relative $L^2$ error for different choices of Sobolev norms as the loss function}
    \label{tab: norm comparison results}
\end{table}

\subsection{Time-dependent Stochastic PDEs}
\label{sec: time dependent spdes}
In this section, we address the problem of solving time-dependent semi-linear SPDEs of the form:
\begin{align}
\label{eq: time-dependent PDE}
\begin{split}
\partial_t u(t,x)  = A u(t,x) + f(u)(t,x) + \dot{\xi}(t,x) & \quad\text{in } \Omega, \\
        u = 0 & \quad\text{on } \partial \Omega \times [0, T], \\
        u = g &\quad\text{on } \Omega \times \{t = 0\},
\end{split}
\end{align}
where $A$ is a linear operator,  \( \dot{\xi} \) is a stochastic forcing term and $f$ is a (nonlinear) Nemytskii operator. This is a general form encountered in many applications \cite{lord2014introduction}.

In \cref{app: gaussian measures}, we define white noise, which serves as our rough forcing term due to its relevance in practical applications. However, our methodology is equally applicable to more regular forcing terms, such as colored noise, which also appear in areas like stochastic fluid dynamics. 

We begin by describing the semi-implicit Euler time discretization for \eqref{eq: time-dependent PDE} in \cref{sec: implicit euler}, and then apply our method to two specific cases: the stochastic heat equation in \cref{sec: heat equation} and the stochastic Allen-Cahn equation in \cref{sec: allen-cahn}. A detailed discussion of time-dependent SPDEs and numerical methods for solving them is beyond the scope of this section. For a more comprehensive treatment, we refer the reader to \cite{lord2014introduction}.

\subsubsection{Semi-implicit Euler Method}
\label{sec: implicit euler}
Given a sample path of \( \xi(t) \), the SPDE \cref{eq: time-dependent PDE} can be discretized in time using a semi-implicit Euler scheme. At the \( k \)-th time step, with the solution \( u_k \), the value of \( u_{k+1} \) is computed as follows (boundary conditions omitted):
\begin{align}
\label{eq:siem}
    u_{k+1}  &= u_k + \delta t f(u_k) + \delta t A u_{k+1} + \delta \xi_k,
\end{align}
where \( \delta t \) is the step-size in time and \( \delta \xi_k \) represents the increment of the stochastic forcing term. Rewriting \eqref{eq:siem}, we have:
\begin{align*}\label{eq: semi-implicit euler pde problem}
      \Big( I - \delta t A \Big) u_{k+1} &= u_k + \delta t f(u_k) + \delta \xi_k.
\end{align*}
Thus, given \( u_k \), the solution at the next time step \( u_{k+1} \) is obtained by solving a linear PDE with a rough right-hand side, which almost surely does not belong to \( L^2(\Omega) \). The linear PDE can be solved using the algorithm described in \cref{sec: numerical method}. Although the fully implicit Euler method could also be used, we do not consider it here.

\subsubsection{Stochastic Heat Equation}
\label{sec: heat equation}
We now consider the stochastic heat equation over the spatial domain $\Omega = (0,1)$ and time interval $[0,T]$, which is a form of infinite dimensional Ornstein–Uhlenbeck process \cite{lord2014introduction}:
\begin{align*}
\begin{split}
 \partial_t u &= \nu\Delta u + \sigma\dot{\xi}, \quad \text{in } \Omega, \\
        u &= 0, \quad \text{on } \partial \Omega \times [0, T], \\
        u &= g, \quad \text{on } \Omega \times \{t = 0\},
\end{split}
\end{align*}
where \( \dot{\xi} \) represents space-time white noise, and parameters \( \nu = 0.025 \) and \( \sigma = 0.1 \). We first compute a reference solution \( u^* \), using a spectral Galerkin method on a fine mesh, with a time resolution defined by step-size \( \Delta t = 2^{-13} \) and spatial resolution  \( L = 2^{11} \), up to the final time \( T = 1.0 \). This fine mesh solution is used as the `ground truth' to compute errors and error rates of the kernel solution. 

\begin{figure}[h]
\begin{subfigure}[b]{\textwidth}
\centering
{{\includegraphics[width=0.31\columnwidth]{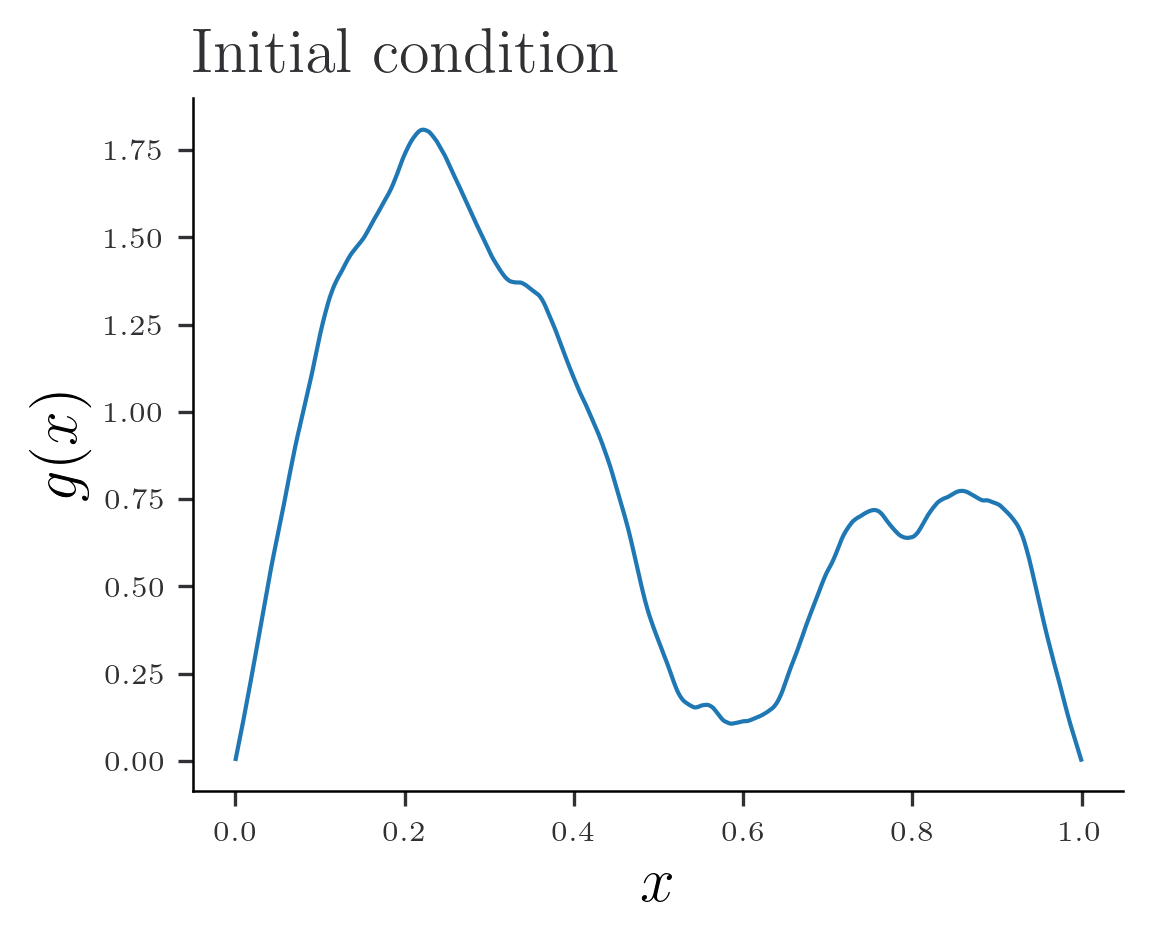} }}%
{{\includegraphics[width=0.31\columnwidth]{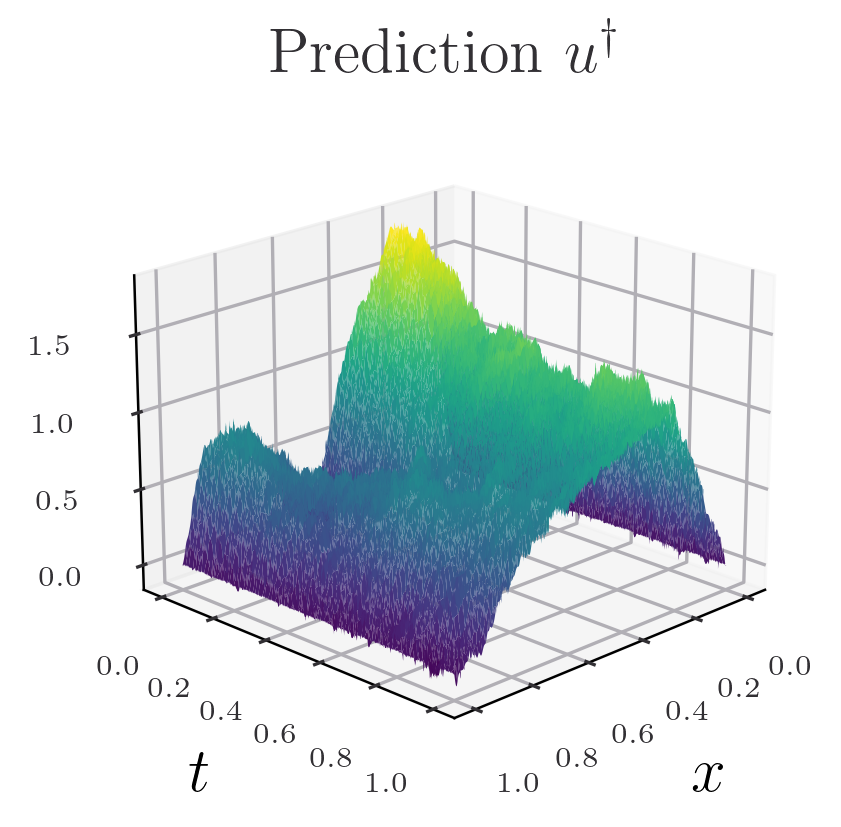} }}%
{{\includegraphics[width=0.31\columnwidth]{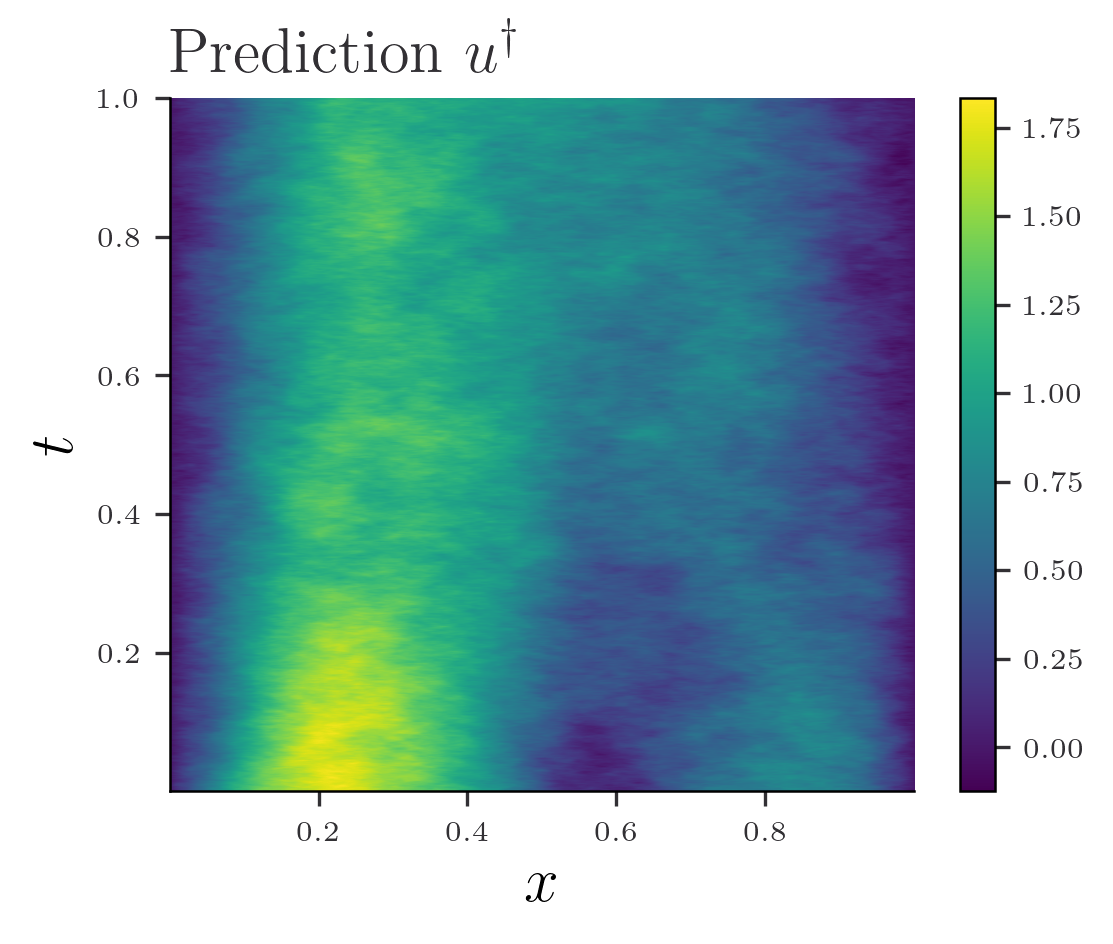} }}%
\caption{Kernel Solution computed with $N= 250$ measurements and $\Delta t =2^{-11}$.}
\end{subfigure}

\begin{subfigure}[b]{\textwidth}{{\includegraphics[width=0.31\columnwidth]{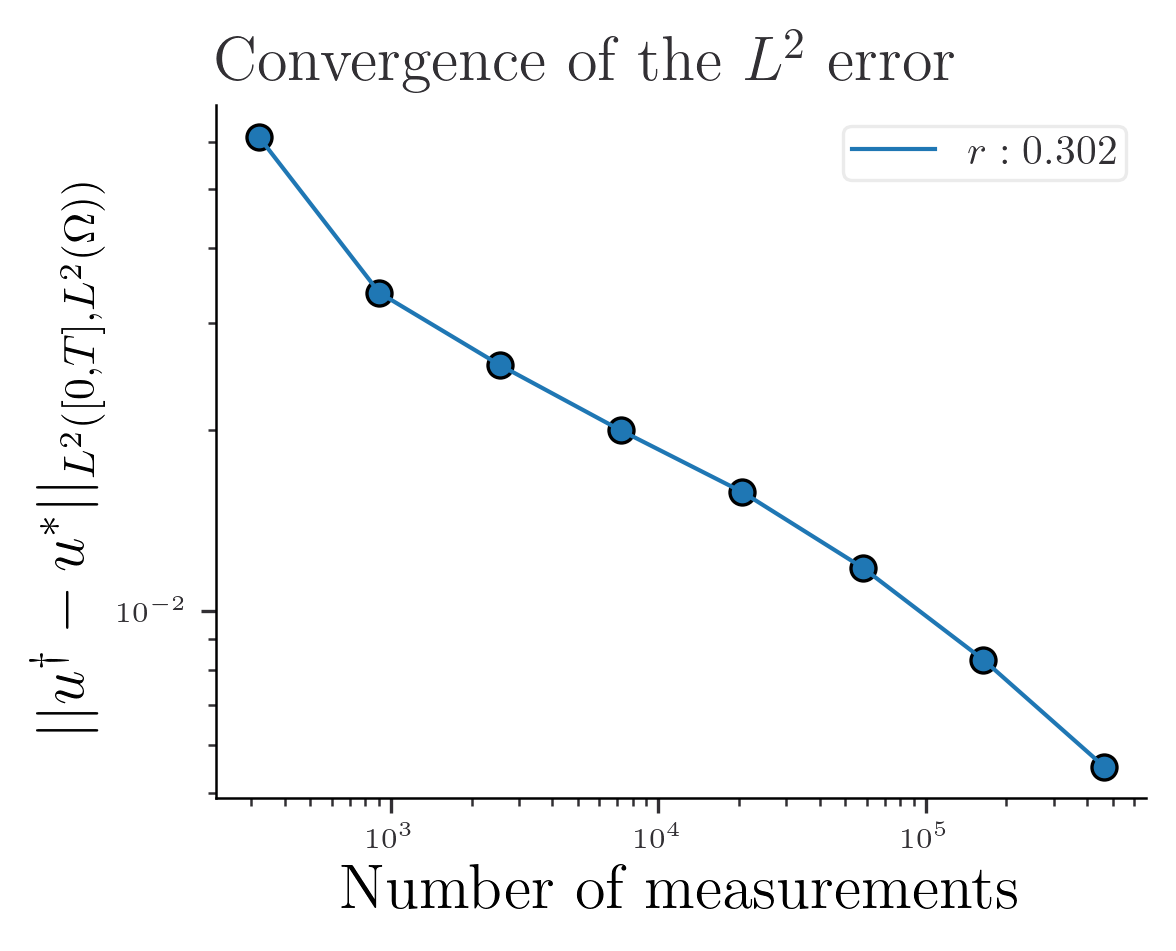} }}%
{{\includegraphics[width=0.31\columnwidth]{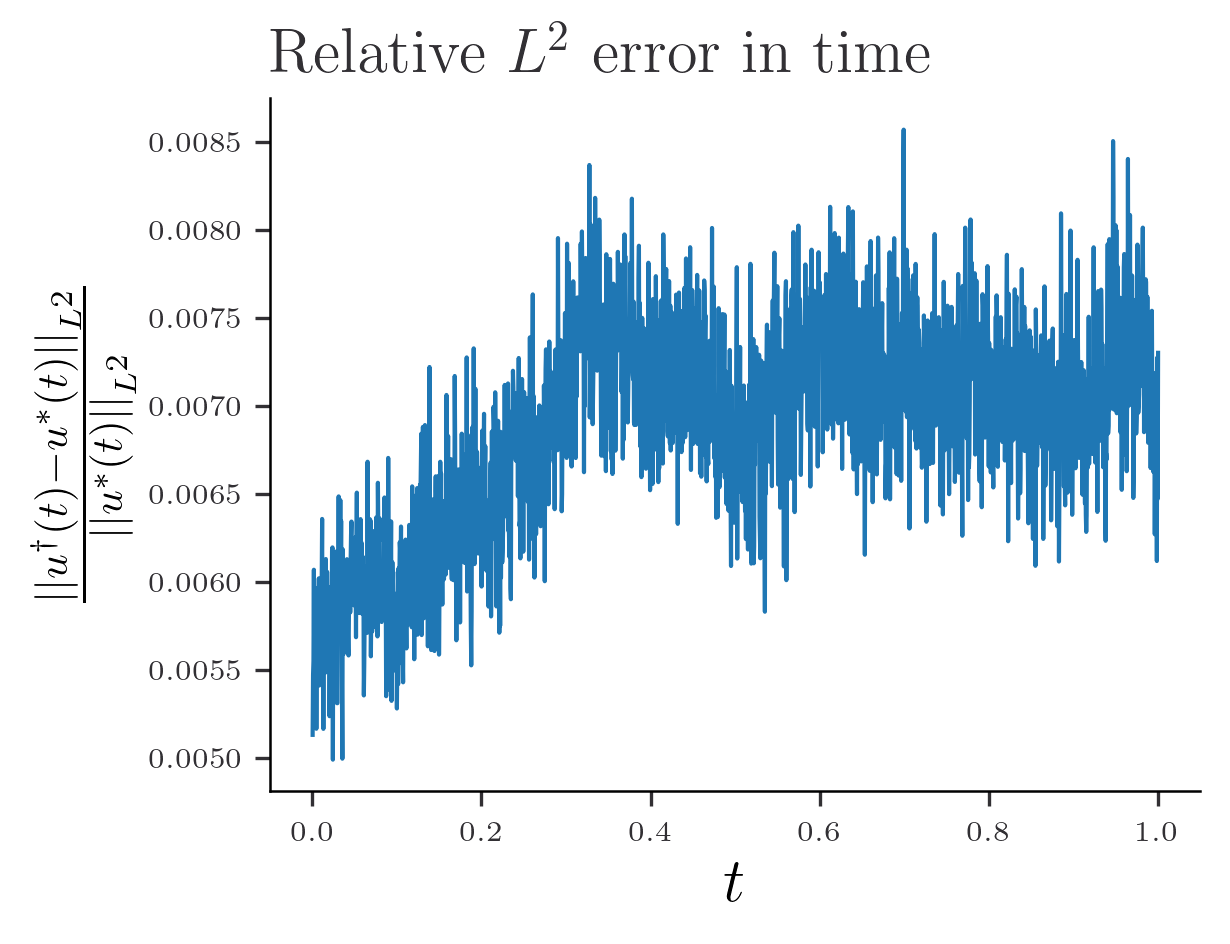} }}%
{{\includegraphics[width=0.31\columnwidth]{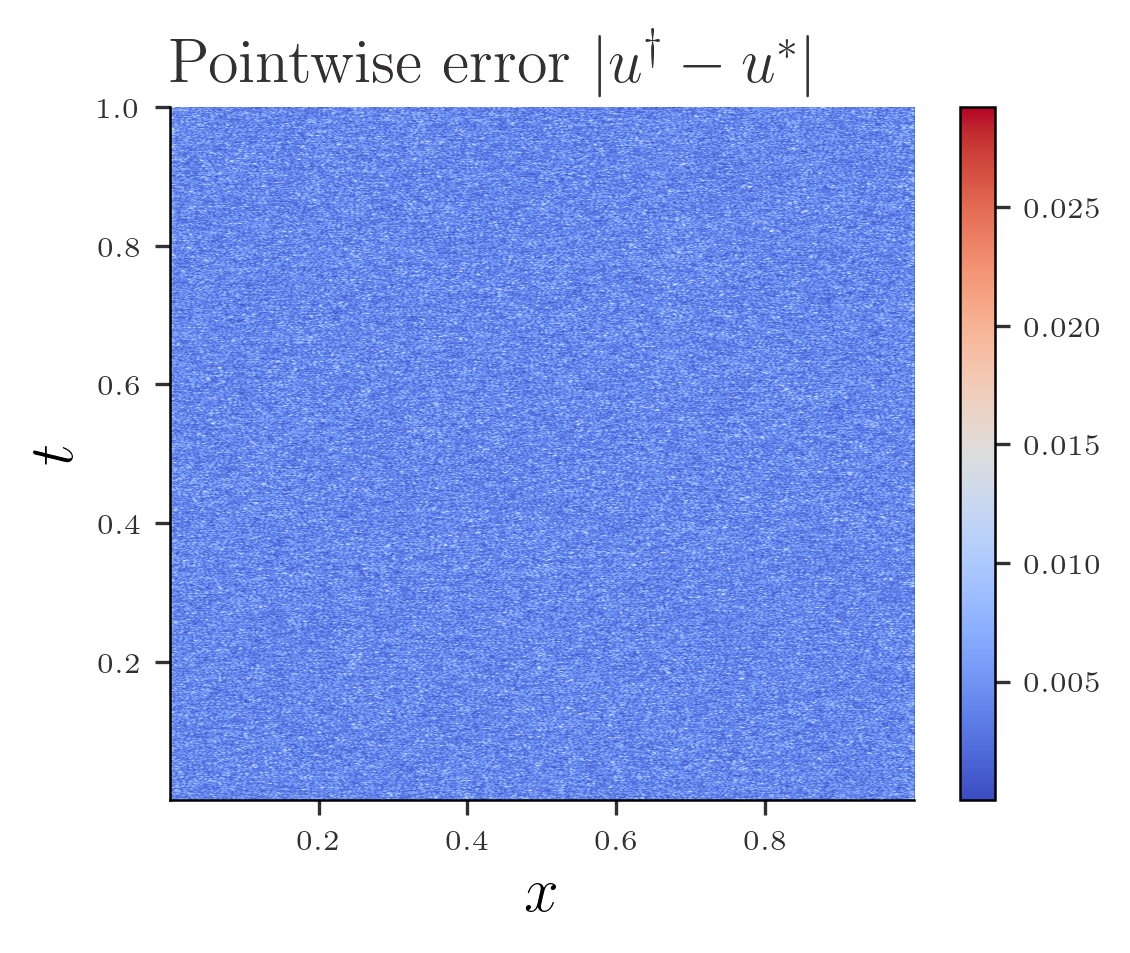} }}%
\caption{Convergence of the $L^2$ error and the pointwise error for the kernel-based prediction as compared to the fine mesh solution.}
\end{subfigure}
    \caption{Numerical solutions for the stochastic heat equation.}
    \label{fig: stochastic heat}
\end{figure}

Next, we compute the kernel-based prediction by constructing the measurement space \( \Phi^N \) using piecewise linear finite element bases, where \( N \) represents the number of basis functions. To evaluate the error, we compute the \( L^2([0, T];L^2(0,1)) \) error of the solution in time and space, which is defined as:
\begin{align*}
    ||u^\dagger - u^*||^2_{L^2([0, T];L^2(0,1))} = \int_{0}^T ||u^\dagger(t) - u^*(t)||^2_{L^2(0,1)} \, dt,
\end{align*}
where \( u^* \) is the fine mesh solution. We employ the semi-implicit integrator in~\eqref{eq:siem} to evolve the solution in time with the CFL condition \( N^2 \Delta t = 5 \). The finest discretization for the kernel solution is set to \( \Delta t = 2^{-11} \) and \( N = 250 \), which is coarser than the discretization used for the fine mesh solution $u^*$. 
The corresponding results and errors are illustrated in~\cref{fig: stochastic heat}, and the results for the finest discretization are recorded in \cref{tab: time results}. The results demonstrate that our method offers an accurate and stable approximation of the solution to the stochastic heat equation.

\begin{table}[h]
    \centering
    \begin{tabular}{|l|c|}
    \hline
    &Kernel  method \\ 
         \hline 
        Stochastic heat equation  & $6.84 \times 10^{-3}$ \\
        \hline
         Stochastic Allen-Cahn equation  & $2.07 \times 10^{-2}$ \\
         \hline 
    \end{tabular}
    \caption{Relative $L^2([0, T];L^2(0,1))$ error of the kernel-based prediction for the solutions to the time-dependent PDEs.}
    \label{tab: time results}
\end{table}

\subsubsection{Stochastic Allen-Cahn Equation}
\label{sec: allen-cahn}
Finally, we consider the stochastic Allen-Cahn equation on the spatial domain $\Omega=(0,1)$ and time interval $[0,T]$:
\begin{align*}
\begin{split}
\partial_t u &= \nu\Delta u + u - u^3 + \sigma\dot{\xi}, \quad \text{in } \Omega, \\
    u &= 0, \quad \text{on } \partial \Omega \times [0, T], \\
    u &= g, \quad \text{on } \Omega \times \{t = 0\},
\end{split}
\end{align*}
where \( \xi \) is space-time white noise, \( \nu = 10^{-4} \) and \(\sigma = 0.01\). As in the previous case, we first compute a reference solution \( u^* \) on a fine mesh using a spectral Galerkin method, integrated with a semi-implicit Euler scheme. The time step is \( \Delta t = 2^{-12} \), and the spatial resolution is \( L = 2^{11} \), with computations carried out up to the final time \( T = 10.0 \). This fine mesh solution is used as the `ground truth' to compute errors and error rates of the kernel solution. 

\begin{figure}[h]
\begin{subfigure}[b]{\textwidth}
\centering
{{\includegraphics[width=0.31\columnwidth]{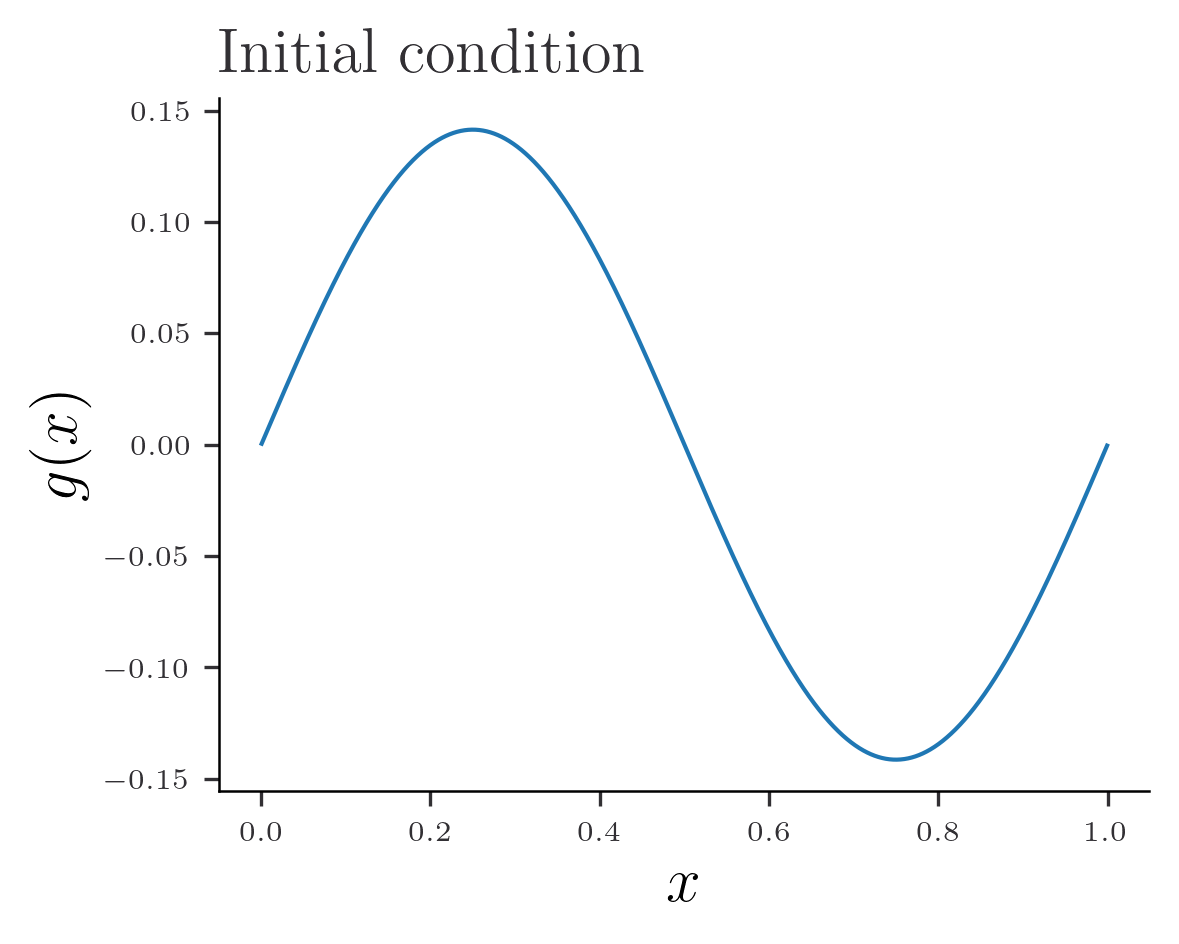} }}%
{{\includegraphics[width=0.31\columnwidth]{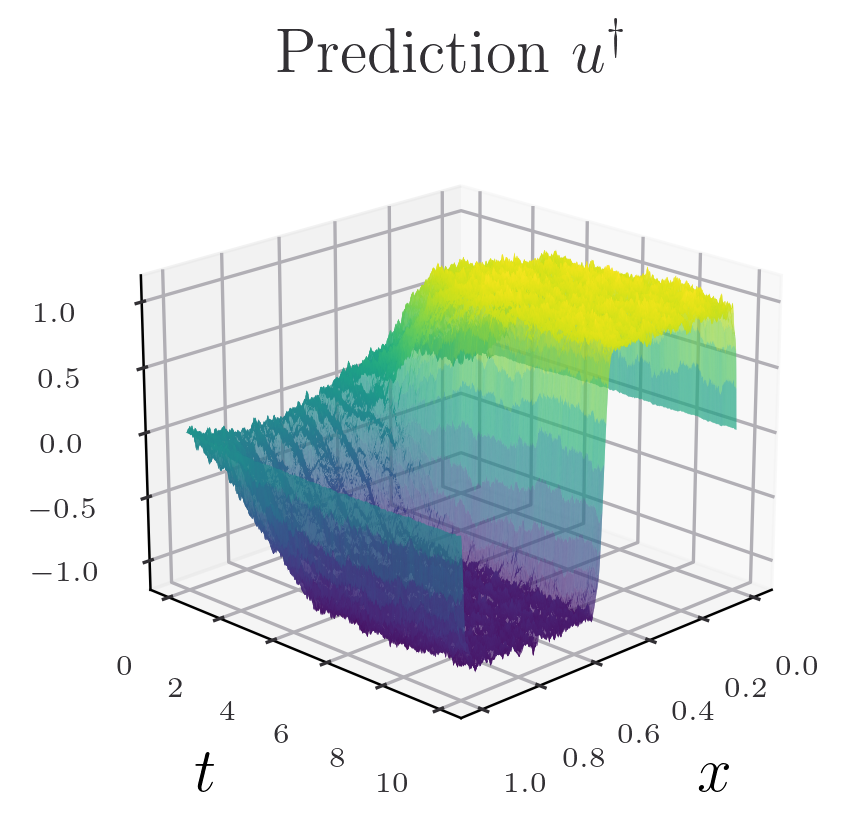} }}%
{{\includegraphics[width=0.31\columnwidth]{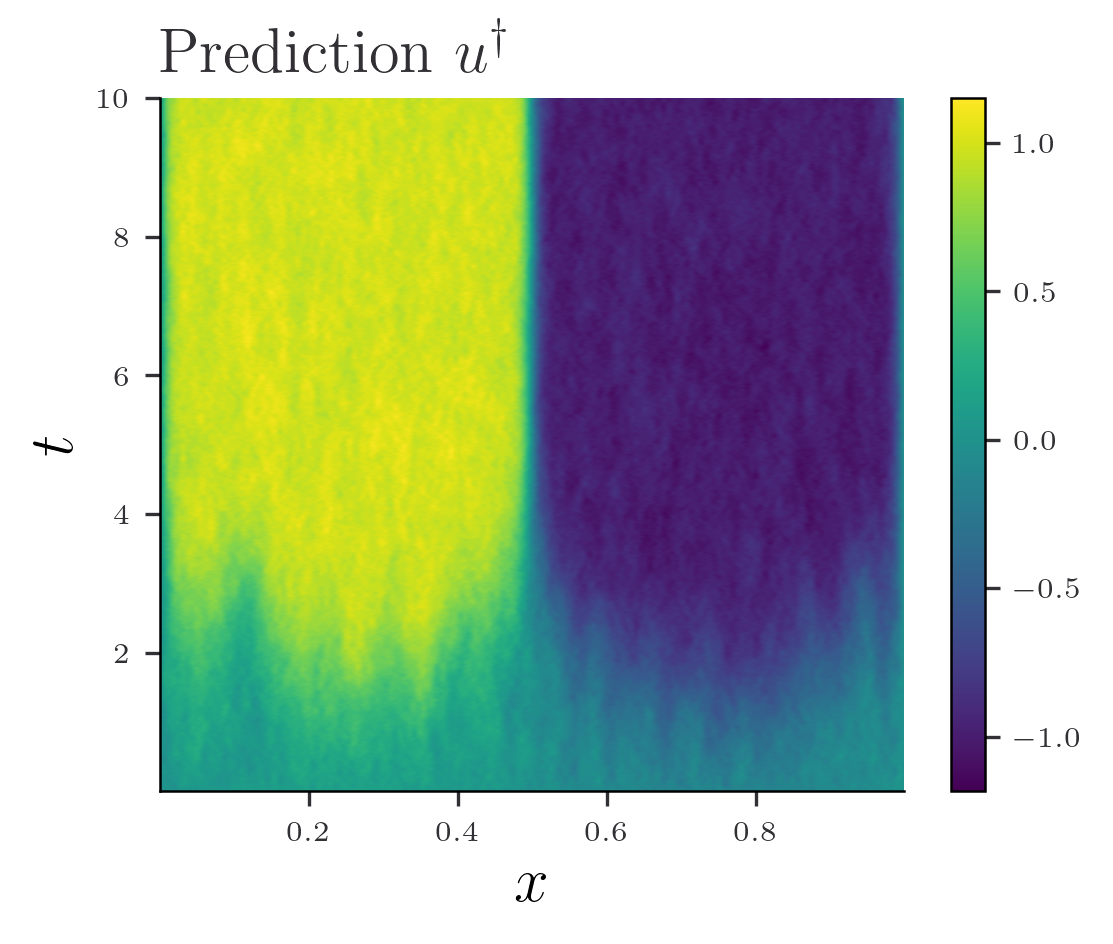} }}%
\caption{Kernel Solution computed with $n= 160$ and $\delta t = 2^{-10}$.}
\end{subfigure}

\begin{subfigure}[b]{\textwidth}    {{\includegraphics[width=0.31\columnwidth]{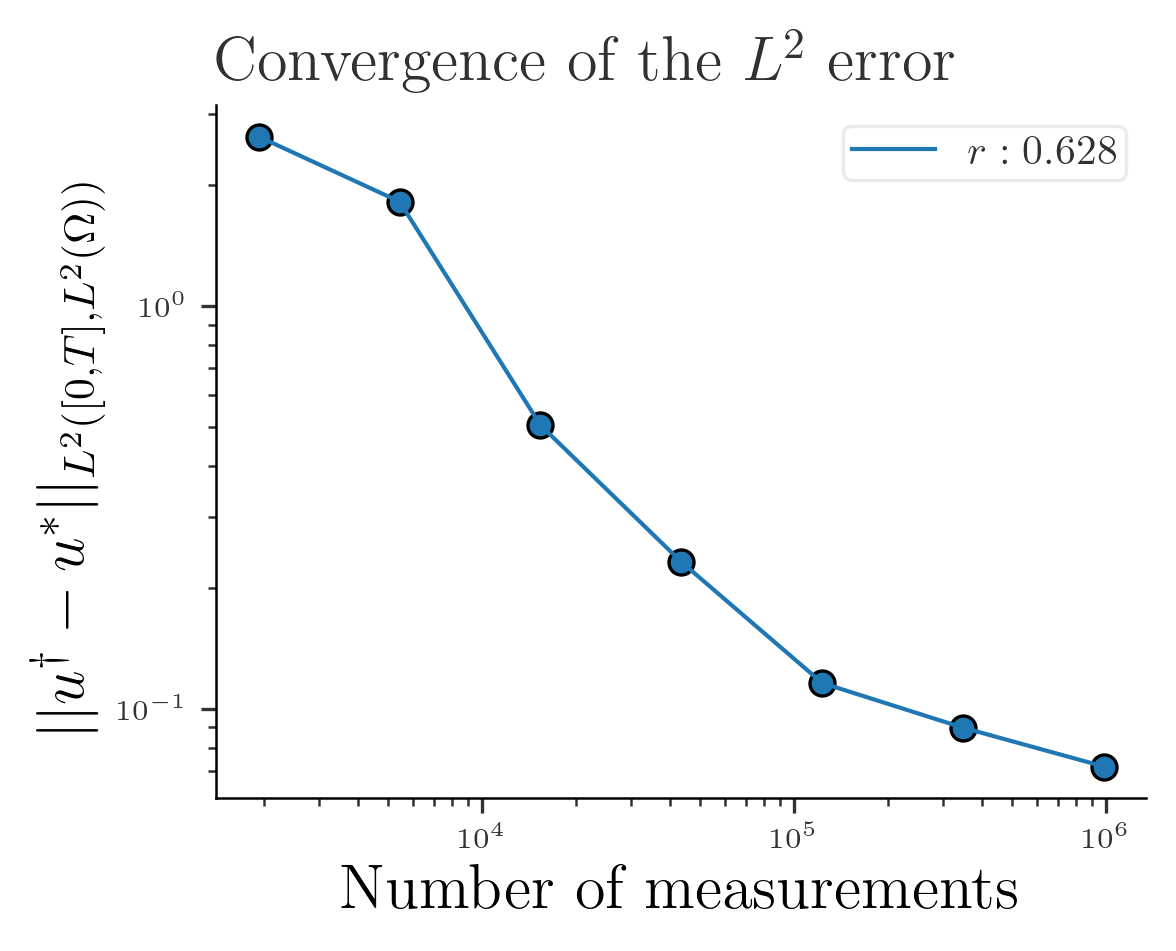} }}%
{{\includegraphics[width=0.31\columnwidth]{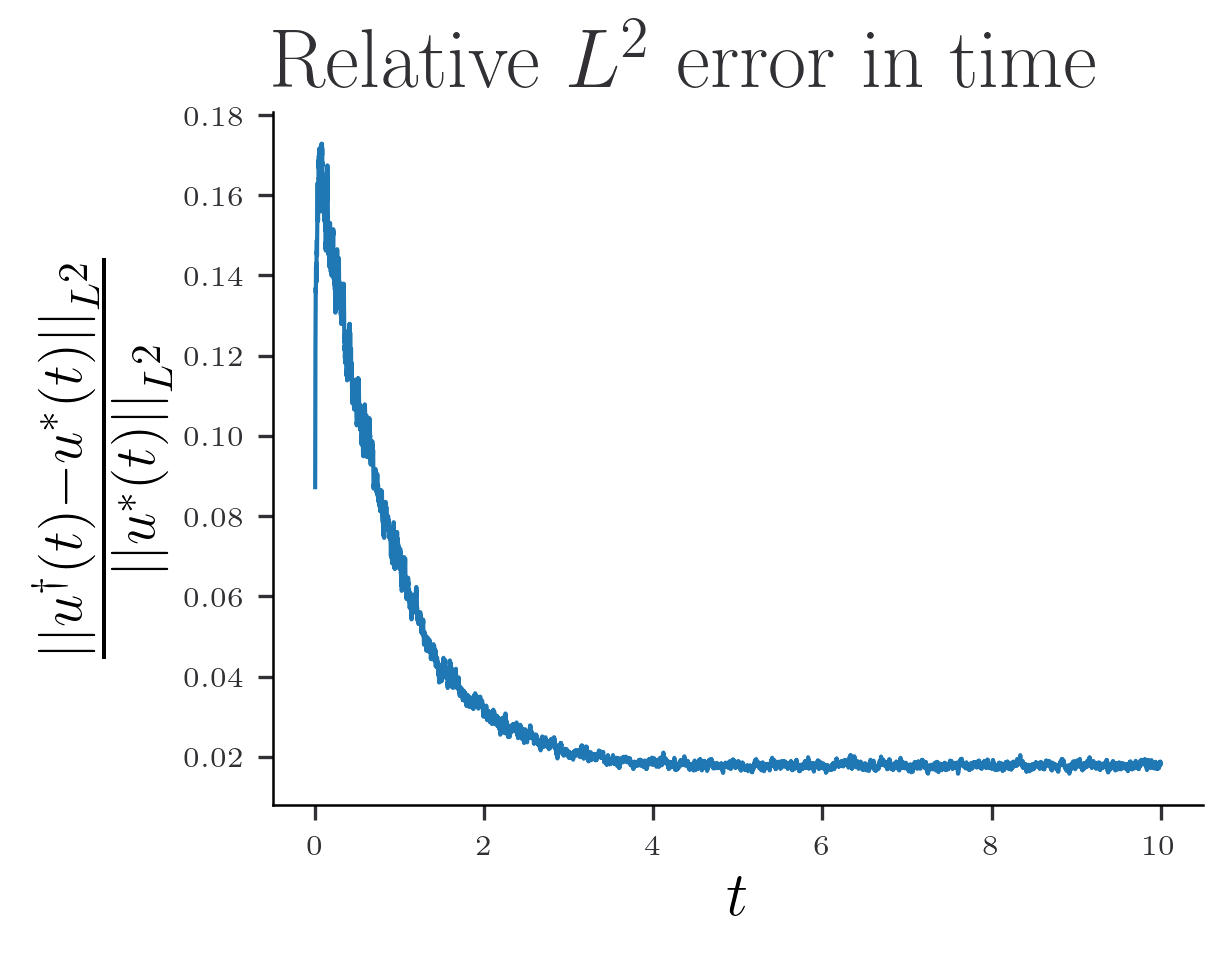} }}%
{{\includegraphics[width=0.31\columnwidth]{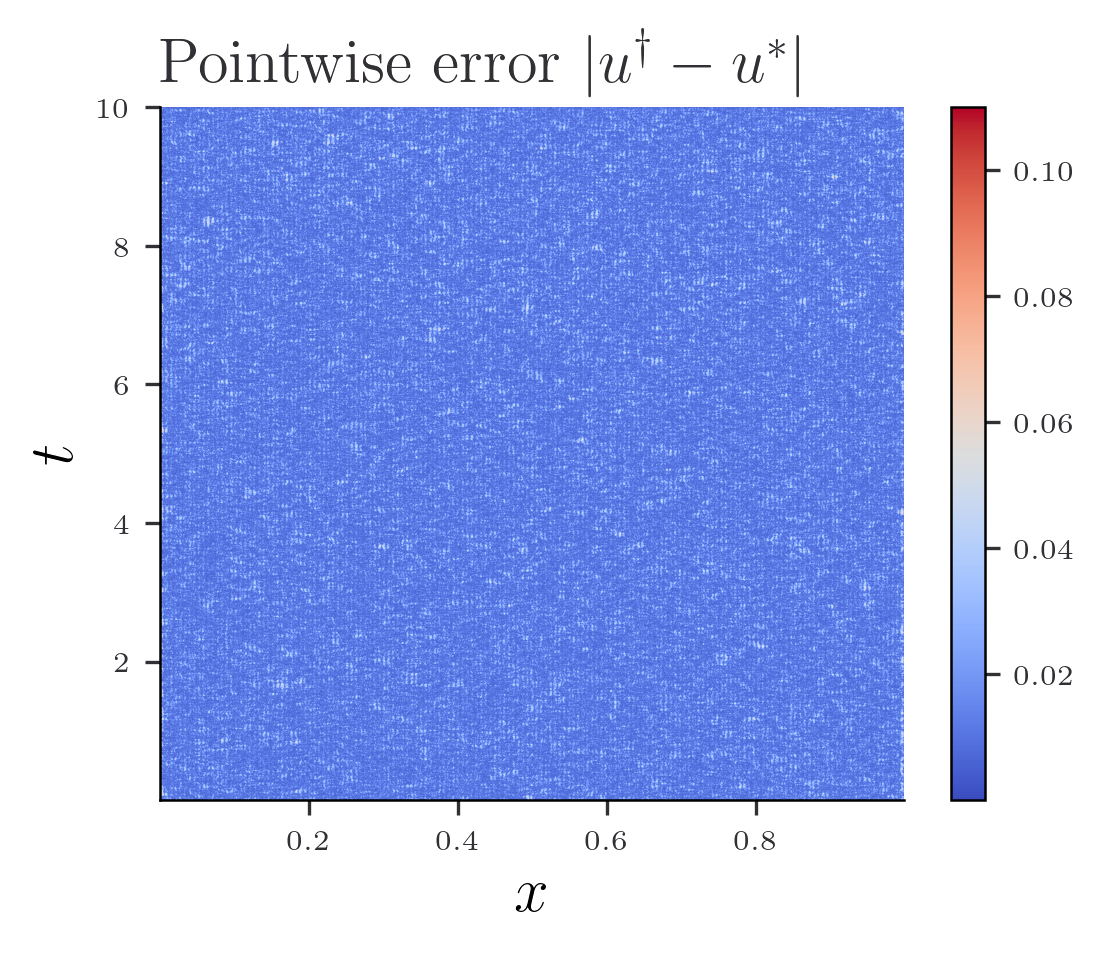} }}%
\caption{Error compared to the fine mesh solution.}
\end{subfigure}
    \caption{Stochastic Allen-Cahn equation.}
    \label{fig: allen-cahn}
\end{figure}

Similarly to the previous section, we compute the kernel-based prediction using piecewise linear finite element basis functions and assess the error relative to the reference solution \( u^* \).
To satisfy the CFL condition, we set \( N^2 \Delta t = 5 \), with the finest discretization parameters \( \Delta t = 2^{-10} \) and \( N = 160 \), which is coarser than for the fine mesh solution $u^*$. The results are illustrated in \cref{fig: allen-cahn}, and the quantitative results. for the finest discretization are presented in \cref{tab: time results}.

\section{Defining a New Solution Concept for Singular PDEs/SPDEs}
\label{sec:new_concept}
The proposed approach to solving PDEs with singular data can be extended to define a new solution concept for singular PDEs and SPDEs.  In this section we explain the idea, leaving detailed development for future work. Some PDEs/SPDEs, such as the stochastic Kardar–Parisi–Zhang (KPZ) equation \cite{kardar1986dynamic} and the $\Phi^4_3$ model \cite{glimm2012quantum}, do not satisfy our stability assumptions and are ill-posed in the usual sense. An emerging successful methodology for the interpretation of these equations, which formally contain an infinite divergent term, is the theory of regularity structures \cite{hairer2012solvingkpzequation, Hairer_2014}. In this section we provider an alternate approach, based on the generalization of \eqref{eq: nonlinear reg problem}, which
is amenable to the computational tools developed in this work. To be concrete, we describe
our approach in the context of nonlinear PDEs of the form:
    \begin{align}\label{eq: nonlinear pde new}
    \begin{split}
    \mathcal{P}(u, Du, D^2u) &= \xi,  \quad x \in\Omega, \\
            u &= 0,  \quad x \in \partial\Omega.
    \end{split}
    \end{align}
In the context of singular stochastic PDEs, the differential operator $\mathcal{P}$ involves nonlinearities that are not well-defined classically. For example, the Kardar–Parisi–Zhang (KPZ) equation includes the nonlinearity \( (Du)^2 \). When $\xi$ is space-time white noise, $Du$ has to be interpreted in a distributional sense, and the non-linear term $(Du)^2$ does not have an a priori valid interpretation.

To establish a well-defined notion of solutions to \eqref{eq: nonlinear pde new}, we introduce continuous regularizers \( (v_0, v_1, v_2) \in \mathcal{H}_{K_0} \times \mathcal{H}_{K_1} \times \mathcal{H}_{K_2} \) (a product of three RKHSs) in order to regularize the irregular terms \((u, Du, D^2u) \). This ensures that the regularization of the operator \( \mathcal{P} \) is well-defined in terms of \( (v_0, v_1, v_2) \). The solutions of \eqref{eq: nonlinear pde new} are then defined as the minimizers of the following optimization problem:
\begin{align}\label{eq: unconstrained min prob}
    \inf_{u\in \mathcal{H}_K, v_i \in \mathcal{H}_{K_i}} \Bigg( \norm{\mathcal{P}(v_0, v_1, v_2) - \xi}^2_{H^{-s}} + \sum_{i=0}^2 \norm{v_i - D^i u}^2_{\mathcal{H}_{\Gamma_i}} + \gamma \Big( \norm{u}^2_{\mathcal{H}_K} + \sum_{i=0}^2 \norm{v_i}^2_{\mathcal{H}_{K_i}} \Big) \Bigg). 
\end{align} 
Here, the spaces \( \mathcal{H}_{\Gamma_i}, \, i = 0, 1, 2 \), are regularization spaces chosen to account for the differences in regularity between the irregular terms \( (u, Du, D^2u) \) and their regularizers \( (v_0, v_1, v_2) \). For each choice of \(\gamma > 0\), a solution \(u^\gamma\) exists, as established using a similar approach to the previous sections; see Appendix~\ref{app: new concept of solution}, \cref{thm: existence of a minimizer new} for details. The question of interest becomes whether the solution $u^\gamma$ converges to a well-defined limit (and in which sense), i.e. whether an equivalent of \cref{thm:convergence to solution} exists for singular SPDEs. 
 

\section{Conclusion}

We have introduced a novel loss for solving partial differential equations with rough forcing terms. This loss is adapted to the irregularity of the right-hand side and allows us to extend the methodology used in \cite{chen2021solving}  beyond the setting of classical smooth solutions. We have shown that under mild assumptions on the PDE and the RKHS, the numerical solution provided by this kernel approach asymptotically converges to the true solution (in the $H^t$ norm). These theoretical results allow us to address the problem of kernel misspecification and show that even when the kernel is not adapted to the solution (i.e., if the solution does belong to the RKHS), we still can still guarantee convergence in the continuum limit of the discretization. Kernel methods can also be viewed from the Bayesian perspective of Gaussian Process Regression \cite{chen2024gaussian}. From this point of view, our approach addresses problems arising from employing misspecified priors, as defined by the choice of kernel. We have not, however, addressed the question of the convergence rates associated with a misspecified kernel and, empirically, the choice of the kernel does have an impact on performance. 

We have shown that the corresponding optimization problem for computing a numerical solution can be effectively solved by the Gaussian-Newton method in the functional space. Moreover, we have presented numerical experiments showing that both the kernel method and the neural network approach can effectively recover the solution of partial differential equations with rough forcing terms. Our proposed loss, based on a weak Sobolev norm, significantly improves the performance of both approximation methods as compared to minimizing the classical pointwise empirical $L^2$ loss. 

Additionally, we have explored the impact of loss choice on numerical solutions. We find that misspecified losses lead to high error, while overly weak losses result in overly smoothed solutions. Our findings underscore the importance of loss selection for enhancing the performance of machine learning solvers, indicating that alternatives to the $L^2$ loss can yield substantial improvements. The relevance of the choice of loss is also considered
in an operator learning context in \cite{li2021learning} and is integral to the Deep Ritz method
\cite{DeepRitz}. Our results demonstrate that both the use of test functions (rather than pointwise evaluations) and selecting the appropriate Sobolev norm are crucial for accurate solution recovery.
Finally, we have shown that the the proposed approach can also be employed to identify new solution concepts for singular PDEs and SPDEs.

\section*{Acknowledgments}

The authors acknowledge support from the Air Force Office of Scientific Research under MURI award number FA9550-20-1-0358 (Machine Learning and Physics-Based Modeling and Simulation).  HO acknowledges support by the Department of Energy under award number DE-SC0023163 (SEA-CROGS: Scalable, Efficient and Accelerated Causal Reasoning Operators, Graphs and Spikes for Earth and Embedded Systems). Additionally, AS and HO acknowledge support from the DoD Vannevar Bush Faculty Fellowship Program.  EC acknowledges support from the Resnick Sustainability Institute and the Vannevar Bush Faculty Fellowship held by AMS.

\bibliographystyle{siamplain}
\bibliography{references}


\appendix

\section{Auxiliary proofs}


\subsection{Convergence Rates for the Continuum Problem}
\label{app: convergence rates}

Consider the Poisson equation of the form \eqref{eq:Poisson}. For fixed regularizer $\gamma \in \R^+$, we consider the continuum problem as in \eqref{eq: nonlinear reg problem}, denoting its solution by $u_\gamma$. We recall that $u^*$ denotes the solution to \eqref{eq:Poisson}. Letting $\mathcal{H}_k \equiv H^q$ for some $q>-s+2$, 
in the following we show for which $\alpha \in \R$, the following two conditions hold:
\begin{equation*}
    \norm{u_\gamma}_{H^\alpha}<\infty \quad \text{and} \quad \lim_{\gamma \rightarrow 0 }\norm{u_\gamma-u_0}_{H^\alpha} = 0.
\end{equation*}
 To show this, we first find an expression for solution $u_\gamma$ in terms of the eigenvalues $\{\lambda_j\}_{j\in\mathbb{N}}$ of the negative Laplacian operator. We recall that any $u \in H^q_0$ can be expressed as  
\begin{equation*}
    u = \sum_{j = 1}^{\infty} u_j \varphi_j,\quad \text{with}\; u_j=\langle u,\varphi_j\rangle, 
\end{equation*}
where $\{\varphi_j\}_{j\in\mathbb{N}}$ is the set $L^2$ orthonormal generalized eigenvectors of the negative Laplacian operator, subject to Dirichlet boundary conditions. Similarly we recall that it is possible to express $\xi\in H^{-s}$ from \eqref{eq:Poisson} as
\begin{equation*}
    \xi = \sum_{j = 1}^{\infty} \xi_j \varphi_j,\quad \text{with}\; \xi_j=\langle \xi,\varphi_j\rangle_{L^2}.
\end{equation*}
Given these expressions, we have by \cref{eq: energy norm}
\begin{equation*}
\norm{r}^2_{H^{-s}} = \sum_{j=1}^\infty (\xi_j-\lambda_ju_j)^2\lambda_j^{-s}, \quad \norm{u}^2_{H^q_0} = \sum_{j=1}^\infty u_j^2\lambda_j^q,
 \end{equation*}
 from which we can reformulate the optimization problem as
 \begin{equation}
 \label{eq:reg_problem_basis}
     (u_\gamma)_j = \argmin_{v\in\R}(\xi_j-\lambda_jv)^2\lambda_j^{-s}+\gamma v^2\lambda_j^q \quad \text{for each}\;j\in\mathbb{N}.
 \end{equation}
It is readily found that the quadratic optimization problem \eqref{eq:reg_problem_basis} yields
\begin{equation*}
u_\gamma = \sum_{j=1}^\infty\frac{\xi_j}{\lambda_j\big(1+\gamma\lambda_j^{q+s-2} \big)}\varphi_j.
\end{equation*}
\begin{proposition}
\label{prop:conv_linear}
For any $\alpha \leq 2-s$, it holds that 
\begin{equation}
\label{eq:convstatement}
    \norm{u_\gamma}_{H^\alpha}<\infty \quad \text{and} \quad \lim_{\gamma \rightarrow 0 }\norm{u_\gamma-u_0}_{H^\alpha} = 0.
\end{equation}
\end{proposition}
\begin{proof}
    We divide the proof into two parts. Indeed we first consider for which $\alpha \in \R$ it holds that $\norm{u_\gamma}_{H^\alpha}<\infty$. Using the growth rate of the eigenvalues of the negative Laplacian operator $\lambda_j\sim j^{2/d}$, we have that
    \begin{subequations}
    \label{eq:alphanorm}
    \begin{align}
        \norm{u_\gamma}^2_{H^\alpha} &= \sum_{j=1}^\infty\left(\frac{\xi_j}{\lambda_j\big(1+\gamma\lambda_j^{q+s-2} \big)}\right)^2\lambda_j^\alpha\\
        &\lesssim\sum_{j=1}^\infty \xi_j^2\lambda_j^{-s}\frac{j^{2(s+\alpha -2)/d}}{\big(1+\gamma j^{2(q+s-2)/d} \big)^2}\\
        &\lesssim \gamma^{-2}\sum_{j=1}^\infty \xi_j^2\lambda_j^{-s}\cdot j^{2(\alpha -2q-s+2)/d},
    \end{align}
    \end{subequations}
    since clearly,
    \[
    \big(1+\gamma j^{2(\alpha -2q-2s+2)/d}\big)^{-2}\leq   \big(\gamma j^{2(\alpha -2q-2s+2)/d}\big)^{-2}.
    \]
    Therefore, with the knowledge that $\norm{\xi}^2_{H^{-s}}<\infty$, (\ref{eq:alphanorm}c) gives that $\norm{u_\gamma}_{H^\alpha}<\infty$ if $\alpha \leq 2q+s-2$.

    We now turn our attention to $\norm{u_\gamma-u_0}_{H^\alpha}$. Again, using the growth rate of the eigenvalues of the negative Laplacian operator $\lambda_j\sim j^{2/d}$, we have that
    \begin{subequations}
    \begin{align*}
       \norm{u_\gamma-u_0}_{H^\alpha}^2 &= \sum_{j=1}^\infty\left(\frac{\xi_j}{\lambda_j\big(1+\gamma\lambda_j^{q+s-2} \big)}-\frac{\xi_j}{\lambda_j}\right)^2\lambda_j^\alpha\\
       &\lesssim \sum_{j=1}^\infty \frac{\gamma^2\xi_j^2j^{2(2q+2s-6+\alpha)/d}}{\big( 1+\gamma j^{2(q+s-2)/d}\big)^2}\\
       &\lesssim \sum_{j=1}^\infty \gamma^2\xi_j^2\lambda_j^{-s}\frac{j^{2(2q+3s-6+\alpha)/d}}{\big( 1+\gamma j^{2(q+s-2)/d}\big)^2}\\
       &\lesssim\sum_{j=1}^\infty \gamma^2\xi_j^2\lambda_j^{-s}\frac{j^\kappa}{\big( 1+\gamma j^{\tau}\big)^2},
    \end{align*}
    \end{subequations}
where we have defined
\begin{equation*}
\kappa:= \frac2d(2q+3s-6+\alpha)\quad\text{and}\quad \tau:= \frac2d(q+s-2).
\end{equation*}
For ease of exposition we define 
\begin{equation*}
\label{eq:Kjgamma}
B_j^\gamma := \frac{\gamma^2 j^\kappa}{\big( 1+\gamma j^{\tau}\big)^2},
\end{equation*}
so that 
\begin{subequations}
\begin{align*}
    \norm{u_\gamma-u_0}_{H^\alpha}^2 &= \sum_{j=1}^\infty \xi_j^2\lambda_j^{-s}B_j^\gamma\\
    &= \sum_{j=1}^M \xi_j^2\lambda_j^{-s}B_j^\gamma + \sum_{j>M}^\infty \xi_j^2\lambda_j^{-s}B_j^\gamma,
\end{align*}
\end{subequations}
for any $M\in \mathbb{N}$. Clearly $B_j^\gamma \rightarrow 0$ as $\gamma\rightarrow0$, therefore 
\[
\sum_{j=1}^M \xi_j^2\lambda_j^{-s}B_j^\gamma \rightarrow 0 \quad \text{as}\;\;\gamma\rightarrow0.
\]
On the other hand, noting that
\[
\big(1+\gamma j^\tau \big)^{-2}\leq \big(\gamma j^\tau \big)^{-2},
\]
we have that $B_j^\gamma \leq j^{\kappa -2\tau}$, and hence
\begin{equation*}
    \sum_{j>M}^\infty \xi_j^2\lambda_j^{-s}B_j^\gamma \leq \sum_{j>M}^\infty \xi_j^2\lambda_j^{-s}j^{\kappa -2\tau}.
\end{equation*} 
By the fact that $\norm{\xi}^2_{H^{-s}}<\infty$, it holds that $ \sum_{j>M}^\infty \xi_j^2\lambda_j^{-s}B_j^\gamma<\infty$ if $\kappa -2\tau \leq 0$. Indeed, for small enough $\gamma$, the sequence of partial sums converges uniformly, thus allowing the interchangeability of limits. Furthermore, if this holds, then 
\begin{equation*}
    \sum_{j>M}^\infty \xi_j^2\lambda_j^{-s}B_j^\gamma \rightarrow0 \quad \text{as} \quad \gamma \rightarrow0,
\end{equation*}
since $B_j^\gamma \rightarrow 0$ as $\gamma\rightarrow0$. Since $\kappa-2\tau \leq 0$ if and only if $\alpha \leq 2-s$, we conclude that for $\alpha \leq 2-s$ it holds that 
\begin{equation*}
\lim_{\gamma \rightarrow 0 }\norm{u_\gamma-u_0}_{H^\alpha} = 0.
\end{equation*}
To finish the proof, we thus note that for $\alpha \leq \min\{2q+s-2,2-s\}$, \eqref{eq:convstatement} holds, but since $q>2-s$, we have that $2q+s-2>2-s$, hence the conclusion. 
\end{proof}

\begin{remark}
    The form of $B_j^\gamma$ in the proof of~\cref{prop:conv_linear}, leads to the conjecture that $u_\gamma$ converges to $u_0$ as $\gamma\to 0$ in $H^{2-s}$ as fast as $\gamma^2$, which is confirmed by the numerical experiments.
\end{remark}

\subsection{Numerical method}
\label{app: numerical method}
Proof of \cref{prop: expression n norm}
\begin{proof}
    We begin by noting the following equivalent characterizations:
\begin{align*}
    |f|_{\Phi^N}&=  \sup_{u \in \Phi^N}\frac{\int f u}{\norm{u}_{H^s_0}} \\
    &=  \sup_{u \in \Phi^N, \norm{u}_{H^s_0} = 1} \int fu\,.
\end{align*}
Moreover, observe that 
\begin{equation*}
     \int fu = \sum_{i=1}^N c_i\int f \varphi_i = c^\intercal [f, \varphi],
\end{equation*}
for some $c\in\R^N$, and for $u\in V^N$
\begin{equation*}
    \norm{u}_{H^s_0} =  \sqrt{\sum_{i=1}^N\sum_{j=1}^N c_i c_j \int \varphi_i (-\Delta)^s \varphi_j} = \sqrt{c^\intercal A c}.
\end{equation*}
Therefore we can write 
\begin{align*}
    |f|_{\Phi^N}= & \sup_{c \in \R^N}   c^\intercal [f, \varphi] \\
    &\text{s.t.} \quad  c^\intercal A c = 1\,,
\end{align*}
whose solution is:
\begin{equation*}
     |f|_{\Phi^N}=  \sqrt{[f, \varphi]^\intercal A^{-1}[f, \varphi]}.
\end{equation*}
Note that the matrix $A$ is invertible because the $\varphi_i$ are linearly independent.
\end{proof}
Proof of \cref{prop: norm consistency}.
\begin{proof}
    Under \cref{assumption: density of test space} for all $\varepsilon >0$, there is a $N = N(\varepsilon, u)$ such that 
    \begin{equation*}
        \norm{v^N - u}_{H^s_0(\Omega)} < \varepsilon,
    \end{equation*}
    for some $v^N \in \Phi^N$. Now, by the definition of the dual norm, for all $\varepsilon >0$, there is a $u \in H^s_0(\Omega)$ with $\norm{u}_{H^s_0({\Omega})}  = 1$ such that
    \begin{equation*}
        \Big| \int f u - \norm{f}_{H^s_0(\Omega)} \Big| < \frac{\varepsilon}{2}.
    \end{equation*}
    By the density of $\Phi^N$, there is a $N$ and a $v^N$ such that 
    \begin{align*}
        \norm{u - v^N}_{H^{s}_0(\Omega)} \leq \frac{\varepsilon}{2\norm{f}_{H^{-s}(\Omega)}}.
    \end{align*}
    Without loss of generality, we may assume that $v^N$ has a unit norm since $u$ has a unit norm. Therefore
        \begin{align*}
        \Big|\int fu - \int f v^N \Big| &\leq \norm{f}_{H^{-s}(\Omega)} \norm{u - v^N}_{H^{s}_0(\Omega)} \\
        &\leq  \frac{\varepsilon}{2}.
    \end{align*}
Hence for all $\varepsilon > 0$, there is a $N = N(\varepsilon)$ and a corresponding $v^N \in \Phi^N$ with norm 1 such that 
\begin{align*}
      \Big| \int f v^N - \norm{f}_{H^s_0(\Omega)} \Big| &\leq \Big| \int f u - \norm{f}_{H^s_0(\Omega)} \Big| + \Big|\int fu - \int f v^N \Big| \\
      & < \varepsilon.
\end{align*}
The assertion then follows.
\end{proof}

Proof of \cref{thm: representer}
\begin{proof}
     Let $S:=\textrm{span}\{K(\cdot,\chi_1),\dots,K(\cdot,\chi_N),K(\cdot,x_1),\dots,K(\cdot,x_M)\}$. Clearly, it holds that $\mathcal{H}_K = S \bigoplus S^\perp$. Therefore, we have that 
    \begin{equation*}
        \widehat{u} = \sum_{i=1}^N\alpha_iK(\cdot,\chi_i) + \sum_{j=1}^M\beta_jK(\cdot,x_j) + s,
    \end{equation*}
    where $s\in S^\perp$. By the reproducing kernel property, for any $1\leq k \leq N$ we have that
    \begin{subequations}
    \begin{align*}
        \widehat{u}(\chi_k) &= \langle \sum_{i=1}^N\alpha_iK(\cdot,\chi_i) + \sum_{j=1}^M\beta_jK(\cdot,x_j) , K(\cdot,\chi_k)\rangle + \langle s, K(\cdot,\chi_k)\rangle  \\
        &= \langle \sum_{i=1}^N\alpha_iK(\cdot,\chi_i) + \sum_{j=1}^M\beta_jK(\cdot,x_j) , K(\cdot,\chi_k)\rangle,
    \end{align*}
    \end{subequations}
    since $s\in S^\perp$. Similarly, for any $1\leq k\leq M$, we have that 
    \begin{subequations}
    \begin{align*}
        \widehat{u}(x_k) &= \langle \sum_{i=1}^N\alpha_iK(\cdot,\chi_i) + \sum_{j=1}^M\beta_jK(\cdot,x_j) , K(\cdot,x_k)\rangle + \langle s, K(\cdot,x_k)\rangle  \\
        &= \langle \sum_{i=1}^N\alpha_iK(\cdot,\chi_i) + \sum_{j=1}^M\beta_jK(\cdot,x_j) , K(\cdot,x_k)\rangle.
    \end{align*}
    \end{subequations}
    Now, since $s$ is orthogonal to $\sum_{i=1}^N\alpha_iK(\cdot,\chi_i) + \sum_{j=1}^M\beta_jK(\cdot,x_j)$,
    we have that 
    \begin{equation*}
        \|\widehat{u}\|^2_{\mathcal{H}_K} = \Bigl\|\sum_{i=1}^N\alpha_iK(\cdot,\chi_i) + \sum_{j=1}^M\beta_jK(\cdot,x_j) \Bigr\|^2_{\mathcal{H}_K} + \|s\|^2_{\mathcal{H}_K}.
    \end{equation*}
    Hence, while $s$ does not affect the first term in (\ref{eq:infinite_dim_problem thm}a), an $s$ of smaller $\mathcal{H}_K$ norm decreases the second term. Therefore, any solution $\widehat{u}$ must have $s=0$.
\end{proof}

Proof of Proposition \ref{prop: representer, explicit solution}.
\begin{proposition}
     Using \cref{thm: representer} we reformulate the problem in \cref{eq:infinite_dim_problem} as a quadratic minimization problem with linear inequality constraint. Each iteration $n$ necessitates solution of 
    \begin{align*}
         \arginf_{c \in \R^{N+M}} &\,\frac{1}{2}\bigg[\sum_{i=1}^{M+N}c_i\mathcal{P}'(u)K(\cdot,\varphi)- r_n  , \varphi\bigg] ^\intercal A^{-1} \bigg[\sum_{i=1}^{M+N}c_i\mathcal{P}'(u)K(\cdot,\varphi_i)- r_n , \varphi\bigg]
         + \frac{\gamma}{2} \norm{\sum_{i=1}^{M+N}c_iK(\cdot,\varphi_i)}^2_{\mathcal{H}_K}\\
   &\text{s.t. } \sum_{i=1}^{M+N}c_iK(x_j,\varphi_i) = g(x_j) \quad\text{for }j=1,\dots,M.
    \end{align*}
    We may transform the above into a matrix expression:
    \begin{align*}
         \arginf_{c \in \R^{N+M}} &\frac{1}{2}\Big( K(\chi, \varphi)c - [r_n, \varphi]\Big)^\intercal A^{-1} \Big(K(\chi, \varphi)c - [r_n, \varphi] \Big)
         + \frac{\gamma}{2} c^\intercal K(\varphi, \varphi)c\\
   &\text{s.t. }K(X, \varphi)c = g(X) \quad\text{for }j=1,\dots,M.
    \end{align*}
    which may be formulated a standard convex quadratic minimization problem with equality constraints:
    \begin{align*}
         \argmin_{c \in \R^{N+M}} & \frac{1}{2}c^\intercal \Big(K(\varphi, \chi)A^{-1} K(\chi, \varphi) + \gamma\ K(\varphi, \varphi) \Big) c  
         - \Big( K(\varphi, \chi)A^{-1}[g_n, \varphi] \Big)^\intercal c \\
         &+ \frac{1}{2}[g_n, \varphi]^\intercal A^{-1} [g_n, \varphi]
         \\
   &\text{s.t. }K(X, \varphi)c = g(X) \quad\text{for }j=1,\dots,M.
    \end{align*}
The KKT system gives the explicit solution to this problem \cite[p. 522]{boyd2004convex}:
\begin{align*}
\begin{bmatrix}
 K(\chi, \varphi)^\intercal A^{-1} K(\chi, \varphi) + \gamma K(\varphi, \varphi) & K(X, \varphi)  \\
K(X, \varphi) &  0
\end{bmatrix}
\begin{bmatrix}
    c \\
    \nu
\end{bmatrix}
= 
\begin{bmatrix}
    K(\chi, \varphi)^\intercal A^{-1}[r_n, \varphi]  \\
    g(X),
\end{bmatrix}
\end{align*}
with $\nu \in \R^{N}$ are the Lagrange multipliers. Note that this system has a unique solution when the KKT matrix is non-singular, which is the case when the $\varphi_i$ are linearly independent and the collocation points $x_j$ are distinct. 
\end{proposition}


\subsection{New Concept of Solution}
\label{app: new concept of solution}
We now present a proof of the existence of a minimizer to the generalized problem \eqref{eq: nonlinear pde new} under suitable assumptions.

\begin{assumption}\label{assumption: compact emb and cont}
   We assume the following compact embeddings: 
   \begin{subequations}
       \begin{align*}
           \mathcal{H}_K \subset  H^{s+i} &\hookrightarrow \mathcal{H}_{\Gamma_i}\quad \text{for } i = 0,1,2. \\
          \mathcal{H}_{K_i}  &\hookrightarrow  \mathcal{H}_{\Gamma_i} \quad \text{for } i = 0,1,2.
       \end{align*}
   \end{subequations}
   Moreover, we assume that the nonlinear differential operator $\mathcal{P}:\mathcal{H}_{\Gamma_0}\times \mathcal{H}_{\Gamma_1}\times\mathcal{H}_{\Gamma_2} \to\mathcal{H}_{\Gamma_3}$    
   is continuous. 
\end{assumption}

\begin{assumption}[Continuity of $\mathcal{P}$]
    We assume that the nonlinear differential operator $\mathcal{P}:\mathcal{H}_{\Gamma_0}\times \mathcal{H}_{\Gamma_1}\times\mathcal{H}_{\Gamma_2} \to\mathcal{H}_{\Gamma_3}$   
    is continuous. 
\end{assumption}
The next theorem is the equivalent of \cref{thm: existence of a minimizer} for the generalized problem.

\begin{theorem}\label{thm: existence of a minimizer new}
   Let $\mathcal{H} \coloneqq \mathcal{H}_K \times \mathcal{H}_{K_0}\times \mathcal{H}_{K_1} \times \mathcal{H}_{K_2}$ for some appropriate choices of RKHS $\mathcal{H}_K,\mathcal{H}_{K_i}$ satisfying Assumption \ref{assumption: compact emb and cont}. Let $m$ be the optimal value of problem \eqref{eq: unconstrained min prob}:
   \begin{equation*}
       m \coloneqq   \inf_{u, v} L(u, v; \gamma).
   \end{equation*}
    Then $0 \leq  m  < \infty$ and there is a $(u^*, v^*) \in \mathcal{H}$ such that $m = L(u^*, v^*)$. Moreover, if $(u^n, v^n)$ is a minimizing sequence, then there is a subsequence which converges strongly to $(u^*, v^*)$ in $\mathcal{H}$ 
    \begin{subequations}, i.e.
        \begin{align*}
        u^n \rightarrow u^*  \quad &\text{in } \mathcal{H}_K\\
            v^n_i \to v^{*}_i \quad &\text{in } \mathcal{H}_{K_i}\quad \text{for } i = 0,1,2.
        \end{align*}
    \end{subequations}
\end{theorem}

\begin{proof}
    To prove existence of a minimizer, we first show the boundedness of the minimizing sequence and then proceed by showing weak lower semi-continuity of the objective, thus proving the claim. 
    
    \textit{Step 1: existence of weak limit}.
    Indeed, we first observe that
    \begin{equation*}
        0 \leq L(u, v; \gamma) <\infty \quad \text{for all } (u, v) \in \mathcal{H},
    \end{equation*}
    and therefore $0 \leq m < \infty$. For ease of presentation, we will now omit the dependence on $\gamma$. Consider some minimizing sequence $(u^n, v^n)$:
    \begin{equation*}
        L(u^n, v^n) \to m \quad \text{as } n \to \infty. 
    \end{equation*}
    Then observe that by the definition of $L$, this implies that the minimizing sequence is bounded
    \begin{subequations}
        \begin{alignat*}{2}
            &\norm{u^n}_{\mathcal{H}_K} &&\leq C \\
            &\norm{v^n_i}_{K_i} &&\leq C \quad \text{for } i= 0,1,2.
        \end{alignat*}
    \end{subequations}
    By the Banach–Alaoglu theorem \cite{alaoglu1940weak}, since each $\mathcal{H}_{K_i}$ and $\mathcal{H}_{K}$ are reflexive, there is a subsequence (which we denote with the same index) such that $(u^{n}, v^{n}) \rightharpoonup (u^{\infty}, v^{\infty})$ in $\mathcal{H}$. 

    \textit{Step 2: weak lower semicontinuity of the objective.} By the compact embedding assumption \ref{assumption: compact emb and cont}, we deduce that
    \begin{subequations}
        \begin{align*}
            D^iu^n \to D^iu^{\infty} \quad \text{in } \mathcal{H}_{\Gamma_i} \quad \text{for } i = 0,1,2 \\
            v^n_i \to v^{\infty}_i \quad \text{in } \mathcal{H}_{\Gamma_i} \quad \text{for } i = 0,1,2.
        \end{align*}
    \end{subequations}
    Therefore
    \begin{equation*}
        \norm{v^n_i - D^i u^n}^2_{\Gamma_i} \to \norm{v^\infty_i - D^iu^\infty}^2_{\mathcal{H}_{\Gamma_i}} \quad \text{for } i = 0,1,2.
    \end{equation*}
    Moreover by the continuity of $\mathcal{P}$ from $\mathcal{H}_{\Gamma_0} \times \mathcal{H}_{\Gamma_1} \times \mathcal{H}_{\Gamma_2}$ into $\mathcal{H}_{\Gamma_3}$, we obtain
    \begin{equation*}
        \norm{\mathcal{P}(v^n_0, v^n_1, v^n_2) - f}^2_{\mathcal{H}_{\Gamma_3}} \to \norm{\mathcal{P}(v^\infty_0, v^\infty_1, v^\infty_2) - f}^2_{\mathcal{H}_{\Gamma_3}}.
    \end{equation*}
    Hence the first part of the term of $L$ is weakly continuous on $\mathcal{H}$. Since norms are lower-semicontinuous \cite[p. 61]{brezis_2010}, it follows that $L: \mathcal{H} \to \R$ is lower-semicontinuous. Therefore
    \begin{equation*}
        (u^{n}, v^{n}) \rightharpoonup (u^{\infty}, v^{\infty}) \Rightarrow \liminf_{n \to \infty} L(u^{n}, v^{n}) \geq L(u^{\infty}, v^{\infty}),
    \end{equation*}
    and consequently
    \begin{equation*}
        m = L(u^{\infty}, v^{\infty}, \mathbf{r}^{\infty}).
    \end{equation*}
    
    \textit{Step 3: strong convergence of the minimizing sequence.} We now turn our attention to showing the existence of a subsequence of $(u^n,v^n)$ that converges strongly to a minimizer. For ease of presentation, let 
    \begin{subequations}
        \begin{align*}
        \Phi(u, v) &\coloneqq \norm{\mathcal{P}(v_0, v_1, v_2) - f}^2_{\mathcal{H}_{\Gamma_3}} + \sum_{i=0}^2 \norm{v_i - D^i u}^2_{\mathcal{H}_{\Gamma_i}} \\
        \norm{(u, v)}^2 &\coloneqq \norm{u}^2_{\mathcal{H}_k} +  \sum_{i=0}^2 \norm{v_i}^2_{\mathcal{H}_{K_i}},
    \end{align*}
    \end{subequations}
    so that 
    \begin{equation*}
        L(u, v) =  \Phi(u, v) + \gamma \norm{(u, v)}^2.
    \end{equation*}
    Using the polarization equality
    \begin{align*}
        \frac{\gamma}{4} \norm{(u, v)^l - (u, v)^k}^2 &=   \gamma \Big( \frac{1}{2}\norm{(u, v)^l}^2 + \frac{1}{2}\norm{(u, v)^k}^2 - \norm{\frac{1}{2}\big((u,v)^l + (u, v)^k \big)}^2\Big) \\
        &= \frac{1}{2}L((u, v)^k) + \frac{1}{2}L((u, v)^l) - L\bigg(\frac{1}{2}\big((u,v)^l + (u, v)^k \big)\bigg) \\ & \quad - \frac{1}{2}\Phi((u, v)^l) - \frac{1} {2}\Phi((u, v)^k) + \Phi\bigg(\frac{1}{2}\big((u,v)^l + (u, v)^k \big)\bigg). 
    \end{align*}
    Letting $l,k \geq N_1(\delta)$ so that $L(u^k, v^k)\leq m+\delta$ and $L(u^k, v^k)\leq m+\delta$, we obtain
    \begin{subequations}
        \begin{align*}
             \frac{\gamma}{4} \norm{(u, v)^l - (u, v)^k}^2 &\leq m + \delta - m- \frac{1}{2}\Phi((u, v)^l) - \frac{1} {2}\Phi((u, v)^k) + \Phi\bigg(\frac{1}{2}\big((u,v)^l + (u, v)^k \big)\bigg) \\
             & = \delta - \frac{1}{2}\Phi((u, v)^l) - \frac{1} {2}\Phi((u, v)^k) + \Phi\bigg(\frac{1}{2}\big((u,v)^l + (u, v)^k \big)\bigg).
        \end{align*}
    \end{subequations}
    Since $(u,v)^l$, $(u,v)^k$ and $\frac{1}{2}\big((u,v)^l + (u, v)^k \big)$ converge to $(u, v)^\infty$ in  each $\mathcal{H}_{\Gamma_i}$, and since $\Phi$ is continuous on $\mathcal{H}_{\Gamma_0} \times \dots \times \mathcal{H}_{\Gamma_3}$, we deduce that there is a $N_2(\delta)$ such that for $k,l \geq N_2(\delta)$:
    \begin{align*}
        \frac{\gamma}{4} \norm{(u, v)^l - (u, v)^k}^2 &\leq \delta +\delta = 2\delta.
    \end{align*}
    Therefore the sequence $(u, v)^n$ is Cauchy in $\mathcal{H}_K \times\mathcal{H}_{K_0}\times \mathcal{H}_{K_1}\times \mathcal{H}_{K_3}$ and as a consequence converges strongly in those spaces: 
    \begin{subequations}
        \begin{align*}
            u^n \rightarrow u^\infty  \quad &\text{in } \mathcal{H}_K \\
            v^n_i \to v^{\infty}_i \quad &\text{in } \mathcal{H}_{K_i} \quad \text{for } i = 0,1,2.
        \end{align*}
    \end{subequations}
\end{proof}

\section{Neural Network Architecture}\label{sec: NN architecture}
\begin{table}[H]
    \centering
    \begin{tabular}{|l|c|c|c|c|c|}
        \hline
        Problem & RFF $m$ & RFF Std $\sigma$ & N$^{\circ}$ Layers & Neurons/layer & Total Parameters \\
        \hline
        1D Elliptic Linear & 24 & 5.0 & 4 & 16 & 1,073 \\
        \hline
        2D Semilinear Elliptic & 32 & 5.0 & 4 & 64 & 8,385 \\
        \hline
    \end{tabular}
    \caption{Architecture for the neural networks used for the NeS-PINN approach.}
    \label{tab: nn architecture}
\end{table}
\cref{tab: nn architecture} describes the architecture of our deep neural networks with $N$ layers. 
\paragraph{Random Fourier Features (RFF)}

As described in \cite{rff_pinn}, the random Fourier features add a non-trainable embedding at the start of the network:
\begin{align}
    \gamma(x) = \begin{bmatrix}
        \cos(Bx) \\
        \sin(Bx)
    \end{bmatrix}
\end{align}
where $B\in \R^{m\times d}$ has entries sampled $i.i.d.$ from a normal distribution $\mathcal{N}(0, \sigma^2)$. This embedding mitigates the spectral bias of neural networks. 
\section{Gaussian Measures and White Noise}\label{app: gaussian measures}
We now briefly introduce the concepts of Gaussian measures, Gaussian spatial noise, and Gaussian space-time white noise. Our exposition will be informal and far from complete. For a more detailed explanation, see \cite{lord2014introduction}. 
\subsection{Gaussian Measures}
We introduce Gaussian measures on $L^2(\Omega)$.
\begin{definition}[$L^2$ valued Gaussian measure]
    An $L^2(\Omega)$-valued random variable $\xi$ is Gaussian if $\inp{X}{\phi}$ is a Gaussian random variable for all $\phi \in L^2(\Omega)$.
\end{definition}
To each Gaussian random variable on $L^2(\Omega)$ one can associate a covariance operator. 
\begin{definition}
    A linear operator $\mathcal{Q}: L^2(\Omega) \rightarrow L^2(\Omega)$ is the covariance of an $L^2$-valued random variable $\xi$ if 
    \begin{align*}
        \inp{\mathcal{Q}\phi}{\psi}_{L^2} = \text{Cov}\Big( \inp{\xi}{\phi}_{L^2}\inp{\xi}{\psi}_{L^2}\Big).
    \end{align*}
\end{definition}
The above definitions imply that for any $\phi \in L^2(\Omega)$
\begin{align*}
    \inp{\xi}{\phi}_{L^2} \sim \mathcal{N}(\inp{\mu}{\phi}_{L^2},  \inp{\mathcal{Q}\phi}{\phi}_{L^2}
    ).
\end{align*}
We will always consider centered Gaussian variables $\E[\xi] = 0$ so that $\inp{\mu}{\phi} = 0$ for all $\phi \in L^2(\Omega)$. One can further show that the covariance operator $\mathcal{Q}$ is of trace class.
\begin{proposition}[Corollary 4.41 \cite{lord2014introduction}]
    Let $\xi$ be a $L^2(\Omega)$-valued Gaussian random variable, then its covariance operator $\mathcal{Q}$ is of trace-class.
 \end{proposition}
This allows for the construction of $\xi$ via the spectral decomposition of $\mathcal{Q}$. 

\begin{proposition}
    Let $\xi$ be a centered $L^2(\Omega)$- valued Gaussian random variable with covariance operator $\mathcal{Q}$. Then $\xi$ admits the representation
    \begin{align*}
        \xi = \sum_{i=1}^\infty \sqrt{\lambda_i} Z_i e_i \quad Z_i \sim \mathcal{N}(0,1) \; \text{i.i.d.},
    \end{align*}
    where $\{\lambda_{i}, e_i\}_{i=1}^\infty$ are the eigenvalues and eigenvectors of $\mathcal{Q}$.  
\end{proposition}
\subsection{Spatial White Noise}
We now briefly introduce the concept of space-time with noise. Formally, this extends the above definitions to the operator $\mathcal{Q} = I$, which is not a trace-class operator. First, we consider the formal power series
    \begin{align}
    \label{eq: wn power series app}
        \xi = \sum_{i=1}^\infty  Z_i e_i \quad Z_i \sim \mathcal{N}(0,1) \; \text{i.i.d.},
\end{align}
where $\{ e_i\}_{i=1}^\infty$ is any orthonormal basis of $L^2(\Omega)$. For each $\varphi \in L^2(\Omega)$, we may formally compute 
\begin{align*}\
    \inp{\xi}{\varphi}_{L^2} = \sum_{j=1}^\infty Z_j \inp{e_j}{\varphi}_{L^2}. 
\end{align*}
Ignoring issues of convergence, we can view the above as normal random variable with mean and covariance
\begin{align*}
    &\E[\inp{\xi}{\varphi}_{L^2} ] =  \sum_{j=1}^\infty \E[Z_j] \inp{e_j}{\varphi}_{L^2} = 0 \\\
     &\E[\inp{\xi}{\varphi}_{L^2} \inp{\xi}{\psi}_{L^2}] =  \sum_{j=1}^\infty \sum_{i=1}^\infty \E[Z_j Z_i] \inp{e_j}{\varphi}_{L^2} \inp{e_j}{\psi}_{L^2} =  \inp{\varphi}{\psi}_{L^2}.
\end{align*}
From this, we may instead  define $\xi$ to be a stochastic process acting on elements of $L^2(\Omega)$ such that  
\begin{align*}
   [\xi, \varphi] \sim \mathcal{N}(0, \norm{\varphi}^2_{L^2}).
\end{align*} 
The series \eqref{eq: wn power series app} is not convergent in $L^2(\Omega)$. This can be fixed by considering the inclusion
\begin{align*}
    i: L^2(\Omega) \rightarrow H^{-1}(\Omega),
\end{align*}
which is Hilbert-Schmidt\footnote{One can make choices other than $H^{-1}(\Omega)$ and any space $V$ such that the inclusion $i: L^2(\Omega) \rightarrow V$ is Hilbert-Schmidt is appropriate.}. The series \ref{eq: wn power series app} does converge in this space and the associated process is not a full Gaussian measure but instead a cylindrical measure on $L^2(\Omega)$. We refer to this process as spatial white process which formally acts on $\varphi \in L^2(\Omega)$ by 
\begin{align*}
   [\xi, \varphi] \sim \mathcal{N}(0, \norm{\varphi}^2_{L^2}).
\end{align*}
\subsection{Space-time White Noise: the Cylindrical Wiener Process}
We now define the cylindrical Wiener process. Let $(\mathcal{S}, \mathcal{F}, \mathcal{F}_t, \mathbb{P})$ be a filtered probability space. 

\begin{definition}
    The cylindrical Wiener-process on $L^2(\Omega)$ if 
    \begin{align}\label{eq: wiener process expansion app}
        \xi(t) = \sum_{j=1}^\infty \beta_j(t) e_j,
    \end{align}
    for any orthonormal basis $\{ e_j\}_{j=1}^\infty$ and $\beta_j(t)$ are i.i.d $\mathcal{F}_t$-Brownian motions. 
\end{definition}
One can show that 
\begin{equation}\label{eq: increment property app}
     \xi(t) -  \xi( s) \sim \mathcal{N}(0, (t-s)I) \quad \text{for all } 0 \leq s \leq t,
\end{equation}
so that the increments of the cylindrical Wiener process are spatial white-noise processes on $L^2(\Omega)$ (with a variance defined by time-step of the increments). Moreover, since the (formal) time-derivative of Brownian motion is white noise in time, this justifies the choice of $\dot{\xi}(t)$ as space-time white noise. Equivalently, one can see $\xi(t)$ as space-time white noise integrated in time.

\subsection{Sampling From a Cylindrical Wiener Process}
In order to solve SPDEs of the form \cref{eq: time-dependent PDE}, we consider integration in time through an implicit or semi-implicit Euler-scheme (see \cref{sec: implicit euler}). To this end, we need to sample increments of the cylindrical Wiener process $\xi(t) -  \xi(s)$. The first option is to directly use the representation \cref{eq: wiener process expansion app} to obtain a truncated expansion
\begin{align}\label{eq: truncated expansion app}
     \Delta \xi^N := \sum_{j=1}^N \Delta\beta_j e_j,
\end{align}
where $\Delta \beta_j \sim \Delta W_j = \mathcal{N}(0, \Delta t)$ are i.i.d. increments of Brownian motion. In the case where $L^2([0,1]) \hookrightarrow H^{-1}([0,1])$ we would select $e_j(x) = \sin(2\pi j x)$. Other boundary conditions (such as Neumann and periodic boundary conditions) can also be considered by selecting the full Fourier basis $e^{i2\pi jx}$.

Another option is to use \cref{eq: increment property app} to show that for any collection $\Phi = \{\varphi_i\}_{i=1}^N \subset L^2(\Omega)$, the vector 
\begin{align} \label{eq: cylindrical_two}
    [\Delta \xi,\varphi] \sim \mathcal{N}(0, \Delta t A),
\end{align}
where $A_{ij} = \inp{\varphi_i}{\varphi_j}_{L^2}$ is the stiffness matrix. If the $\varphi_i$ are selected to be approximations to the identity around a point $x$, then one can informally view $[\Delta \xi,\varphi]$ as approximating the value of $ \Delta \xi$ at the point $x$ (although we stress that action of  $ \Delta \xi$ against $\delta_{x}$ are not well defined as it does not have the required regularity). Likewise, while the truncated series \cref{eq: truncated expansion app} has well-defined pointwise values, the limiting series \cref{eq: wiener process expansion app} does not.

\begin{figure}[htp]%
\begin{subfigure}[t]{1.0\textwidth}
    \centering
    \subfloat{{\includegraphics[width=0.3\columnwidth]{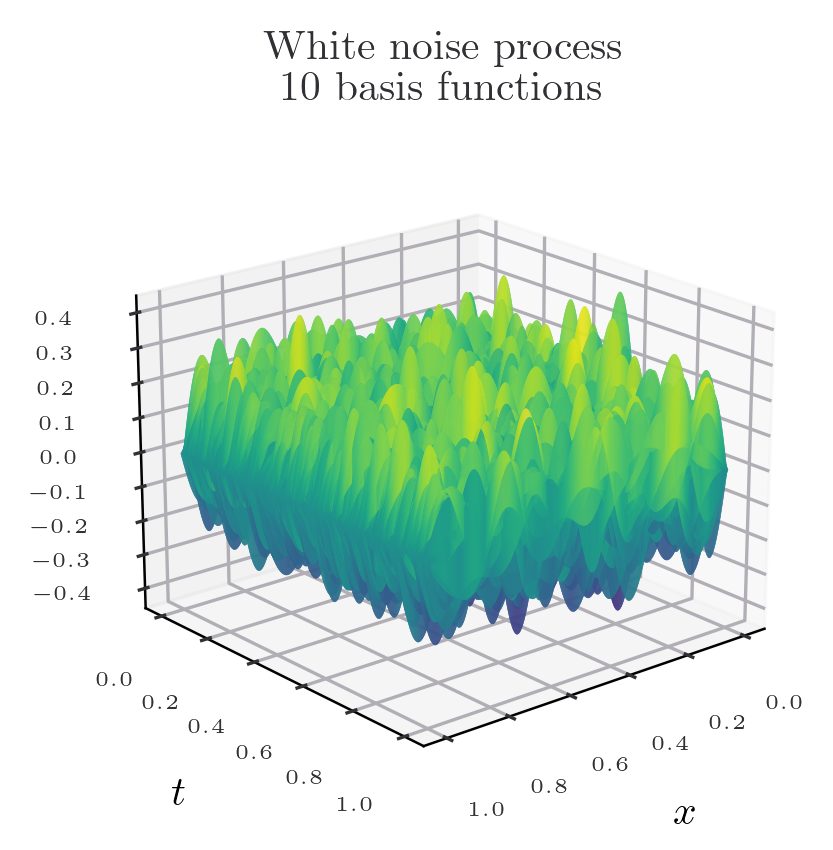} }}%
    \subfloat{{\includegraphics[width=0.3\columnwidth]{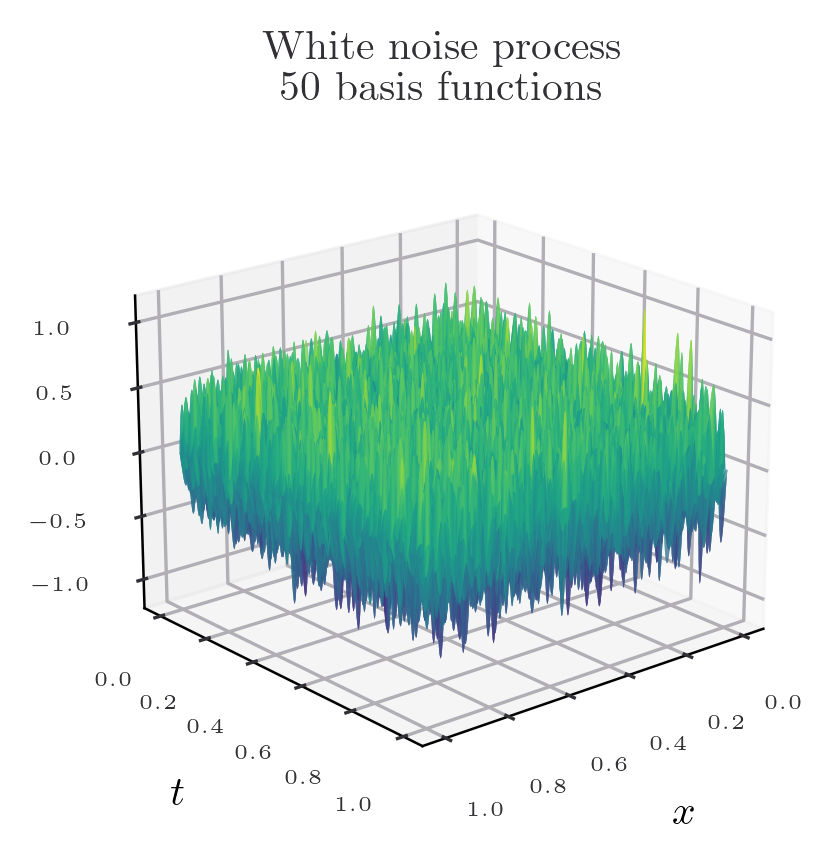} }}%
    \subfloat{{\includegraphics[width=0.3\columnwidth]{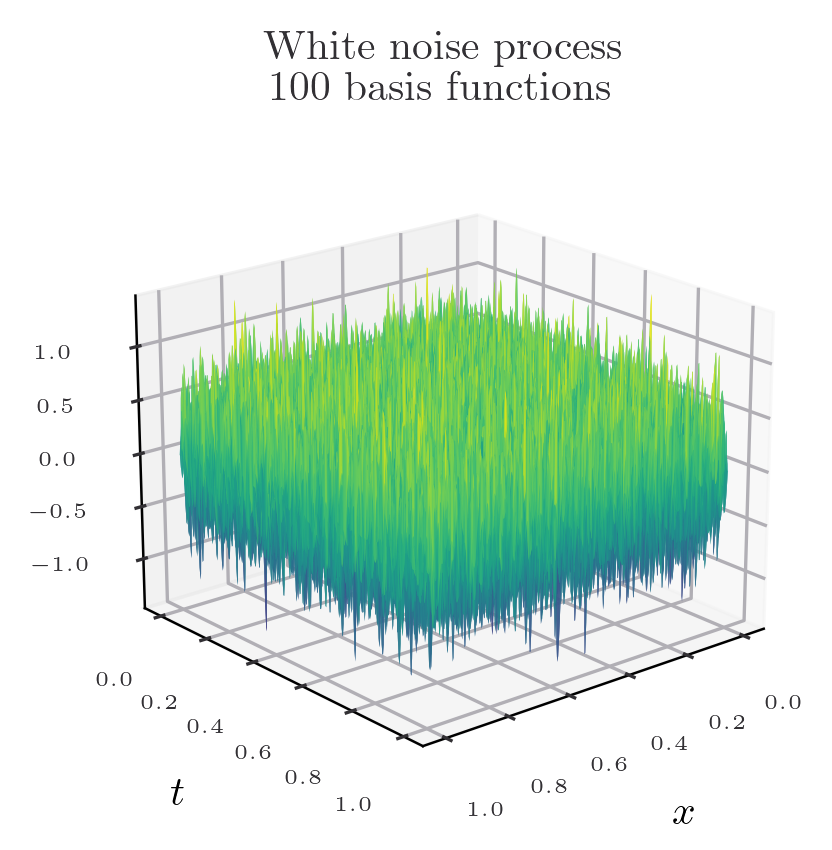} }}%
    \caption{Spectral basis}%
\end{subfigure}
\begin{subfigure}{1.0\textwidth}
    \centering
    \subfloat{{\includegraphics[width=0.3\columnwidth]{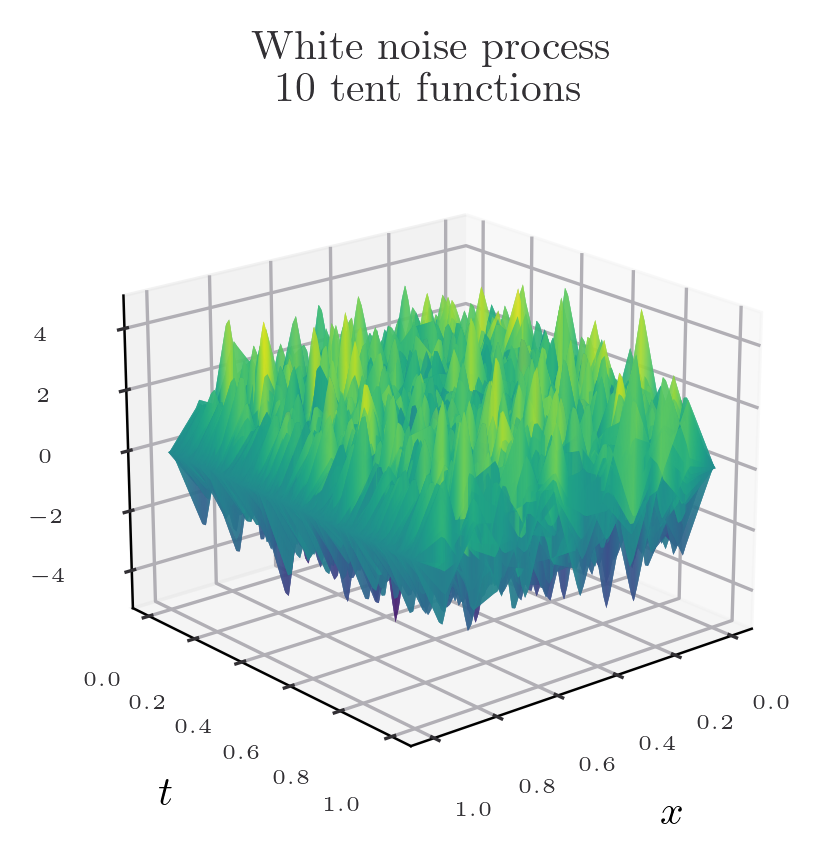} }}%
    \subfloat{{\includegraphics[width=0.3\columnwidth]{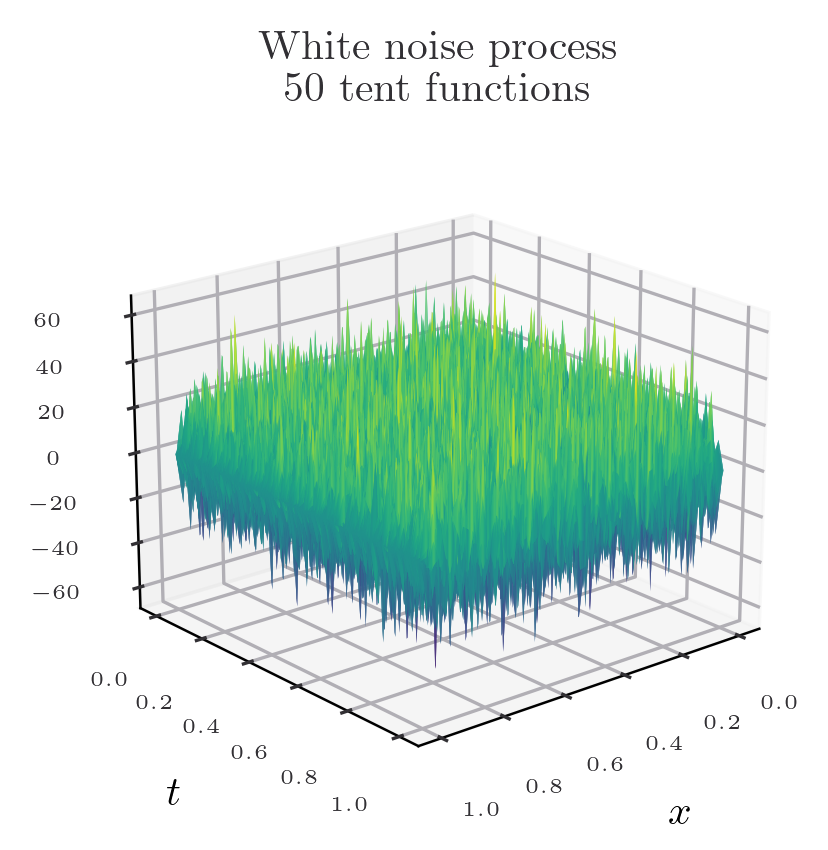} }}%
    \subfloat{{\includegraphics[width=0.3\columnwidth]{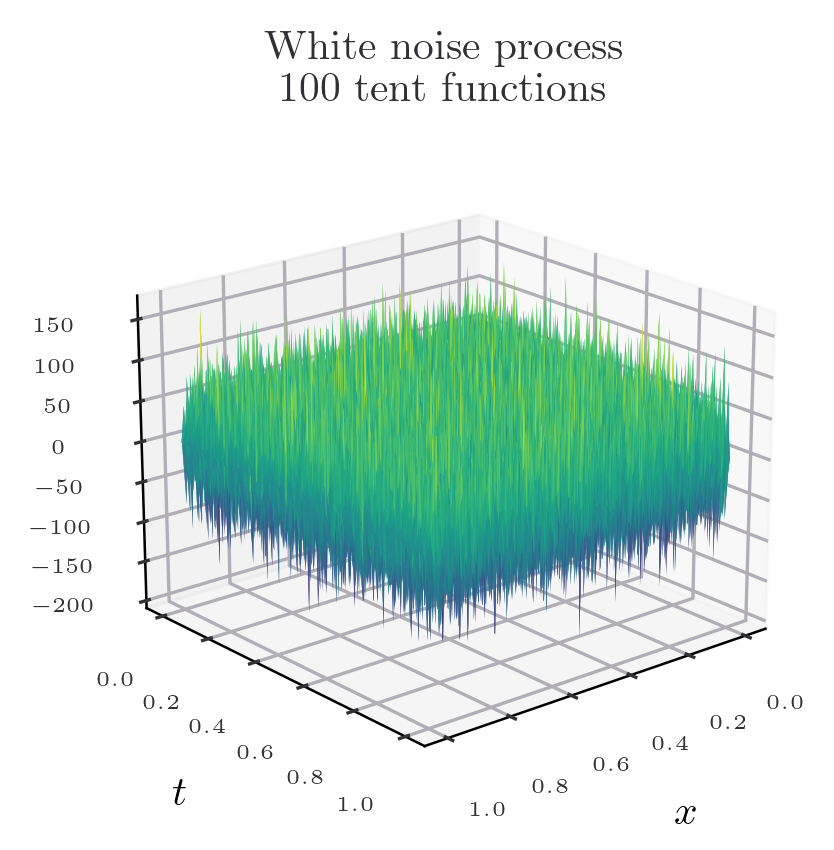} }}%
    \caption{Tent function basis}%
    \end{subfigure}
    \caption{Space-time white noise generated by two methods.}
    \label{fig: space-time white noise}
\end{figure}

An illustration of the two sampling methods based on expansion in~\eqref{eq: truncated expansion app} and the projection in~\eqref{eq: cylindrical_two} is shown in \cref{fig: space-time white noise}. 

\end{document}